\documentclass[preprint,3p]{elsarticle}
\journal{J. Comput. Phys.}

\usepackage{stmaryrd}
\usepackage{amsmath,amssymb}
\usepackage{tabularx,epsfig,color,tabularx}
\usepackage{epstopdf}
\usepackage{algorithm,algorithmic,multirow}
\usepackage{graphicx}
\usepackage{caption}
\usepackage{subcaption}
\usepackage{mathtools}
\usepackage{epsfig}

\newcommand{\be}{\begin{equation}}
\newcommand{\ee}{\end{equation}}
\newcommand{\ba}{\begin{array}}
\newcommand{\ea}{\end{array}}
\newcommand{\bea}{\begin{eqnarray}}
\newcommand{\eea}{\end{eqnarray}}
\newcommand{\beas}{\begin{eqnarray*}}
\newcommand{\eeas}{\end{eqnarray*}}

\newcommand{\bx}{{\bf x}}

\newtheorem{lemma}{Lemma}[section]

\newcommand{\nn}{\nonumber}
\newcommand{\ap}{\alpha}

\newcommand{\lmd}{\lambda}

\begin{document}

\begin{frontmatter}

\title{A Jacobi spectral method for computing eigenvalue gaps and their distribution statistics of the fractional Schr\"{o}dinger operator}
\author[nus]{Weizhu Bao}
 \ead{matbaowz@nus.edu.sg}

\address[nus]{Department of Mathematics, National University of
Singapore, Singapore 119076, Singapore}
\ead[url]{http://blog.nus.edu.sg/matbwz/}

\author[csrc,nus]{Lizhen Chen\corref{5}}
\ead{lzchen@csrc.ac.cn}
\address[csrc]{Beijing Computational Science Research Center,
Beijing, 100193, P.R. China}
\address[sd]{School of Mathematics, Shandong University, Jinan 250100, P.R. China}

\cortext[5]{Corresponding author.}
\author[sd]{Xiaoyun Jiang}
\ead{wqjxyf@sdu.edu.cn}
\author[csrc]{Ying Ma}
\ead{yingma@csrc.ac.cn}

%%%%% Begin Abstract %%%%%%%%%%%
\begin{abstract}
We propose a spectral method by using the Jacobi functions for computing
eigenvalue gaps and their distribution statistics of the fractional Schr\"{o}dinger operator (FSO). In the problem, in order to get reliable gaps distribution statistics, we have to calculate accurately and efficiently a very large number of eigenvalues, e.g. up to thousands or even millions eigenvalues,
of an eigenvalue problem related to the FSO. For simplicity, we start
with the eigenvalue problem of FSO in one dimension (1D), reformulate
it into a variational formulation and then discretize it by using the Jacobi spectral method. Our numerical results
demonstrate that the proposed Jacobi spectral method has several advantages
over the existing finite difference method (FDM) and finite element method (FEM) for the problem: (i) the Jacobi spectral method is spectral accurate,
while the FDM and FEM are only first order accurate; and more importantly (ii) under a fixed number of degree of freedoms $M$, the Jacobi spectral method
can calculate accurately a large number of eigenvalues with the number
proportional to $M$, while the FDM and FEM perform badly when a large number of eigenvalues need to be calculated. Thus the proposed Jacobi spectral method is extremely suitable and demanded for the discretization of an eigenvalue problem when a large number of eigenvalues need to be calculated. Then the Jacobi spectral method is applied to study numerically
the asymptotics of the nearest neighbour gaps, average gaps, minimum gaps,
normalized gaps and their distribution statistics in 1D. Based on our numerical results, several interesting numerical observations (or conjectures)
about eigenvalue gaps and their distribution statistics of the FSO in 1D are formulated. Finally, the Jacobi spectral method is extended to the directional fractional Schr\"{o}dinger operator in high dimensions and extensive numerical results about eigenvalue gaps and their distribution statistics are reported.
\end{abstract}

\begin{keyword}
fractional Schr\"{o}dinger operator, Jacobi spectral method, nearest neighbour gaps, average gaps, minimum gaps, normalized gaps, gaps distribution statistics.
\end{keyword}

\end{frontmatter}

\section{Introduction}\setcounter{equation}{0}
Consider the eigenvalue problem of the
fractional Schr\"{o}dinger operator (FSO) (or time-independent
fractional Schr\"{o}dinger equation) in one dimension (1D):

Find $\lambda \in \mathbb{R}$ and a nonzero real-valued function $u(x)\ne 0$ such that
\begin{equation}\label{fproblem}
\begin{split}
L_{\rm FSO}\; u(x)&:=\left[(-\partial_{xx})^{\alpha/2}+V(x)\right]u(x)
=\lambda\; u(x), \qquad x\in \Omega:=(a,b), \\
u(x)&=0, \qquad  x \in \Omega^c:=\mathbb{R} \backslash \Omega,
\end{split}
\end{equation}
where $0<\alpha\le 2$, $V(x)\in L^2(\Omega)$ is a given real-valued function and the fractional Laplacian operator (FLO) $(-\partial_{xx})^{\alpha/2}$ is defined via
the Fourier transform (see \cite{ZRK07,DG13,LPGSGZMCMA18} and references therein) as
  \begin{equation}\label{foperator}
(-\partial_{xx})^{\alpha/2}\;u(x)={\mathcal F}^{-1} (|\xi|^{\alpha}(\mathcal{F}u)(\xi))\qquad x,\xi\in \mathbb{R},
\end{equation}
with ${\mathcal F}$ and ${\mathcal F}^{-1}$ the Fourier transform
and inverse Fourier transform \cite{CS07,LPGSGZMCMA18,G15}, respectively. We
remark here that an alternative way to define $(-\partial_{xx})^{\alpha/2}$
is through the principle value integral (see \cite{SKM93,V09,DGLZ12,L13,DG13} and references therein)
as
\be\label{pvi}
(-\partial_{xx})^{\alpha/2}\;u(x):=C_1^\alpha \int_{\mathbb R} \frac{u(x)-u(y)}{|x-y|^{1+\alpha}}dy, \qquad x\in {\mathbb R},
\ee
where $C_1^\alpha$ is a constant whose value can be computed explicitly as
\[C_1^\alpha =\frac{2^\alpha \Gamma((1+\alpha)/2)}{\pi^{1/2}|\Gamma(-\alpha/2)|}=\frac{\alpha \Gamma((1+\alpha)/2)}{2^{1-\alpha}\pi^{1/2}\Gamma(1-\alpha/2)}.\]
Another remark here is that the problem \eqref{fproblem} is equivalent
to the problem defined on the whole $x$-axis  ${\mathbb R}$ by taking the potential
$V(x)=+\infty$ for $x\in \Omega^c$.
When $\alpha=2$, \eqref{fproblem} collapses to the (classical)
time-independent Schr\"{o}dinger equation (or a standard Sturm-Liouville
eigenvalue problem) which has been widely used for determining
energy levels and their corresponding stationary states of a quantum particle within an external potential
$V(x)$ in quantum physics and chemistry \cite{CT95} and many other areas \cite{LB02,CCJ07,CS05}. When $\alpha=1$, the FLO $(-\Delta)^{1/2}$
and its variation $(\beta-\Delta)^{1/2}$ with $\beta>0$ a constant
have been widely adopted in representing Coulomb interaction and dipole-dipole interaction in two dimensions \cite{BaoB,BaoJ,Cai,Jiang} and modeling relativistic quantum mechanics for boson star \cite{Elg,BaoD}.
When $0<\alpha<2$, \eqref{fproblem} is usually referred as the time-independent fractional Schr\"{o}dinger equation (or  fractional eigenvalue problem)
which has been widely adopted for computing energy levels and their stationary states
in fractional quantum mechanics \cite{L13,BaoJ,Cai}, polariton condensation and quantum fluids of lights \cite{Car,Pins15}, while the FSO can be interpreted via the Feynman path integral approach over Brownian-like quantum
paths or over the L\'{e}vy-like quantum paths, see \cite{V09,L13,JLPR13} and references
therein.

Without loss of generality, we assume that $V(x)$ is non-negative, i.e.
$V(x)\ge0$ for $x\in\Omega$.
Since all eigenvalues of  \eqref{fproblem} are distinct (or all spectrum are discrete and no continuous spectrum), we can rank (or order) all eigenvalues of \eqref{fproblem} as
\be\label{eigenorder}
0<\lambda_1^\alpha<\lambda_2^\alpha\le \ldots\le \lambda_n^\alpha\le\ldots\,,
\ee
where the times that an eigenvalue $\lambda$ of \eqref{fproblem} appears in the above sequence
\eqref{eigenorder} is the same as its algebraic multiplicity. When $V(x)\equiv 0$ for $x\in\Omega$, all eigenvalues of \eqref{fproblem} are simple eigenvalues,
i.e. their algebraic multiplicities are all equal to $1$, then all $\le$ in
\eqref{eigenorder} can be replaced by $<$.
Define the {\sl nearest neighbor gaps} as \cite{JMRR99}
\be\label{nngap}
\delta_{\rm nn}^\alpha(N):= \lambda_{N+1}^\alpha-\lambda_{N}^\alpha, \qquad
N=1,2,3,\ldots,
\ee
where when $N=1$, i.e., $\delta_{\rm nn}^\alpha(1)= \lambda_{2}^\alpha-\lambda_{1}^\alpha:=\delta_{\rm fg}(\alpha)$ (i.e. the difference between the first two smallest eigenvalues) is called as the {\sl fundamental gap} of the FSO \eqref{fproblem}, which has been studied analytically and/or numerically for $\alpha=2$ \cite{AB89,AC10,BR18} and
$0<\alpha\le 2$ \cite{BRSS18,BBROA18}; the {\sl minimum gaps} as \cite{BBR16,R18}
\be\label{mingap}
\delta_{\rm min}^\alpha(N):=\min_{1 \leq n \leq N}\ \delta_{\rm nn}^\alpha(n)
=\min_{1 \leq n \leq N}\ \lambda_{n+1}^\alpha-\lambda_{n}^\alpha,
\qquad N=1,2,3,\ldots;
\ee
the {\sl average gaps} as \cite{JMRR99}
\be\label{avegap}
\delta_{\rm ave}^\alpha(N):=\frac{1}{N}\sum_{n=1}^N\delta_{\rm nn}^\alpha(n)=\frac{1}{N}\sum_{n=1}^N\left(\lambda_{n+1}^\alpha-
\lambda_{n}^\alpha\right)=\frac{\lambda_{N+1}^\alpha-\lambda_1^\alpha}{N},
\qquad N=1,2, \cdots\; .
\ee
In addition, if there exist constants $\gamma>0$ and $C>0$ such that
\be\label{nngap1}
\lim_{n\to+\infty} \frac{\lambda_{n}^\alpha}{n^\gamma}=C>0,
\ee
then the {\sl normalized  gaps} (or ``unfolding'' local
statistics in the physics literature) are defined as \cite{JMRR99,RZ02}
\be\label{normgap}
\delta_{\rm norm}^\alpha(N):= y_{N+1}^\alpha-y_{N}^\alpha,
\qquad N=1,2,\ldots\;,
%\frac{\lambda_{N+1}^\alpha-\lambda_{N}^\alpha}{\delta_{\rm ave}^\alpha(N)}=
%\frac{\lambda_{N+1}^\alpha-\lambda_{N}^\alpha}{\left(
%\lambda_{N+1}^\alpha-\lambda_1^\alpha\right)/N}=
%\frac{N(\lambda_{N+1}^\alpha-\lambda_{N}^\alpha)}{
%\lambda_{N+1}^\alpha-\lambda_1^\alpha}, \qquad
%N=1,2,3,\ldots;
\ee
where
\be\label{nngap2}
y_{n}^\alpha:=\left(\frac{\lambda_n^\alpha}{C}\right)^{1/\gamma},
\qquad n=1,2,\ldots\;.
\ee
Then an interesting question is to study their asymptotics, i.e.
the behaviour of $\delta_{\rm nn}^\alpha(N)$, $\delta_{\rm min}^\alpha(N)$, $\delta_{\rm ave}^\alpha(N)$ and $\delta_{\rm nrom}^\alpha(N)$  when $N\to +\infty$,
and another interesting and very challenging question is to study
the level spacing distribution $P_\alpha(s):=$ limiting distribution
of the normalized  gaps $\delta_{\rm norm}^\alpha(N)$,
which is defined as \cite{JMRR99,RZ02}
\be\label{Pas}
\frac{\#\left\{1\le n\le N\ | \ \delta_{\rm norm}^\alpha(n)< x\right\}}{N}
\stackrel{{N\to+\infty}}{\to} \int_0^x P_\alpha(s)ds, \qquad 0\le x <+\infty,
\ee
where $\#S$ denotes the number of elements in the set $S$.

 When $\alpha=2$ and $V(x)\equiv 0$ in \eqref{fproblem}, it collapses to a standard Sturm-Liouville eigenvalue problem of the Laplacian operator as
\begin{equation}\label{SO}
\begin{split}
L_{\rm SO}\; u(x)&:=-\partial_{xx}\, u(x)=-u^{\prime\prime}(x)
=\lambda\; u(x), \qquad x\in \Omega=(a,b), \\
u(a)&=u(b)=0.
\end{split}
\end{equation}
The eigenvalues and their corresponding eigenfunctions of \eqref{SO}
can be obtained analytically via the sine series as
\be
\lambda_n^{\alpha=2}=\left(\frac{n\pi}{b-a}\right)^2,
\qquad u_n(x)=\sin\left(\frac{n\pi(x-a)}{b-a}\right), \qquad n=1,2,\ldots \ .
\ee
These results immediately imply that the fundamental gap $\delta_{\rm fg}(\alpha=2)=\frac{3\pi^2}{(b-a)^2}$ and
\be\label{gapa2}
\begin{split}
&\delta_{\rm nn}^{\alpha=2}(N)=\left(\frac{(N+1)\pi}{b-a}\right)^2-
\left(\frac{N\pi}{b-a}\right)^2
=\frac{\pi^2}{(b-a)^2}(2N+1), \\
&\delta_{\rm min}^{\alpha=2}(N)\equiv \delta_{\rm nn}^{\alpha=2}(N=1)=\frac{3\pi^2}{(b-a)^2}, \\
&\delta_{\rm ave}^{\alpha=2}(N)=\frac{1}{N}\left[\left(\frac{(N+1)\pi}{b-a}\right)^2-
\left(\frac{\pi}{b-a}\right)^2\right]=\frac{\pi^2}{(b-a)^2}(N+2), \\
&\delta_{\rm norm}^{\alpha=2}(N)=y_{N+1}^{\alpha=2}-y_{N}^{\alpha=2}=N+1-N\equiv 1,
\end{split}
\qquad N=1,2,\ldots\; ;
\ee
where
\[y_{n}^{\alpha=2}=\sqrt{\lambda_n^{\alpha=2}/\left(\frac{\pi}{b-a}\right)^2}=
\sqrt{n^2}= n, \qquad n=1,2,\ldots\ .
\]
From the last equation in \eqref{gapa2}, one can immediately obtain
the level spacing distribution defined in \eqref{Pas} for $\alpha=2$ as
\be\label{pa22}
P_{\alpha=2}(s)=\delta(s-1),\qquad s\ge0,
\ee
where $\delta(\cdot)$ is the Dirac delta function.

When $\alpha=2$ and $V(x)\ne 0$ in \eqref{fproblem}, it collapses to a standard Sturm-Liouville eigenvalue problem, which has been extensively
studied in the literature. For analytical results, we refer to
\cite{LY83,LB02,H07} and references therein. For numerical methods and results, we refer to \cite{Boffi,BBC13,WT88} and references therein.

 When $0<\alpha<2$, in general, one cannot find the eigenvalues of the eigenvalue problem \eqref{fproblem} analytically and/or explicitly.
For mathematical theories of the eigenvalue problem \eqref{fproblem},
we refer to \cite{D12,K12} and references therein.
Some numerical methods have been proposed to solve \eqref{fproblem} numerically, including an asymptotic method was proposed in \cite{ZRK07},
a finite element method (FEM) \cite{BPM18} with piecewise linear element was presented in \cite{JLPR13} and a finite difference method (FDM) was studied in \cite{DZ15}. The FDM and FEM are usually first order
accurate when $0<\alpha<2$  and they can be adapted to compute
the first several eigenvalues \cite{JLPR13,DZ15,BPM18}. However, if we want to
calculate accurately and efficiently a very large number of eigenvalues, e.g. up to thousands or even millions eigenvalues, of
the eigenvalue problem \eqref{fproblem} in order to obtain a reliable
gaps distribution statistics, the FDM and FEM have severe drawbacks. The main aim of this paper is to propose a spectral method by using the generalized Jacobi functions for computing
different eigenvalue gaps and their distribution statistics of
the fractional eigenvalue problem related to FSO \eqref{fproblem}.
The proposed numerical method has at least two advantages:
(i) it is spectral accurate, and more importantly (ii) under a fixed number of degree of freedoms (DOF) $M$, it can calculate accurately a large number of eigenvalues with the number proportional to $M$. Thus this method is a very good
candidate for solving our problem, i.e. to compute eigenvalue gaps and their distribution statistics of the fractional eigenvalue problem \eqref{fproblem}.

\bigskip

Based on our extensive numerical results and observations, we speculate the following:

\smallskip

{\bf Conjecture} (Gaps and their distribution statistics of FSO in \eqref{fproblem} without potential) Assume $0<\alpha<2$ and $V(x)\equiv 0$ in \eqref{fproblem}, then we have the following asymptotics of its eigenvalues:
\be \label{asy-eig-1D}
\lambda_n^\alpha=\left(\frac{n\pi}{b-a}\right)^\alpha-
\left(\frac{\pi}{b-a}\right)^\alpha\frac{\alpha(2-\alpha)}{4} n^{\alpha-1} +O(n^{\alpha-2})=\lambda_{\rm loc}^\alpha(n)\left[1-\frac{\alpha(2-\alpha)}{4n}+O(n^{-2})\right], \ n\ge1,
\ee
where $\lambda_{\rm loc}^\alpha(n)=\left(\frac{n\pi}{b-a}\right)^\alpha$
($n=1,2,\ldots$) are the eigenvalues of the {\sl local fractional Laplacian
operator} on $\Omega=(a,b)$ with homogeneous Dirichlet boundary condition \cite{BR18}.
From \eqref{asy-eig-1D}, we obtain immediately the following approximations of different gaps:
\be\label{asy-gaps-1D}
\begin{split}
&\delta_{\rm nn}^\alpha(N)\approx\left(\frac{\pi}{b-a}\right)^\alpha\left[\alpha N^{\alpha-1}+\frac{\alpha(\alpha-1)(2+\alpha)}{4}N^{\alpha-2}
+O(N^{\alpha-3})\right],\ \  0<\alpha<2,\\
&\delta_{\rm min}^\alpha(N)=\left\{\ba{ll}
\delta_{\rm nn}^\alpha(1)=\lambda_2^\alpha-\lambda_1^\alpha, &1< \alpha<2,\\
\approx \delta_{\rm nn}^\alpha(1)=\lambda_2^\alpha-\lambda_1^\alpha, & \alpha=1,\\
\delta_{\rm nn}^\alpha(N)=\lambda_{N+1}^\alpha-\lambda_N^\alpha \approx \alpha \left(\frac{\pi}{b-a}\right)^\alpha N^{\alpha-1},  &0<\alpha<1,\\
\ea  \right.  \\
&\delta_{\rm ave}^\alpha(N)\approx\left(\frac{\pi}{b-a}\right)^\alpha
\left\{\ba{ll}\vspace{2mm}
\left[N^{\alpha-1}+\frac{\alpha(2+\alpha)}{4}N^{\alpha-2}+O(N^{-1})\right], &1<\alpha<2,\\ \vspace{2mm}
\left[1+\left(\frac{3}{4}-\frac{b-a}{\pi}\lambda_1^{\alpha=1}\right)N^{-1}+O(N^{-2})\right],
 &\alpha=1,\\
\left[N^{\alpha-1}-\left(\frac{b-a}{\pi}\right)^\alpha \lambda_1^\alpha N^{-1} +O(N^{\alpha-2})\right], &0<\alpha<1,\\
\ea\right.\\
&\delta_{\rm norm}^\alpha(N)\approx 1+O(N^{-2}),\qquad 0<\alpha<2,
%\left\{\ba{ll}\vspace{2mm}
%\alpha-\frac{\alpha(2+\alpha)}{4}N^{-1}+O(N^{-2})\le \alpha, &1<\alpha<2,\\ %\vspace{2mm}
%1+\left(\frac{b-a}{\pi}\lambda_1^{\alpha=1}-\frac{3}{4}\right)
%N^{-1}+O(N^{-2}),
% &\alpha=1,\\
%\alpha+\alpha\left(\frac{b-a}{\pi}\right)^\alpha \lambda_1^\alpha %N^{-\alpha} +O(N^{-1})\ge\alpha, &0<\alpha<1.\\
%\ea\right.
\end{split}
\qquad  N\ge1.
\ee
In addition, for the gaps distribution statistics defined
in \eqref{Pas},  we have
\be\label{gap-stat}
P_{\alpha}(s)=\delta(s-1), \qquad s\ge0, \qquad 0<\alpha\le2.
\ee
%In addition, for the fundamental gap, we have the following lower bound
%which is sharper than that in \cite{BRSS18}
%\be\label{fdgap1}
%\delta_{\rm fd}(\alpha):=\lambda_2^\alpha-\lambda_1^\alpha\ge %\frac{(2^\alpha-1)\pi^\alpha}{(b-a)^\alpha},
%\qquad 0<\alpha\le 2.
%\ee

\bigskip

The paper is organized as follows.
In Section \ref{sec:method}, we begin with some scaling properties of \eqref{fproblem} and propose a spectral-Galerkin method
by using the generalized Jacobi functions to discretize
the fractional eigenvalue problem \eqref{fproblem}.
In Section \ref{comp}, we test the accuracy and resolution capacity (or trust region) with respect to the DOF $M$ of the proposed Jacobi spectral method and compare it with the existing numerical methods such as FDM
and FEM. In Section \ref{sec:estimate}, we apply the proposed numerical method to study numerically asymptotics of different eigenvalue gaps and their distribution statistics of  \eqref{fproblem} without potential and formulate several interesting numerical
observations (or conjectures). Similar results are reported in Section \ref{sec:pot} for \eqref{fproblem} with potential.
Extensions of the numerical method and
results to the directional fractional Schr\"{o}dinger operator in high dimensions are presented in Section \ref{sec:2dgaps}.
Finally, some conclusions are drawn in Section \ref{sec:conclusion}.

\section{A Jacobi spectral method} \label{sec:method}
\setcounter{equation}{0}

In this section, we begin with a scaling argument to the problem \eqref{fproblem} so as to reduce it on a standard interval $(-1,1)$,
then reformulate it into a variational formulation and
discretize the problem by using the Jacobi spectral method.

\subsection{Scaling property}
Introduce
\be\label{trans}
x_0=\frac{a+b}{2},\qquad L=\frac{b-a}{2}, \qquad
\tilde{x} =\frac{x-x_0}{L},  \qquad
\tilde V(\tilde x) =L^\alpha V(x), \qquad
x\in \Omega=(a,b),
\ee
and consider the re-scaled fractional eigenvalue problem:

Find $\tilde \lambda \in \mathbb{R}$ and a real-valued function
$\tilde u(\tilde x)\ne0$ such that
\begin{equation}\label{sfep}
\begin{split}
\tilde L_{\rm FSO}\; \tilde u(\tilde x)&:=\left[(-\partial_{\tilde x \tilde x})^{\alpha/2}+\tilde V(\tilde x)\right]\tilde u(\tilde x)
=\tilde \lambda\; \tilde u(\tilde x), \qquad \tilde x\in \tilde \Omega:=(-1,1), \\
\tilde u(\tilde x)&=0, \qquad  \tilde x \in \tilde \Omega^c:=\mathbb{R} \backslash \tilde \Omega;
\end{split}
\end{equation}
then we have

\medskip

\begin{lemma} Let $\tilde \lmd$ be an eigenvalue of \eqref{sfep} and
$\tilde u:=\tilde u(\tilde x)$ be the corresponding eigenfunction, then
$\lmd= L^{-\alpha} \tilde \lmd$ is an eigenvalue of \eqref{fproblem}   and $u:=u(x)=\tilde u(\tilde x)=\tilde u\left(\frac{x-x_0}{L}\right)$
is the corresponding eigenfunction. Assume that
$0<\tilde \lambda_1^\alpha<\tilde \lambda_2^\alpha\le\ldots\le\tilde \lambda_n^\alpha\le\ldots$ are all eigenvalues of \eqref{sfep},
then $0<\lambda_1^\alpha<\lambda_2^\alpha\le\ldots\le \lambda_n^\alpha\le\ldots$ (ranked as in \eqref{eigenorder})
with $\lambda_n^\alpha=L^{-\alpha}\tilde \lambda_n^\alpha$ ($n=1,2,\ldots$) are all eigenvalues of  \eqref{fproblem}. In addition, we have the scaling property on the
different gaps as
\be\label{gapscl}
\begin{split}
&\delta_{\rm nn}^\alpha(N)=L^{-\ap}\tilde \delta_{\rm nn}^\alpha(N),
\qquad \hbox{with} \quad \tilde \delta_{\rm nn}^\alpha(N):=
\tilde\lambda_{N+1}^\alpha-\tilde \lambda_{N}^\alpha, \\
&\delta_{\rm min}^\alpha(N)=L^{-\ap}\tilde \delta_{\rm min}^\alpha(N), \qquad \hbox{with} \quad \tilde \delta_{\rm min}^\alpha(N):=\min_{1 \leq n \leq N}\ \tilde \delta_{\rm nn}^\alpha(n), \\
&\delta_{\rm ave}^\alpha(N)=L^{-\ap}\tilde \delta_{\rm ave}^\alpha(N), \qquad \hbox{with} \quad \tilde \delta_{\rm ave}^\alpha(N):=\frac{1}{N}\sum_{n=1}^N\tilde \delta_{\rm nn}^\alpha(n), \\
&\delta_{\rm norm}^\alpha(N)=\tilde \delta_{\rm norm}^\alpha(N),
\qquad \hbox{with} \quad \tilde \delta_{\rm norm}^\alpha(N):=
\tilde y_{N+1}^\alpha-\tilde y_{N}^\alpha, \quad \tilde y_{N}^\alpha
=\left(\frac{\tilde \lambda_N^\alpha}{L^\alpha C}\right)^{1/\gamma},
\end{split}
\quad N=1,2,\ldots\;;
\ee
which immediately imply that the level spacing distribution $P_\alpha(s)$
of \eqref{fproblem} does not  change under the rescaling \eqref{trans}, i.e.
the problems \eqref{fproblem} and \eqref{sfep} have the
same level spacing distribution.
\end{lemma}

\medskip

\noindent {\bf Proof:}
From \eqref{pvi} and noticing \eqref{trans}, a direct computation implies the scaling property of the  fractional Laplacian operator
\bea\label{dtapsc}
(-\partial_{xx})^{\alpha/2}\,u(x)&=&C_1^{\alpha}
\int_{\mathbb{R}}\frac{u(x)-u(y)}{|x-y|^{1+\alpha}}\,dy
= C_1^{\alpha}\int_{\mathbb{R}}\frac{u(x_0+L \tilde x)
-u(x_0+L \tilde y)}{|x_0+L \tilde x-x_0-L \tilde y|^{1+\alpha}}\,L\,{d\tilde y} \nonumber\\
&=&L^{-\alpha} C_1^{\alpha}\int_{\mathbb{R}}\frac{\tilde u(\tilde x)-
\tilde u(\tilde y)}{|\tilde x-\tilde y|^{1+\alpha}}\,d\tilde y
=L^{-\alpha}\,(-\partial_{\tilde x \tilde x})^{\alpha/2}\,\tilde u(\tilde x), \quad x\in\Omega,
\quad  \tilde x\in \tilde\Omega.
\eea
Noticing
\be
u(x)=0, \quad x\in \Omega^c \qquad \Longleftrightarrow \qquad \tilde u(\tilde x)=0, \quad \tilde x\in \tilde\Omega^c.
\ee
Substituting \eqref{dtapsc} into \eqref{sfep}, noting \eqref{fproblem}, we get
\bea
\tilde \lmd\,u(x)&=&\tilde \lmd \,\tilde u(\tilde x)=\left[(-\partial_{\tilde x\tilde x})^{\frac{\alpha}{2}}+\tilde V(\tilde x)\right]
\tilde u(\tilde x)=\left[L^{\alpha}\,(-\partial_{xx})^{\frac{\alpha}{2}} +\tilde V\left(\frac{x-x_0}{L}\right)\right]u(x)\nonumber\\
&=&L^{\alpha}
\left[(-\partial_{xx})^{\frac{\alpha}{2}} +L^{-\alpha} \tilde V\left(\frac{x-x_0}{L}\right)\right]
u(x)
=L^{\alpha}
\left[ (-\partial_{xx})^{\frac{\alpha}{2}}+V(x)\right]
u(x),\ x\in \Omega, \ \tilde x\in \tilde \Omega,\quad
\eea
which immediately implies that $u(x)$ is an eigenfunction of the operator $(-\partial_{xx})^{\frac{\alpha}{2}}+V(x)$ with the eigenvalue $ \lmd=L^{-\alpha}\tilde \lmd$.

 From the assumption  \eqref{eigenorder} with $\Omega=(-1,1)$ that
 $0<\tilde \lambda_1^\alpha<\tilde \lambda_2^\alpha\le\ldots\le\tilde \lambda_n^\alpha\le\ldots$ are all eigenvalues of \eqref{sfep}, we get immediately that $0<\lambda_1^\alpha<\lambda_2^\alpha\le\ldots\le\lambda_n^\alpha\le\ldots$
 with $\lambda_n^\alpha=L^{-\alpha}\tilde \lambda_n^\alpha$ ($n=1,2,\ldots$) are all eigenvalues of the eigenvalue problem \eqref{fproblem}. Then the scaling property on the different gaps
 \eqref{gapscl} can be obtained straightforward by using $\tilde \lambda_n^\alpha=L^\alpha \lambda_n^\alpha$ ($n=1,2,\ldots$).
\hfill $\Box$

\subsection{A variational formulation}

Following those in the literature \cite{LM72,G85}, we introduce the fractional functional space $H^{\frac{\alpha}{2}}(\mathbb{R})$ through the Fourier transform
\begin{equation}
H^{\frac{\alpha}{2}}(\mathbb{R})=\left\{v\in \mathcal{D}'(\mathbb{R})\ | \ \|v\|_{\frac{\alpha}{2},\mathbb{R}}<\infty\right \},
\ee
where the norms are defined as
\be
|v|_{\frac{\alpha}{2},\mathbb{R}}
=\left(\int_{\mathbb{R}}
\left|\xi|^{\alpha}\right|\; \left|(\mathcal{F}v)(\xi)\right|^2 \,d\xi\right)^{\frac{1}{2}},\qquad
\|v\|_{\frac{\alpha}{2},\mathbb{R}}
=\left(\int_{\mathbb{R}}(1+|\xi|^2)^{\frac{\alpha}{2}}
\left|(\mathcal{F}v)(\xi)\right|^2 d\xi\right)^{\frac{1}{2}};
\end{equation}
and then the fractional functional space $H^{\frac{\alpha}{2}}(\Omega)$ can be obtained from $H^{\frac{\alpha}{2}}(\mathbb{R})$ by extension \cite{LM72,G85}
\begin{equation}
H^{\frac{\alpha}{2}}(\Omega)=\left\{v: \; \Omega \to
\mathbb{R}
\ | \ \hat{v}=E_\Omega v \in H^{\frac{\alpha}{2}}(\mathbb{R})\right\},
\ee
where the norms are defined as
\be
|v|_{\frac{\alpha}{2}}:=|v|_{\frac{\alpha}{2},\Omega}=|E_\Omega v|_{\frac{\alpha}{2},\mathbb{R}},\qquad
\|v\|_{\frac{\alpha}{2}}:=\|v\|_{\frac{\alpha}{2},\Omega}=\|E_\Omega v\|_{\frac{\alpha}{2},\mathbb{R}},
\qquad \forall v\in H^{\frac{\alpha}{2}}(\Omega),
\end{equation}
with $\hat{v}=E_\Omega v: \mathbb{R} \to \mathbb{R}$ (extension of $v$ from
$\Omega$ to $\mathbb{R}$) defined as
\be
\hat{v}(x)=(E_\Omega v)(x)=\left\{\ba{ll}
v(x), & x\in\Omega,\\
0, &x\in \mathbb{R}\backslash \Omega.\\
\ea\right.
\ee

For any $v\in H^{\frac{\alpha}{2}}(\Omega)$, multiplying $v$ to \eqref{fproblem} and integrating
over $\Omega$ and using integration by parts, we
obtain the variational (or weak) formulation of the fractional eigenvalue
problem \eqref{fproblem} as:

find $\lambda\in \mathbb{R}$ and $0\ne u \in
H^{\frac{\alpha}{2}}(\Omega)$, such that
\begin{equation} \label{weakf}
a(u,v)= \lambda\; b(u,v), \qquad \forall v\in H^{\frac{\alpha}{2}}(\Omega),
\ee
where the bilinear forms $a(\cdot,\cdot)$ and $b(\cdot,\cdot)$ are given as
\begin{equation}
\begin{split}
&a(u,v)=\int_{\Omega} \left[(-\partial_{xx})^{\frac{\alpha}{2}}u+
V(x)u\right]v dx=\int_{\Omega} \left[(-\partial_{xx})^{\frac{\alpha}{4}}u\;
(-\partial_{xx})^{\frac{\alpha}{4}}v+
V(x)uv\right] dx,\\
&b(u,v)=\int_{\Omega}u(x) v(x) dx, \qquad \forall u,v\in H^{\frac{\alpha}{2}}(\Omega).
\end{split}
\end{equation}

\subsection{A spectral discretization by using the Jacobi functions}
Since we are mainly interested in gaps and their distribution statistics,
from the results in Lemma 2.1, without loss of generality, from now on,
we always assume that $\Omega=(-1,1)$, i.e. $a=-1$ and $b=1$ in
\eqref{fproblem}.

Let $\{P_n^{\frac{\alpha}{2},\frac{\alpha}{2}}(x)\}_{n=0}^\infty$ denote the classical Jacobi polynomials (or Gegenbauer polynomials) which are orthogonal with
respect to the weight function $\omega^{\frac{\alpha}{2},\frac{\alpha}{2}}(x)=(1-x^2)^{\frac{\alpha}{2}}$ over the interval $(-1,1)$, i.e.
\be
\int_{-1}^1 P_n^{\frac{\alpha}{2},\frac{\alpha}{2}}(x)\,
P_m^{\frac{\alpha}{2},\frac{\alpha}{2}}(x)\,
\omega^{\frac{\alpha}{2},\frac{\alpha}{2}}(x)dx =C_n \delta_{nm},
\qquad n,m=0,1,2,\ldots,
\ee
where $\delta_{nm}$ is the  kronecker delta and
\be
C_n=\frac{2^{\alpha+1}}{2n+\alpha+1}\frac{\Gamma(n+\alpha/2+1)^2}
{\Gamma(n+\alpha+1)n!} \qquad n=0,1,2\ldots\; .
\ee
Define the generalized Jacobi functions
\be \label{gJf11}
\mathcal{J}_n^{-\frac{\alpha}{2},-\frac{\alpha}{2}}(x)
 =(1-x^2)^{\frac{\alpha}{2}}P_n^{\frac{\alpha}{2},\frac{\alpha}{2}}(x)=
 \omega^{\frac{\alpha}{2},\frac{\alpha}{2}}(x)\,
 P_n^{\frac{\alpha}{2},\frac{\alpha}{2}}(x),
 \quad -1\le x\le 1, \qquad n=0,1,2,\ldots,
\ee
then by Theorem 2 in Ref. \cite{MCS16}, we have
\begin{equation} \label{derpol}
(-\partial_{xx})^{\frac{\alpha}{2}} \mathcal{J}_n^{-\frac{\alpha}{2},-\frac{\alpha}{2}}(x)= \frac{\Gamma(n+\alpha+1)}{n!} P^{\frac{\alpha}{2},\frac{\alpha}{2}}_n(x),
\quad -1< x< 1, \qquad n=0,1,2,\ldots\;.
\end{equation}
Combining \eqref{gJf11} and \eqref{derpol}, we obtain
\begin{equation} \label{jnmoth}
\begin{split}
&\int_{-1}^1 (-\partial_{xx})^{\frac{\alpha}{2}} \mathcal{J}_n^{-\frac{\alpha}{2},-\frac{\alpha}{2}}(x)\; \mathcal{J}_m^{-\frac{\alpha}{2},-\frac{\alpha}{2}}(x)\;dx
=\int_{-1}^1 \mathcal{J}_n^{-\frac{\alpha}{2},-\frac{\alpha}{2}}(x)\; (-\partial_{xx})^{\frac{\alpha}{2}} \mathcal{J}_m^{-\frac{\alpha}{2},-\frac{\alpha}{2}}(x)\;dx\\
&\ =\int_{-1}^1 (-\partial_{xx})^{\frac{\alpha}{4}} \mathcal{J}_n^{-\frac{\alpha}{2},-\frac{\alpha}{2}}(x)\; (-\partial_{xx})^{\frac{\alpha}{4}} \mathcal{J}_m^{-\frac{\alpha}{2},-\frac{\alpha}{2}}(x)\;dx
=\int_{-1}^1 \frac{\Gamma(n+\alpha+1)}{n!} P^{\frac{\alpha}{2},\frac{\alpha}{2}}_n(x)\; \mathcal{J}_m^{-\frac{\alpha}{2},-\frac{\alpha}{2}}(x)\;dx\\
&\ =\frac{\Gamma(n+\alpha+1)}{n!} \int_{-1}^1 P^{\frac{\alpha}{2},\frac{\alpha}{2}}_n(x)\;
P_m^{\frac{\alpha}{2},\frac{\alpha}{2}}(x)\,
\omega^{\frac{\alpha}{2},\frac{\alpha}{2}}(x)\;dx \\
&\ =\frac{2^{\alpha+1}\Gamma(n+\alpha/2+1)^2} {(n!)^2(2n+\alpha+1)}\delta_{nm}, \quad n,m=0,1,2\ldots\; .
\end{split}
\end{equation}
Introduce
\be\label{phinx}
\phi_n(x):=\frac{\sqrt{2n+\alpha+1}n!}{2^{\alpha/2+1/2}\Gamma(n+\alpha/2+1)} \mathcal{J}_n^{-\frac{\alpha}{2},-\frac{\alpha}{2}}(x),
\quad -1\le x \le 1, \qquad n=0,1,2,\ldots\; .
\ee
Let $M>0$ be a positive integer and define the
finite dimensional space (which is an approximate
subspace of $H^{\frac{\alpha}{2}}(\Omega)$) as
\begin{equation}
\mathbb{W}_M:
=\mbox{span}\left\{\phi_m(x),\ 0\leq m \leq M-1\right\},
\end{equation}
then a Jacobi spectral method ({\bf JSM}) for \eqref{weakf} is given as:

Find $\lambda_M\in \mathbb{R}$ and  $0\ne u_M \in \mathbb{W}_M$ such that
\begin{equation} \label{weakn}
a(u_M,v_M)=\lambda_M\, b(u_M,v_M), \qquad \forall v_M \in \mathbb{W}_M.
\end{equation}

In order to cast the eigenvalue problem \eqref{weakn} into matrix form,
we express $u_M\in \mathbb{W}_M$ as a combination of the  basis functions as
\begin{equation} \label{expression}
u_M(x)=\sum^{M-1}_{m=0}\hat{u}_{m}\,\phi_m(x), \qquad -1\le x\le 1.
\end{equation}
Plugging \eqref{expression} into \eqref{weakn} and noticing
\eqref{jnmoth}, after some detailed computation, we obtain the following
standard matrix eigenvalue problem:
\begin{equation} \label{matrixf}
\left({\bf I}_{M}+\mathbf{V}\right)\hat{U}=\lambda_M\; \mathbf{B}\;\hat{U},
\end{equation}
where $\hat{U}=(\hat{u}_0,\hat{u}_1,\cdots,\hat{u}_{M-1})^T\in \mathbb{R}^{M}$
is the eigenvector, ${\bf I}_{M}$ is the $M\times M$ identity matrix, and $\mathbf{V}=(v_{nm})_{0\leq n,m\leq M-1}\in \mathbb{R}^{M\times M}$ and $\mathbf{B}=(b_{nm})_{0\leq n,m\leq M-1}\in \mathbb{R}^{M\times M}$ are given as
\begin{equation} \label{matrem1}
\begin{split}
&v_{nm}=\int_{-1}^1 V(x)\phi_n(x)\phi_m(x)dx, \\
&b_{nm}=\int_{-1}^1 \phi_n(x)\phi_m(x)dx,
\end{split}\qquad n,m=0,1,\ldots,M-1.
\end{equation}
Plugging \eqref{phinx} into the second equation
in \eqref{matrem1}, after a detailed computation, we get
\begin{equation}
b_{nm}=\left\{
\begin{array}{ll}
\displaystyle  \frac{(-1)^{\frac{n-m}{2}}\sqrt{\pi(2n+\alpha+1)(2m+\alpha+1)}\Gamma(\alpha+1)(n+m)!}
{2^{\alpha+n+m+1}\Gamma(\alpha+\frac{n+m}{2}+\frac{3}{2})
\Gamma(\frac{\alpha}{2}+\frac{n-1}{2}+1)
\Gamma(\frac{\alpha}{2}+\frac{m-1}{2}+1)(\frac{n+m}{2})!}, & \hbox{$n+m$~even,} \\
0, & \hbox{$n+m$~odd.}
\end{array}
\right.
\end{equation}
If $V(x)\equiv 0$, then $\mathbf{V}={\bf 0}$. Of course, if
$V(x)\ne 0$, then the integrals in the first equation in
\eqref{matrem1} can be computed numerically via numerical
quadratures with spectral accuracy \cite{QV94,BM97}.
Finally the matrix eigenvalue problem \eqref{matrixf}
can be solved numerically by the standard eigenvalue
solvers such as QR-method \cite{MB08}.

\bigskip

We remark here that different numerical methods have been proposed in the literature for discretizing the fractional Laplacian operator $(-\partial_{xx})^{\alpha/2}$ via
the formulation \eqref{pvi} or \eqref{foperator} or their equivalent forms
for numerical simulation of fractional partial differential equations,
see \cite{LCTBA13,ZSP14,ATZ16,MCS16,OA16,R17,CML18} and references therein. In fact, a method to discretize
the fractional Laplacian operator $(-\partial_{xx})^{\alpha/2}$ can
directly generate a method to solve the fractional eigenvalue problem
\eqref{fproblem}. For example, a finite element method ({\bf FEM}) with piecewise linear elements  was proposed and analyzed in \cite{JLPR13,BPM18} for computing the eigenvalues of \eqref{fproblem}. Similarly, if we adopt the
standard finite difference method  to discretize
the fractional Laplacian operator $(-\partial_{xx})^{\alpha/2}$ \cite{CD12,LCTBA13}
in \eqref{fproblem}, we can obtain a finite difference method ({\bf FDM}) for computing
the eigenvalues of \eqref{fproblem}. The details are omitted here for brevity.

\section{Accuracy and comparison with existing methods}\label{comp}
\setcounter{equation}{0}

In this section, we test the accuracy and resolution capacity of the
Jacobi specral method ({\sl JSM}) presented in the previous section and compare it with the fractional centered finite difference method ({\sl FDM}) proposed in \cite{ZSP14,CD12} and the finite element method ({\sl FEM}) with piecewise linear element proposed in \cite{JLPR13} for the eigenvalue problem \eqref{fproblem} with $\Omega=(-1,1)$. The `exact' eigenvalues $\lambda_n^\alpha$ ($n=1,2,\ldots$) are
obtained numerically by using the JSM \eqref{weakn} under a very large DOF $M=M_0$, e.g. $M_0=12800$. Let $\lambda_{n,M}^\alpha$ be the numerical approximation of $\lambda_n^\alpha$ ($n=1,2,\ldots, M$) obtained by a numerical method
with the DOF chosen as $M$.
Define the absolute and relative errors of $\lambda_n^\alpha$ as
\be\label{errors}
e_n^\alpha:=\left|\lambda_n^\alpha-\lambda_{n,M}^\alpha\right|,
\qquad  e_{n,r}^\alpha:=\frac{\left|\lambda_n^\alpha-\lambda_{n,M}^\alpha\right|}
{\lambda_n^\alpha}, \qquad n=1,2,\ldots,
\ee
respectively.

\subsection{Accuracy test} \label{sec:accuracy}

We first test the convergence rates of different numerical methods for
the eigenvalue problem \eqref{fproblem} including the
JSM \eqref{weakn}, FEM \cite{JLPR13,BPM18} and FDM \cite{ZSP14,CD12,DZ15}. %In order to %do so, we take   and $\Omega=(-1,1)$ in \eqref{fproblem}.
Table \ref{errors1} displays the absolute errors of computing
the first eigenvalue of \eqref{fproblem} with $V(x)\equiv0$ and  different $\alpha$ by using our JSM \eqref{weakn},
FEM \cite{JLPR13} and FDM \cite{ZSP14,CD12}; and
Table \ref{errors2} lists the absolute errors of computing
the first, second, fifth and tenth eigenvalues of \eqref{fproblem}
with $\alpha=0.5$ and $V(x)\equiv0$ by using those  methods.
For comparison with existing results, Table \ref{eig} lists
the first three eigenvalues of \eqref{fproblem} with $V(x)\equiv0$ and  different $\alpha$ obtained by using our JSM \eqref{weakn}
under the DOF $M=160$ and the asymptotic method in \cite{ZRK07} under the DOF $M=5000$.
Figure \ref{fig:rate} shows convergence rates
of our JSM \eqref{weakn}  for computing the first, second, fifth and tenth eigenvalues of \eqref{fproblem} with $V(x)\equiv0$ and  different
$\alpha$; and Figure \ref{fig:rate1} lists similar results
of \eqref{fproblem} with $V(x)=\frac{x^2}{2}$ and  different
$\alpha$.

\begin{table}[h!]
\centering
\begin{tabular}{ |c|c|c |c|c|c|c|c|c|l|} \hline
& &$M=2$ &$M=4$ &$M=8$ &$M=16$ &$M=32$ &$M=64$ &$M=128$ &$M=256$\\ \hline
              &JSM &3.63E-5  &8.47E-9 &1.36E-12 & 1.36E-12&1.39E-12&1.40E-12&1.17E-12 &3.62E-12 \\
$\alpha=2.0$  &FEM &5.32E-1  &1.29E-1 &3.18E-2  & 7.92E-3 &1.97E-3 &4.87E-4 &1.16E-4  &2.32E-5\\
              &FDM &4.67E-1  &1.24E-1 &3.15E-2  & 7.90E-3 &1.97E-3 &4.87E-4 &1.16E-4  &2.32E-5\\
\hline
              &JSM &3.18E-5  &1.68E-8 &1.78E-11& 2.49E-12&2.55E-12&2.24E-12&3.08E-12 &2.12E-12 \\
$\alpha=1.95$ &FEM &4.96E-1  &1.16E-2 &2.79E-2 & 6.86E-3 &1.72E-3 &4.49E-4 &1.24E-4  &3.78E-5\\
              &FDM &2.31E-1  &2.86E-2 &5.16E-3 & 5.41E-4 &2.75E-5 &7.56E-6 &3.76E-6  &1.18E-6\\
\hline
              &JSM &2.31E-6  &7.17E-7 &1.57E-8 & 1.72E-10&2.16E-12&1.02E-12&6.64E-13 &1.41E-12 \\
$\alpha=1.5$  &FEM &2.72E-1  &6.86E-2 &2.55E-2 & 1.18E-2 &5.86E-3 &2.96E-3 &1.49E-3  &7.53E-4\\
              &FDM &9.15E-2  &6.78E-2 &5.41E-2 & 3.21E-2 &1.73E-2 &9.01E-3 &4.59E-3  &2.31E-3\\
\hline
              &JSM &2.16E-5  &6.32E-6 &3.56E-7 & 1.15E-8&2.65E-10&4.67E-12&5.94E-13 &5.53E-13 \\
$\alpha=1.0$  &FEM &1.66E-1  &5.97E-2 &2.29E-2 & 1.51E-2&7.83E-3 &4.01E-3 &2.03E-3  &1.01E-3\\
              &FDM &1.15E-1  &1.00E-1 &6.03E-2 & 3.28E-2&1.71E-2 &8.77E-3 &4.44E-3  &2.24E-3\\
\hline
              &JSM &1.22E-4  &3.14E-5 &3.95E-6 & 3.65E-7&2.80E-8 &1.94E-9 &1.26E-10 &7.10E-12 \\
$\alpha=0.5$  &FEM &8.74E-2  &3.93E-2 &2.03E-2 & 1.06E-2&5.54E-3 &2.84E-3 &1.45E-3  &7.35E-4\\
              &FDM &1.08E-1  &7.00E-2 &3.87E-2 & 2.04E-2&1.05E-2 &5.40E-3 &2.74E-3  &1.38E-3\\
\hline
              &JSM &1.29E-4  &4.01E-5 &8.58E-6 & 1.57E-6&2.68E-7 &4.49E-8 &7.36E-9  &1.06E-9 \\
$\alpha=0.1$  &FEM &2.02E-2  &1.01E-2 &5.27E-3 & 2.75E-3&1.42E-3 &7.30E-4 &3.72E-4  &1.89E-4\\
              &FDM &3.12E-2  &1.80E-2 &9.59E-3 & 4.99E-3&2.56E-3 &1.31E-3 &6.65E-4  &3.36E-4\\
\hline
\end{tabular}
\caption{Absolute errors of computing the first eigenvalue of \eqref{fproblem} with $\Omega=(-1,1)$, $V(x)\equiv0$ and different $\alpha$ by using our JSM \eqref{weakn}, FEM \cite{JLPR13} and FDM \cite{ZSP14,CD12}}
\label{errors1}
\end{table}

\begin{table}[h!]
\centering
\begin{tabular}{ |c|c|c |c|c|c|c|c|c|l|} \hline
& &$M=2$ &$M=4$ &$M=8$ &$M=16$ &$M=32$ &$M=64$ &$M=128$ &$M=256$\\ \hline
              &JSM &1.22E-4  &3.14E-5 &3.95E-6 & 3.65E-7&2.80E-8 &1.94E-9 &1.26E-10 &7.10E-12 \\
$e_1^\alpha$  &FEM &8.74E-2  &3.93E-2 &2.03E-2 & 1.06E-2&5.54E-3 &2.84E-3 &1.45E-3  &7.35E-4\\
              &FDM &1.08E-1  &7.00E-2 &3.87E-2 & 2.04E-2&1.05E-2 &5.40E-3 &2.74E-3  &1.38E-3\\
\hline
              &JSM &NA  &1.88E-4 &2.54E-5 & 2.03E-6&1.41E-7 &9.29E-9 &5.90E-10 &3.42E-11 \\
$e_2^\alpha$  &FEM &NA  &8.03E-2 &3.10E-2 & 1.59E-2&8.49E-3 &4.46E-3 &2.31E-3  &1.18E-3\\
              &FDM &NA  &2.54E-2 &4.02E-2 & 2.71E-2&1.55E-2 &8.36E-3 &4.35E-3  &2.23E-3\\
\hline
              &JSM &NA  &NA &2.14E-3 & 7.30E-6&5.89E-7 &4.14E-8 &2.73E-9  &1.16E-10 \\
$e_5^\alpha$  &FEM &NA  &NA &1.26E-1 & 3.05E-2&1.33E-2 &6.91E-3 &3.66E-3  &1.91E-3\\
              &FDM &NA  &NA &1.19E-2 & 3.88E-3&1.13E-3 &3.10E-4 &8.17E-5  &2.10E-5\\
\hline
                &JSM &NA  &NA &NA & 1.02E-2&1.92E-6 &1.31E-7 &8.44E-9  &5.01E-10 \\
$e_{10}^\alpha$ &FEM &NA  &NA &NA & 1.41E-1&2.66E-2 &9.96E-3 &5.00E-3  &2.63E-3\\
                &FDM &NA  &NA &NA & 2.14E-3&5.99E-4 &1.59E-4 &4.14E-5  &1.06E-5\\
\hline
\end{tabular}
\caption{Absolute errors of computing the first, second, fifth and tenth eigenvalues of \eqref{fproblem} with $\Omega=(-1,1)$, $\alpha=0.5$ and $V(x)\equiv0$  by using our JSM \eqref{weakn}, FEM \cite{JLPR13} and FDM \cite{ZSP14,CD12}}
\label{errors2}
\end{table}

\begin{table}[h!]
\centering
\begin{tabular}{ |l |l|l|l|l|r|r|} \hline
 & \multicolumn{2}{|c|}{$\lambda_1^\alpha$}& \multicolumn{2}{|c|}{$\lambda_2^\alpha$}& \multicolumn{2}{|c|}{$\lambda_3^\alpha$} \\ \hline
 &JSM \eqref{weakn} &Ref. \cite{ZRK07} &JSM  \eqref{weakn} &Ref. \cite{ZRK07}&JSM \eqref{weakn} &Ref. \cite{ZRK07} \\ \hline
$\alpha=1.99$  &2.443691434   &2.442 &9.73318159   &9.729& 21.82868373       &21.829 \\
$\alpha=1.9$   &2.244059359   &2.243 &8.59575252   &8.593& 18.71689400       &18.718\\
$\alpha=1.8$   &2.048734983   &2.048 &7.50311692   &7.501& 15.79989416       &15.801\\
$\alpha=1.5$   &1.597503545   &1.597 &5.05975992   &5.059&  9.59430576       & 9.957\\
$\alpha=1.0$   &1.157773883   &1.158 &2.75475474   &2.754&  4.31680106       & 4.320 \\
$\alpha=0.5$   &0.970165419   &0.970 &1.60153773   &1.601&  2.02882105       & 2.031 \\
$\alpha=0.2$   &0.957464477   &0.957 &1.19653989   &1.197&  1.31909097       & 1.320\\
$\alpha=0.1$   &0.972594401   &0.973 &1.09219649   &1.092&  1.14732244       & 1.148\\
$\alpha=0.01$  &0.996634628   &0.997 &1.00871791   &1.009&  1.01374130       & 1.014\\
\hline
\end{tabular}
\caption{The first three eigenvalues of \eqref{fproblem} with $\Omega=(-1,1)$, $V(x)\equiv0$ and  different $\alpha$ obtained numerically by our JSM \eqref{weakn}
under the DOF $M=160$ and the asymptotic method in \cite{ZRK07} under the DOF $M=5000$.}
\label{eig}
\end{table}

\begin{figure}[h!]
\centerline{
\psfig{figure=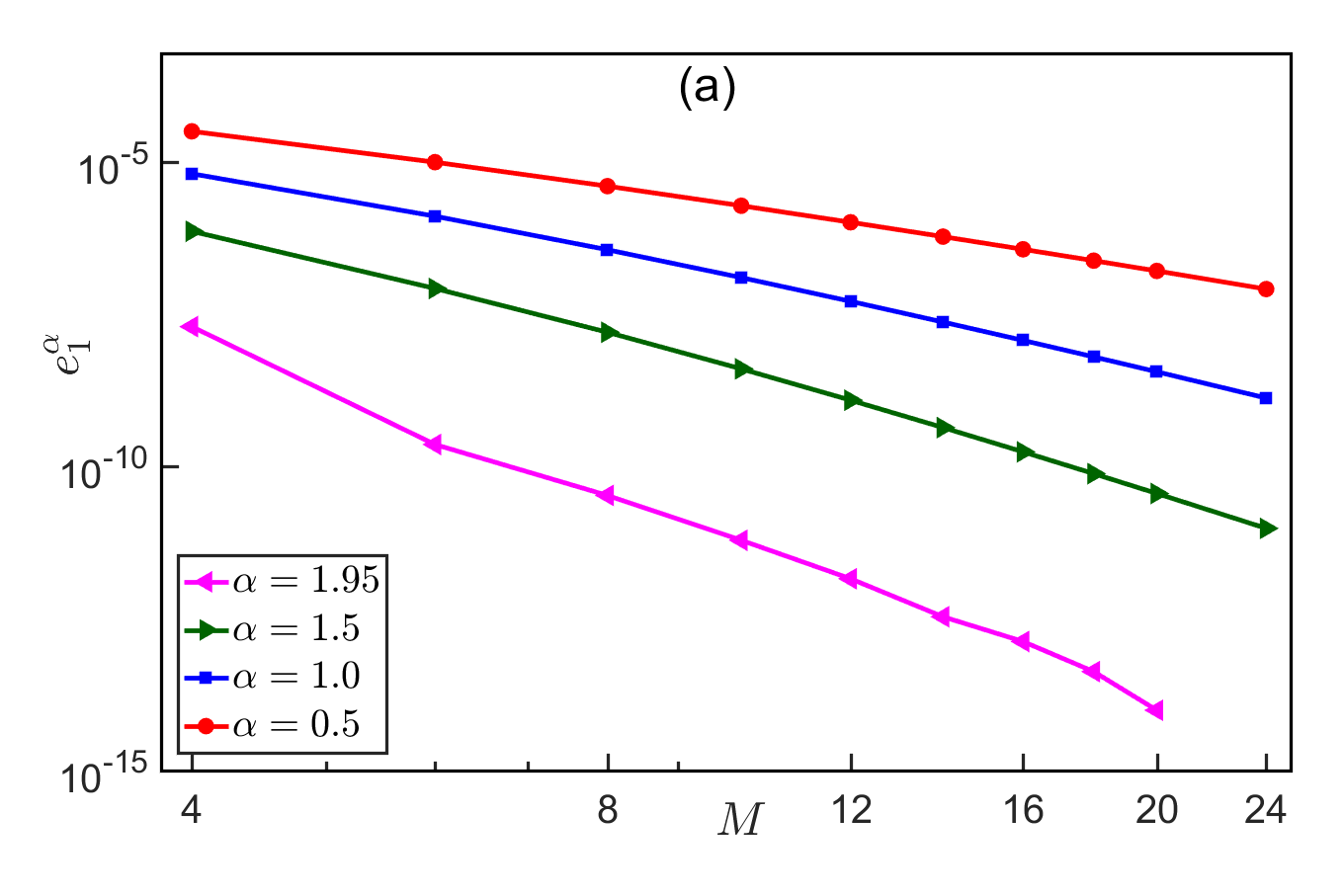,height=5cm,width=6cm,angle=0}\qquad
\psfig{figure=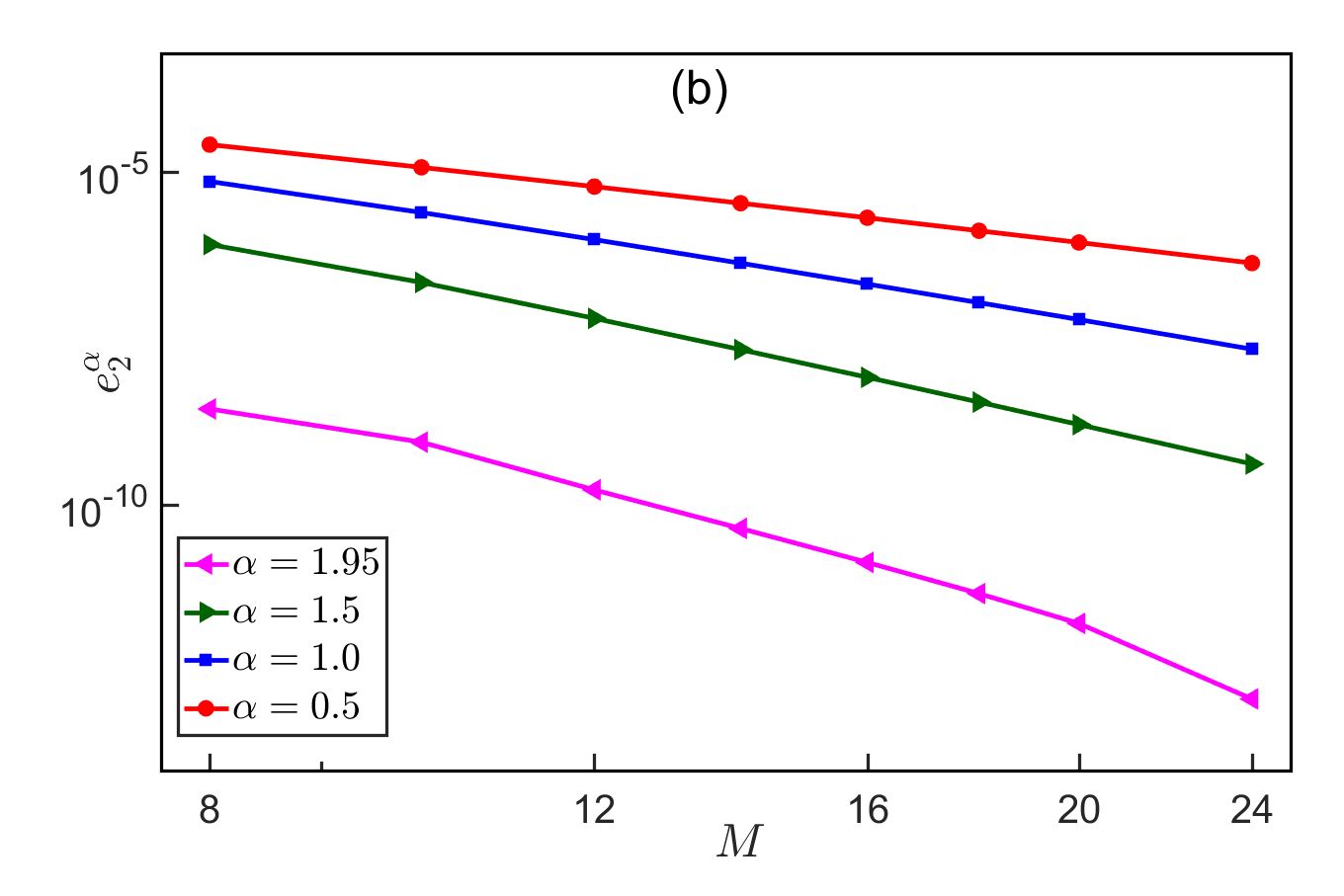,height=5cm,width=6cm,angle=0}}
\centerline{
\psfig{figure=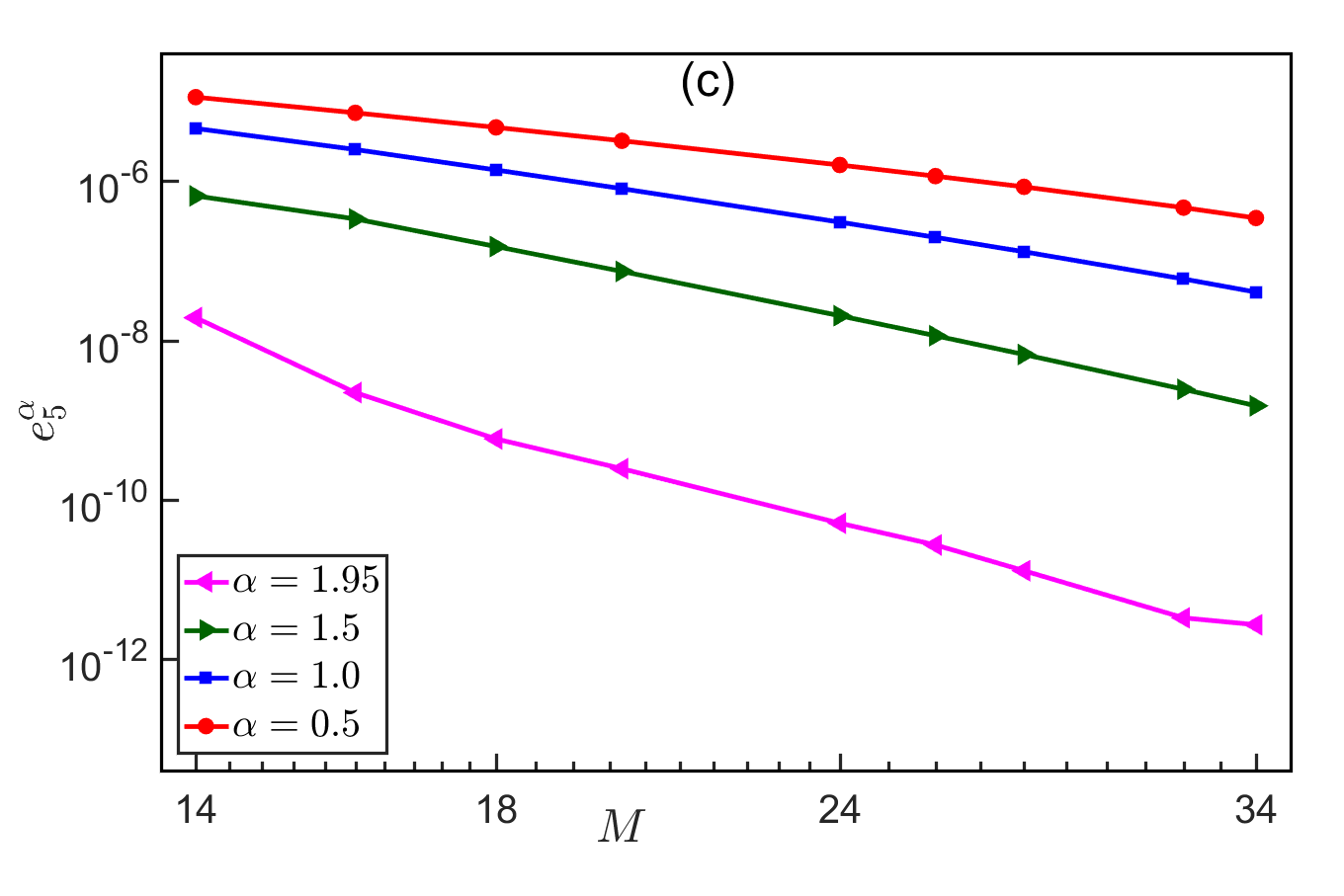,height=5cm,width=6cm,angle=0}\qquad
\psfig{figure=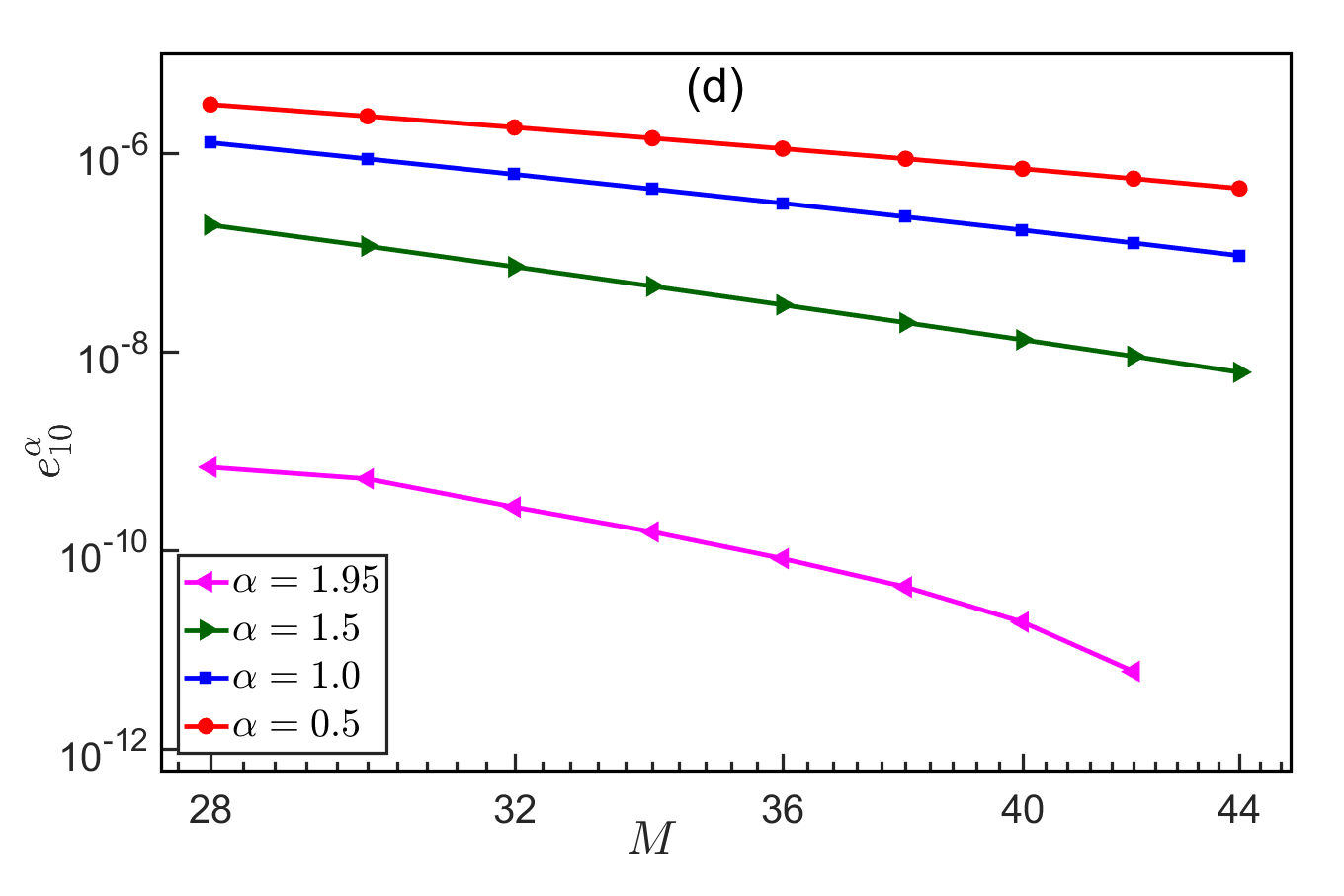,height=5cm,width=6cm,angle=0}}
\caption{Convergence rates of computing different eigenvalues of \eqref{fproblem} with $\Omega=(-1,1)$, $V(x)\equiv0$ and  different $\alpha$
by using our JSM \eqref{weakn} for:
(a) the first eigenvalue $\lambda_1^\alpha$, (b) the second eigenvalue $\lambda_2^\alpha$, (c) the fifth eigenvalue $\lambda_5^\alpha$,
and (d) the tenth eigenvalue $\lambda_{10}^\alpha$.}
\label{fig:rate}
\end{figure}

\begin{figure}[h!]
\centerline{
\psfig{figure=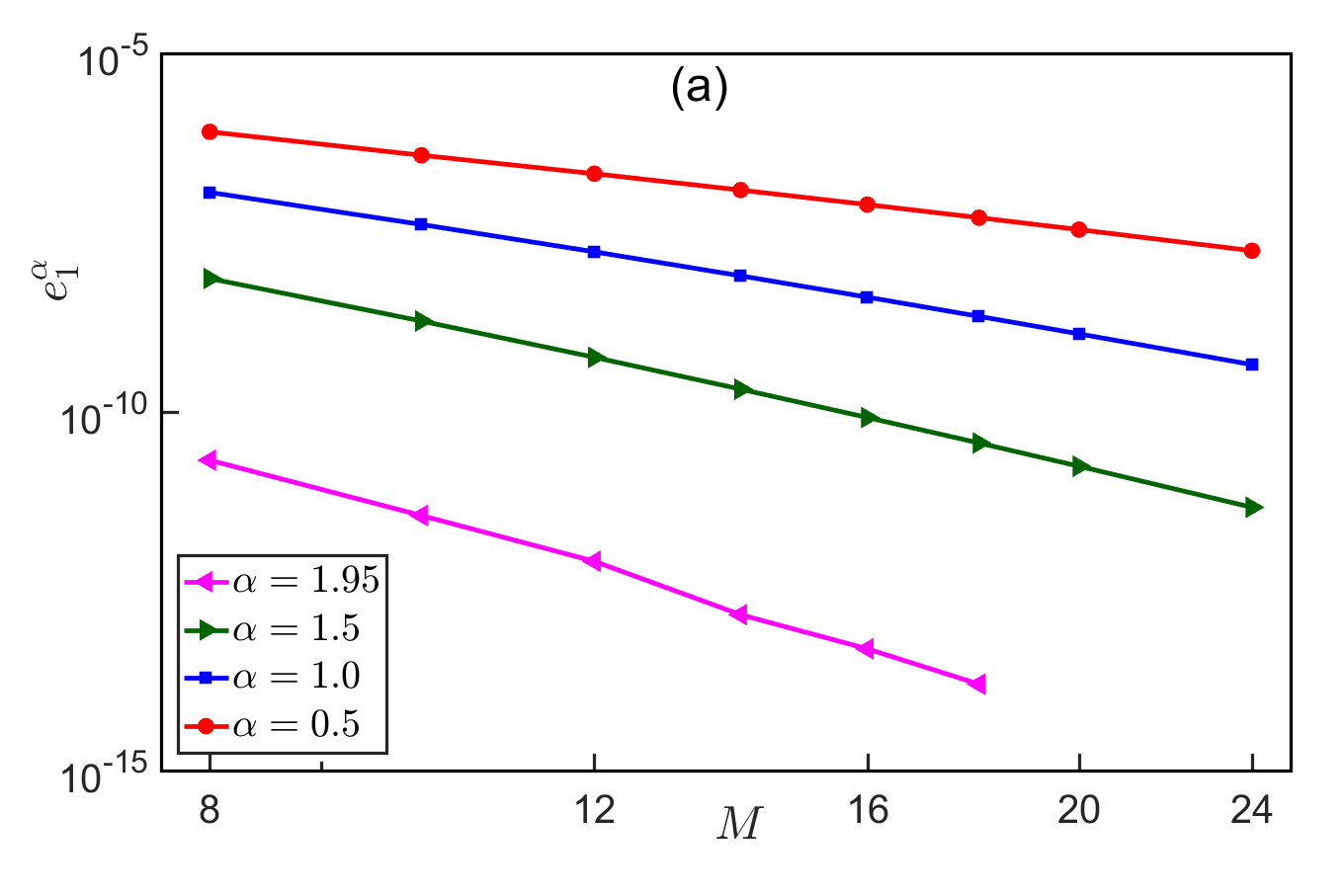,height=5cm,width=6cm,angle=0}\qquad
\psfig{figure=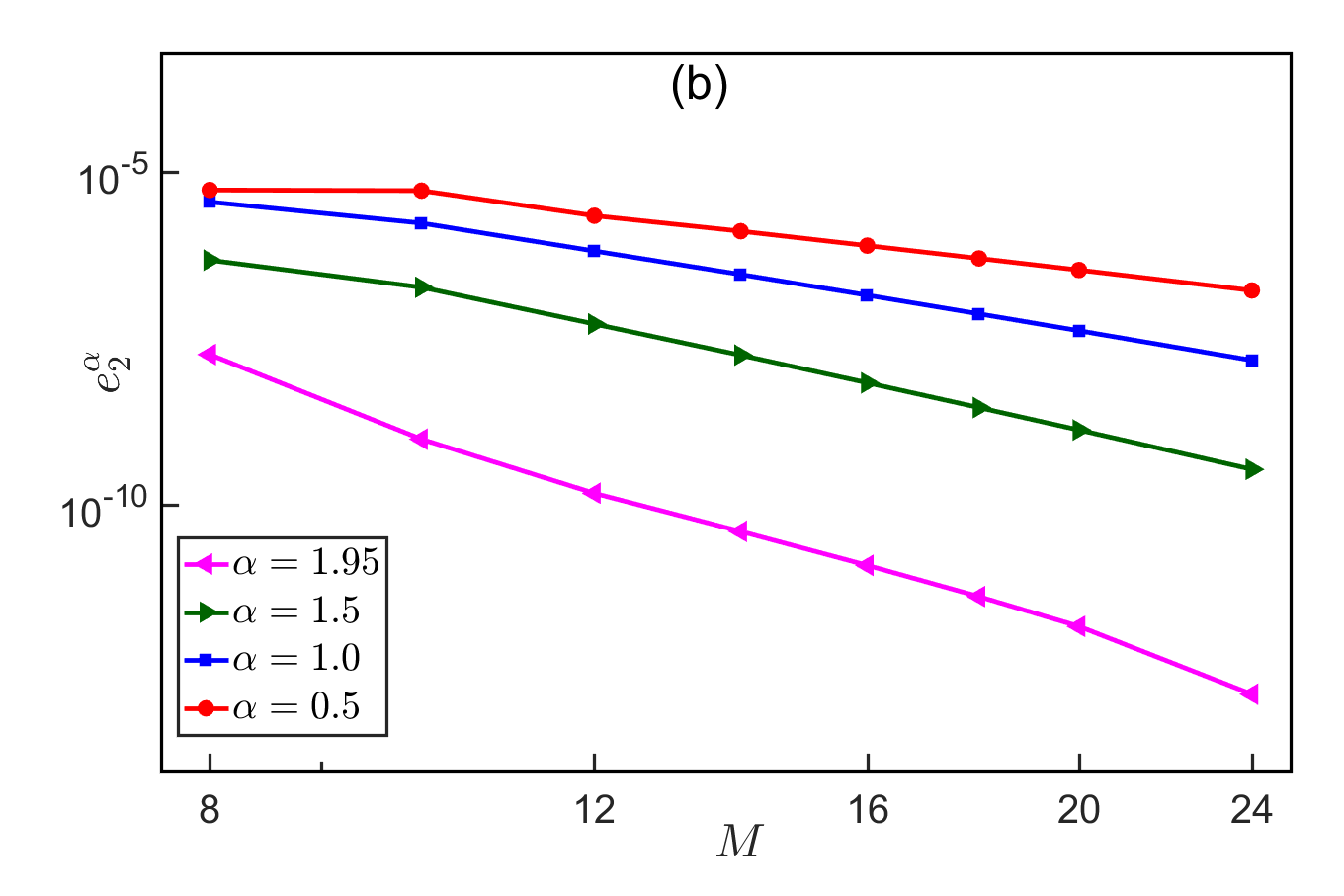,height=5cm,width=6cm,angle=0}}
\centerline{
\psfig{figure=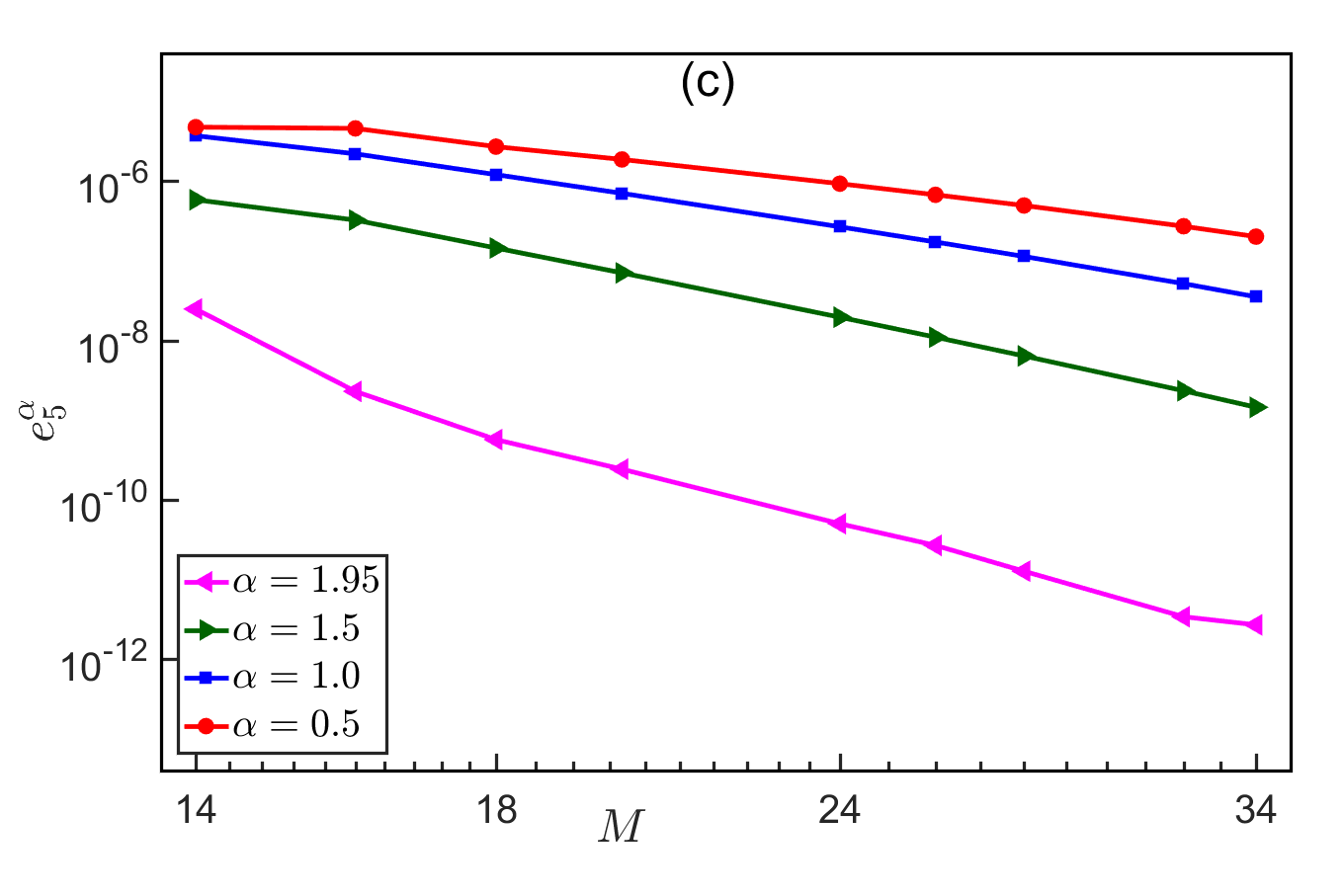,height=5cm,width=6cm,angle=0}\qquad
\psfig{figure=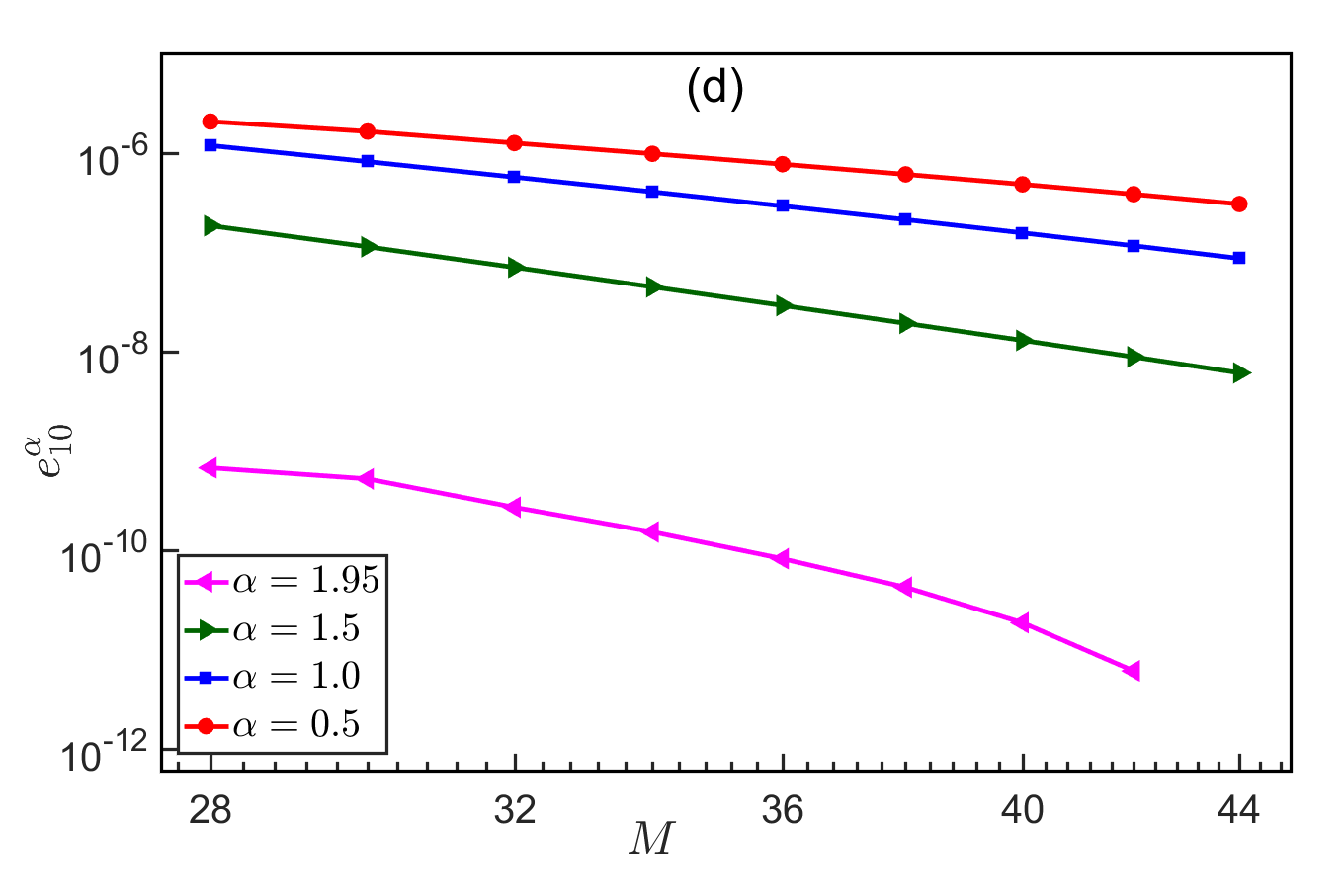,height=5cm,width=6cm,angle=0}}
\caption{Convergence rates of computing different eigenvalues of \eqref{fproblem} with $\Omega=(-1,1)$, $V(x)=\frac{x^2}{2}$ and  different $\alpha$
by using our JSM \eqref{weakn} for:
(a) the first eigenvalue $\lambda_1^\alpha$, (b) the second eigenvalue $\lambda_2^\alpha$, (c) the fifth eigenvalue $\lambda_5^\alpha$,
and (d) the tenth eigenvalue $\lambda_{10}^\alpha$.}
\label{fig:rate1}
\end{figure}

From Tabs. \ref{errors1}~\&~\ref{errors2} and Figs. \ref{fig:rate} \& \ref{fig:rate1} and extensive additional results
not shown here for brevity, we can draw the following conclusions:
(i) For fixed DOF $M$ and $\alpha\in(0,2]$, the errors from our
JSM  \eqref{weakn} are significantly smaller than those from the
FEM \cite{JLPR13} and FDM \cite{ZSP14,CD12} (cf. Tabs. \ref{errors1} \& \ref{errors2}). (ii) Both the
FEM \cite{JLPR13} and FDM \cite{ZSP14,CD12} converge almost quadratically and linearly with
respect to the DOF $M$ when $\alpha=2$ and $0<\alpha<2$, respectively (cf. Tabs. \ref{errors1}~\&~\ref{errors2}).
(iii) Our JSM method \eqref{weakn} converges spectrally and super-linearly
(or sub-spectrally)
with respect to  the DOF $M$ when $\alpha=2$ and $0<\alpha<2$, respectively (cf. Fig. \ref{fig:rate}~\&~\ref{fig:rate1}).
(iv) In Tab. \ref{eig}, the numerical results reported by our JSM \eqref{weakn} have at least eight significant digits when the DOF $M\ge160$,
while the results by the asymptotic method in \cite{ZRK07} have
at most four significant digits even when the DOF $M=5000$!
Thus our JSM method \eqref{weakn} is significantly accurate than those
low-order numerical methods in the literatures for computing eigenvalues of  the eigenvalue problem \eqref{fproblem}.

\subsection{Resolution capacity (or trust region) test}

In order to get reliable gaps and their distribution statistics, we have to calculate accurately and efficiently a very large number of eigenvalues, e.g. up to thousands or even millions eigenvalues. Specifically we need to make sure  that the numerical errors are much smaller than the minimum gap of those gaps which are used to find numerically the distribution statistics.
In general, to solve the eigenvalue problem \eqref{fproblem} by a numerical method with a given DOF $M$, we can obtain $M$ approximate eigenvalues.
A key question is that how many eigenvalues or what fraction among the $M$
approximate eigenvalues can be used to find numerically the distribution statistics, i.e. the errors to them are quite small. We remark here that for the Schr\"{o}dinger operator, i.e. $\alpha=2$ in  \eqref{fproblem}, by using a spectral method, it is proved that about $\frac{2}{\pi}$ fraction of the
$M$ approximate eigenvalues is quite accurate (or the errors are quite small) \cite{WT88}. To see whether this property is still valid for our JSM \eqref{weakn} for the FSO \eqref{fproblem},  Figure \ref{fig:error} displays the relative errors
$e_{n,r}^\alpha$ ($n=1,2,\ldots,6400$) of \eqref{fproblem}
with $V(x)\equiv0$ and different $\alpha$ by using our JSM \eqref{weakn},
FEM \cite{JLPR13} and FDM \cite{ZSP14,CD12} under the DOF $M=8192$.

\begin{figure}[h!]
\centerline{
\psfig{figure=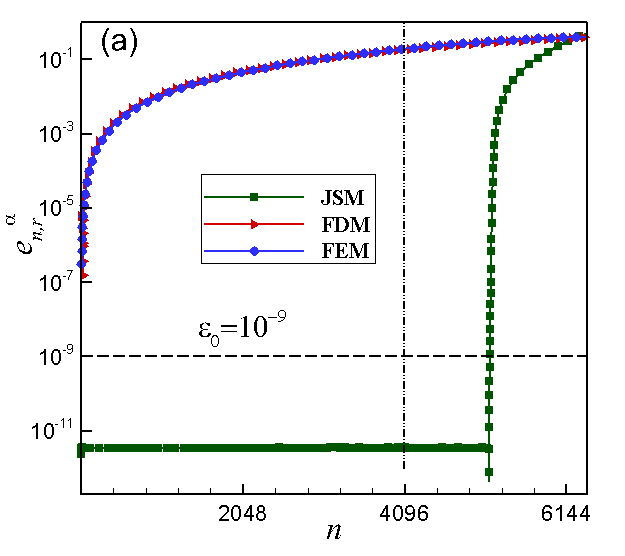,height=5cm,width=8cm,angle=0}
\psfig{figure=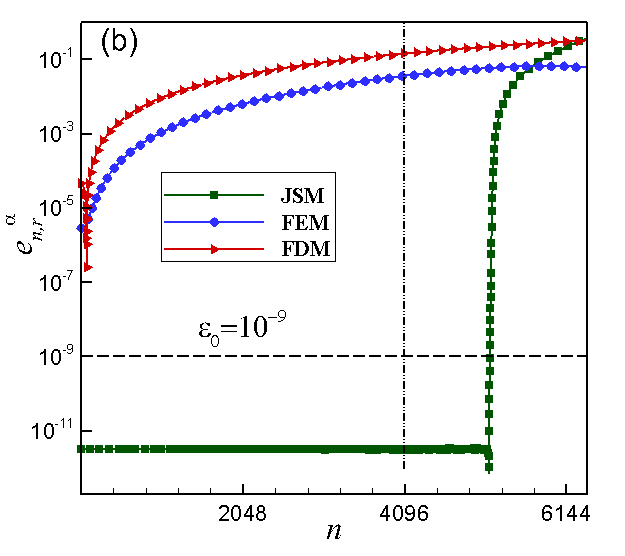,height=5cm,width=8cm,angle=0}}
\centerline{
\psfig{figure=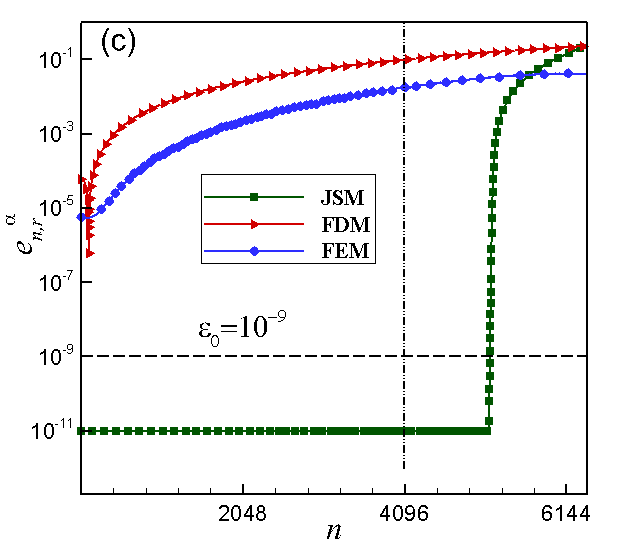,height=5cm,width=8cm,angle=0}
\psfig{figure=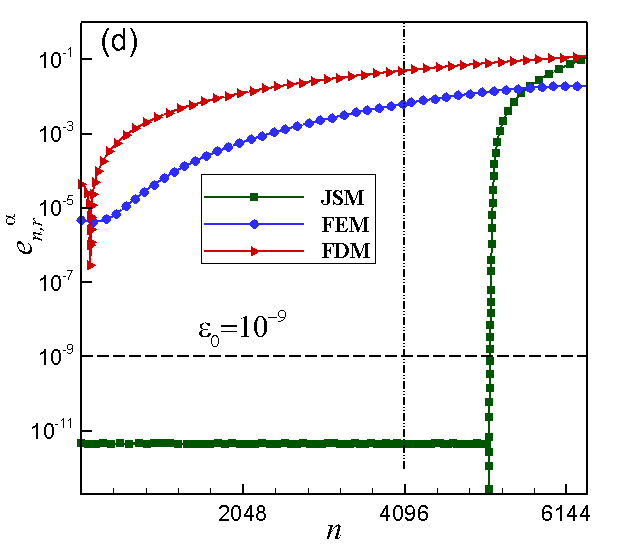,height=5cm,width=8cm,angle=0}}
\caption{Relative errors of the first $6400$ eigenvalue, i.e. $e_{n,r}^\alpha$ ($n=1,2,\ldots,6400$) of \eqref{fproblem} with $\Omega=(-1,1)$ and $V(x)\equiv0$ by using our JSM \eqref{weakn},
FEM \cite{JLPR13} and FDM \cite{ZSP14,CD12} under the DOF $M=8192$
for: (a) $\alpha=1.95$, (b) $\alpha=1.5$, (c) $\alpha=1.0$, and
(d) $\alpha=0.5$. A horizonal (dash) line with $\varepsilon_0:=10^{-9}$ and a
vertical (dash) line with $n:=M/2$ are added in each sub-figure. }
\label{fig:error}
\end{figure}

From Fig. \ref{fig:error}, we can see that our JSM \eqref{weakn}
is significantly better than FEM and FDM when a large number of eigenvalues
are to be computed accurately. In fact, FEM and FDM can be used
to compute the first a few eigenvalues of \eqref{fproblem}.
However, when a large amount of eigenvalues are needed,
one has to adapt a spectral method such as our JSM \eqref{weakn}.

To quantify the resolution capacity of our JSM \eqref{weakn},
Figure \ref{fig:errorall} displays the relative errors
$e_{n,r}^\alpha$ ($n=1,2,\ldots,M$) of \eqref{fproblem}
with $V(x)\equiv0$ and different $\alpha$ under different
DOFs $M$, i.e. $M=512$, $2048$ and $8192$; and
Figure \ref{fig:errorall1} shows similar results
when $V(x)=\frac{x^2}{2}$.

\begin{figure}[h!]
\centerline{
\psfig{figure=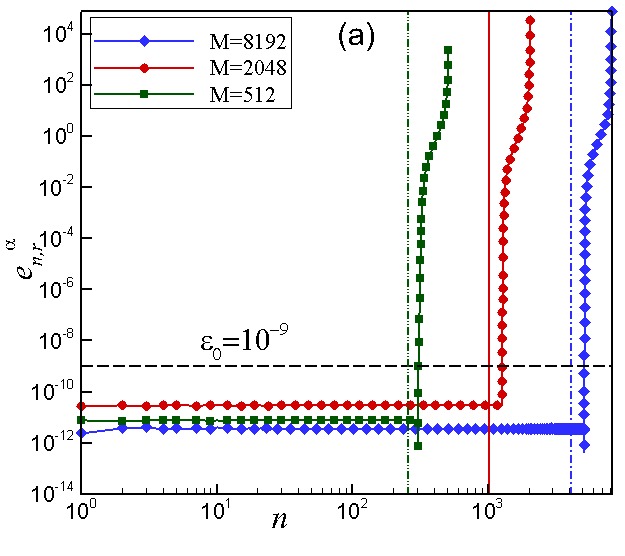,height=5cm,width=8cm,angle=0}
\psfig{figure=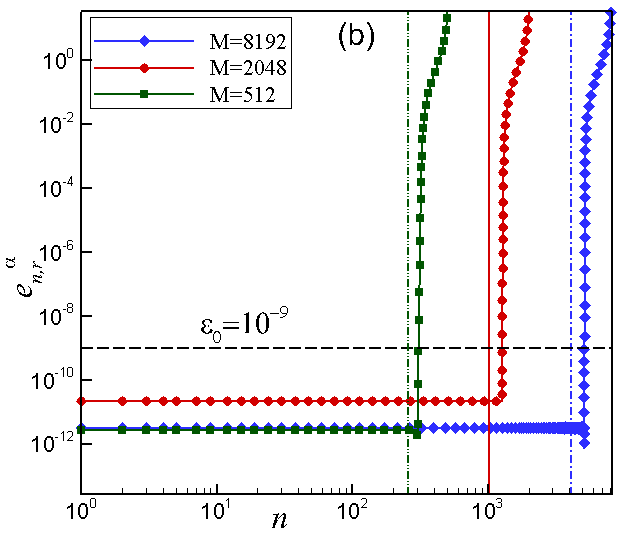,height=5cm,width=8cm,angle=0}}
\centerline{
\psfig{figure=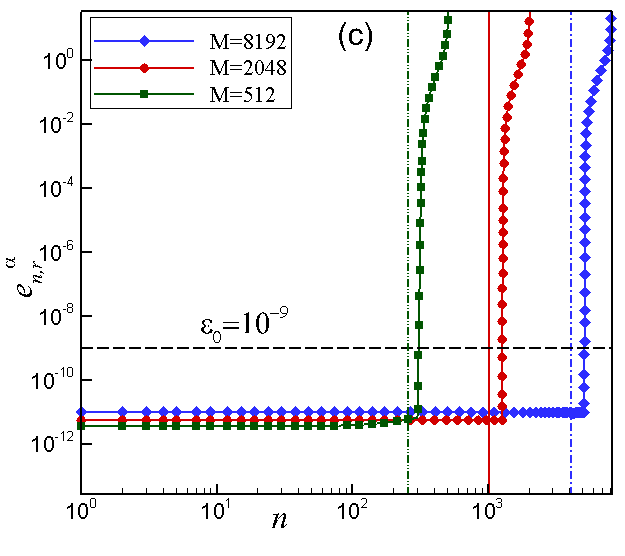,height=5cm,width=8cm,angle=0}
\psfig{figure=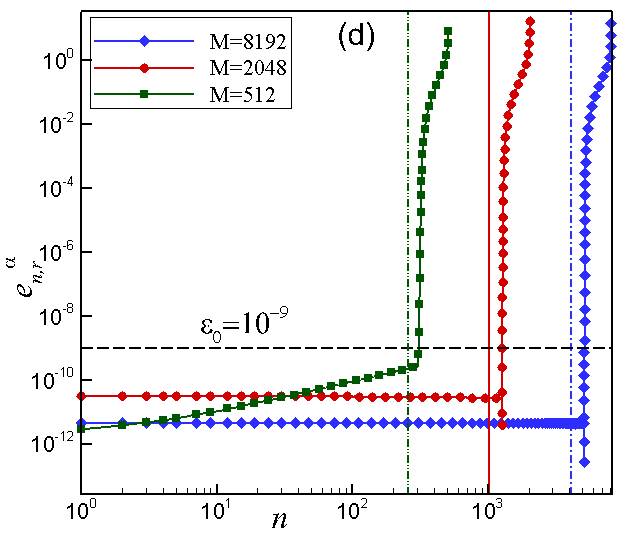,height=5cm,width=8cm,angle=0}}
\caption{Relative errors of the eigenvalues of \eqref{fproblem} with $\Omega=(-1,1)$ and $V(x)\equiv0$ by using our JSM \eqref{weakn} under different DOFs
$M$ for: (a) $\alpha=1.95$, (b) $\alpha=1.5$, (c) $\alpha=1.0$, and
(d) $\alpha=0.5$. A horizonal (dash) line with $\varepsilon_0:=10^{-9}$ and
vertical (dash) lines with $n:=M/2$ are added in each sub-figure. }
\label{fig:errorall}
\end{figure}

\begin{figure}[h!]
\centerline{
\psfig{figure=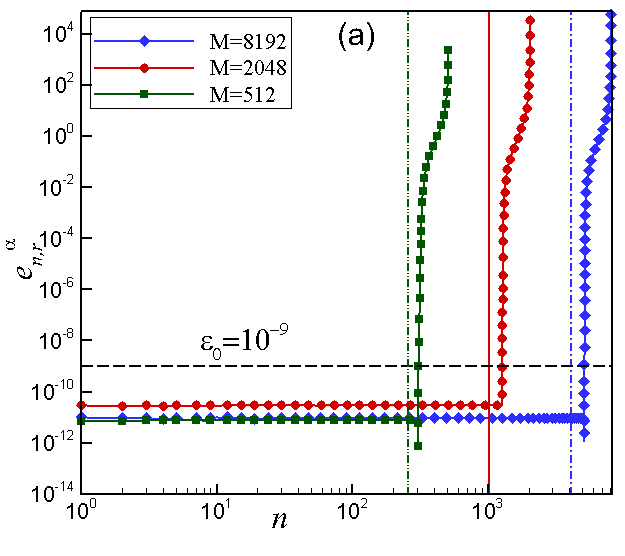,height=5cm,width=8cm,angle=0}
\psfig{figure=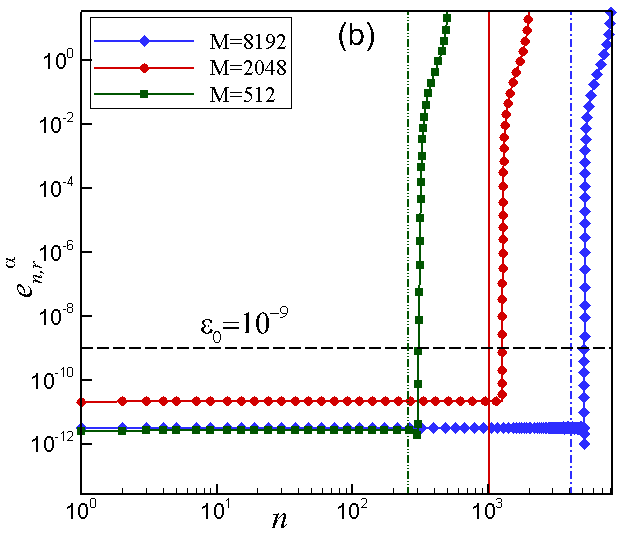,height=5cm,width=8cm,angle=0}}
\centerline{
\psfig{figure=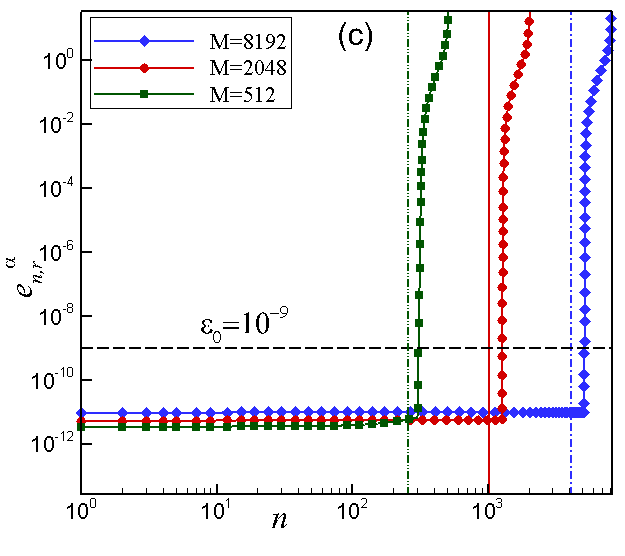,height=5cm,width=8cm,angle=0}
\psfig{figure=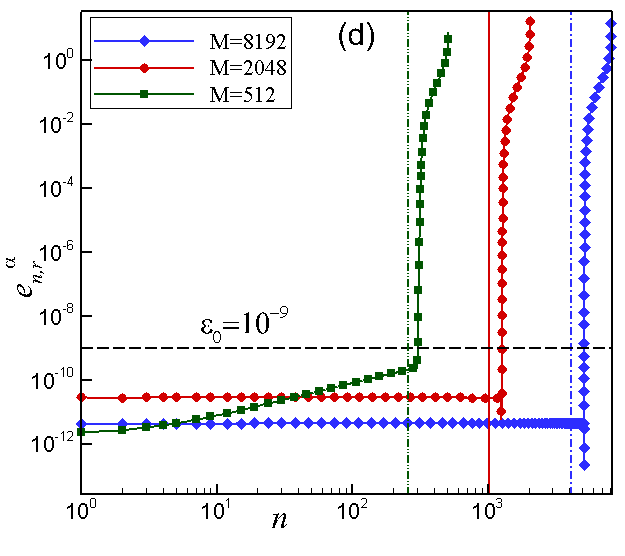,height=5cm,width=8cm,angle=0}}
\caption{Relative errors of the eigenvalues of \eqref{fproblem} with $\Omega=(-1,1)$ and $V(x)=\frac{x^2}{2}$ by using our JSM \eqref{weakn}
under different DOFs $M$ for: (a) $\alpha=1.95$, (b) $\alpha=1.5$, (c) $\alpha=1.0$, and (d) $\alpha=0.5$. A horizonal (dash) line with $\varepsilon_0:=10^{-9}$ and
vertical (dash) lines with $n:=M/2$ are added in each sub-figure.  }
\label{fig:errorall1}
\end{figure}

From  Figs. \ref{fig:errorall} \& \ref{fig:errorall1}, we can see that our JSM \eqref{weakn} under a given DOF $M$ has the following resolution capacity
(or trust region)
\be
e_{n,r}^\alpha:=\frac{\left|\lambda_n^\alpha-\lambda_{n,M}^\alpha\right|}
{\lambda_n^\alpha}\le \varepsilon_0:=10^{-9},
\qquad n=1,2,\ldots, c_r M, \qquad \hbox{with}\quad  c_r\approx \frac{2}{\pi}> \frac{1}{2}.
\ee

\section{Numerical results of FSO in 1D without potential}\label{sec:estimate}
\setcounter{equation}{0}

In this section, we report numerical results on eigenvalues of \eqref{fproblem} with $\Omega=(-1,1)$ and $V(x)\equiv0$ by using our JSM \eqref{weakn} under the DOF $M=8192$. All results are based on the first
$4096$ eigenvalues, i.e. we use half of the eigenvalues obtained numerically
to present the results and to calculate distribution statistics.

\subsection{Eigenvalues and their approximations}

Figure \ref{fig:eig1dnp}a plots eigenvalues $\lambda_n^\alpha$ ($n=1,2,\ldots$) and their leading order approximations as
$\lambda_n^\alpha\approx
\tilde \lambda_n^{\alpha}:=\left(\frac{n\pi}{2}\right)^\alpha$ ($n=1,2,\ldots$), while $\tilde \lambda_n^{\alpha}$ ($n=1,2,\ldots$) are
the eigenvalues of the {\sl local fractional Laplacian
operator} on $\Omega=(-1,1)$ with homogeneous Dirichlet boundary condition \cite{BR18}.
Figure \ref{fig:eig1dnp}b displays the relative errors of the  eigenvalues and their
leading order approximations, i.e. $\tilde e_{n,r}^\alpha:=\left(\tilde \lambda_n^{\alpha}-\lambda_n^\alpha\right)
/\tilde \lambda_n^{\alpha}$, which immediately
suggests a high order approximation at
$\lambda_n^\alpha\approx
\hat \lambda_n^{\alpha}:=\tilde \lambda_n^{\alpha}
\left(1-\frac{C_3^\alpha}{n}\right) $ ($n=1,2,\ldots$). By fitting our numerical results, we can obtain numerically $C_3^{\alpha}=\frac{\alpha(2-\alpha)}{4}$ which is plotted in
Figure \ref{fig:eig1dnp}c. Finally Figure \ref{fig:eig1dnp}d displays the absolute errors
of the  eigenvalues and their
high order approximations, i.e. $\tilde e_n^\alpha:=\left|\lambda_n^\alpha-
\hat \lambda_n^{\alpha}\right|$.

\begin{figure}[h!]
\centerline{
\psfig{figure=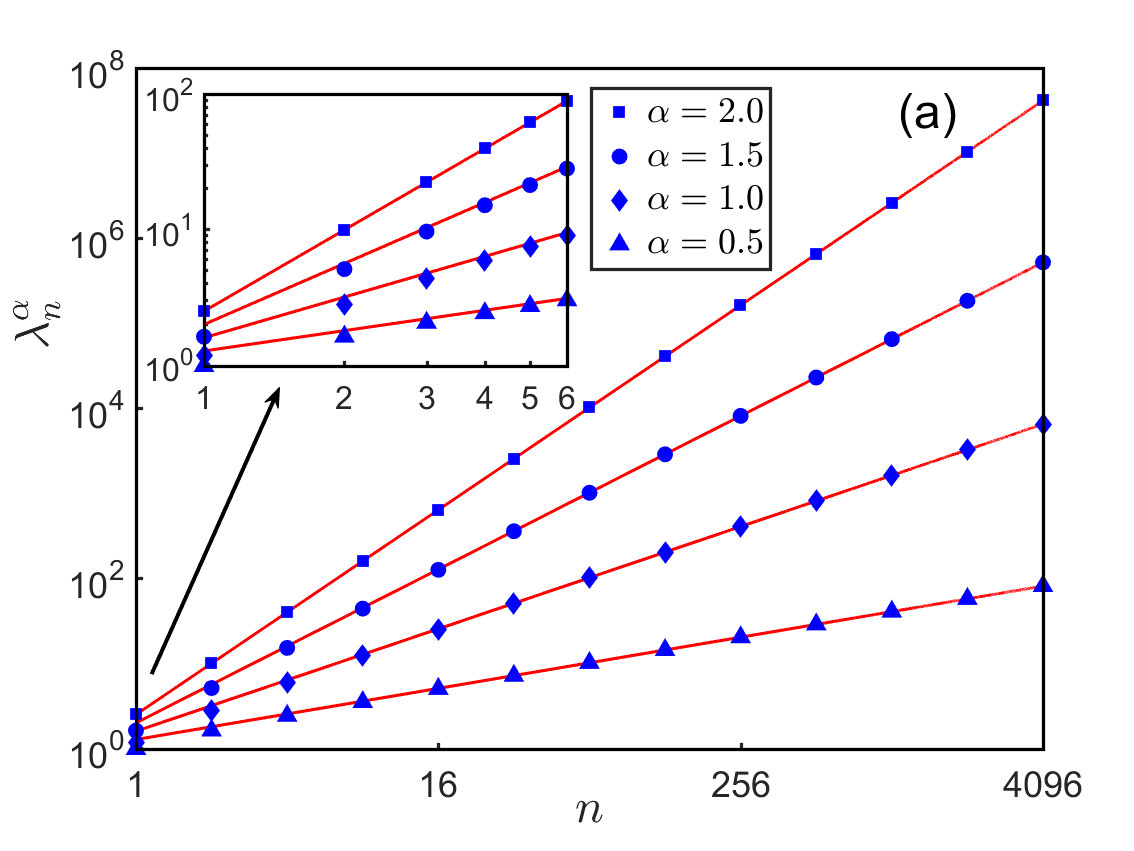,height=6cm,width=7cm,angle=0}
\psfig{figure=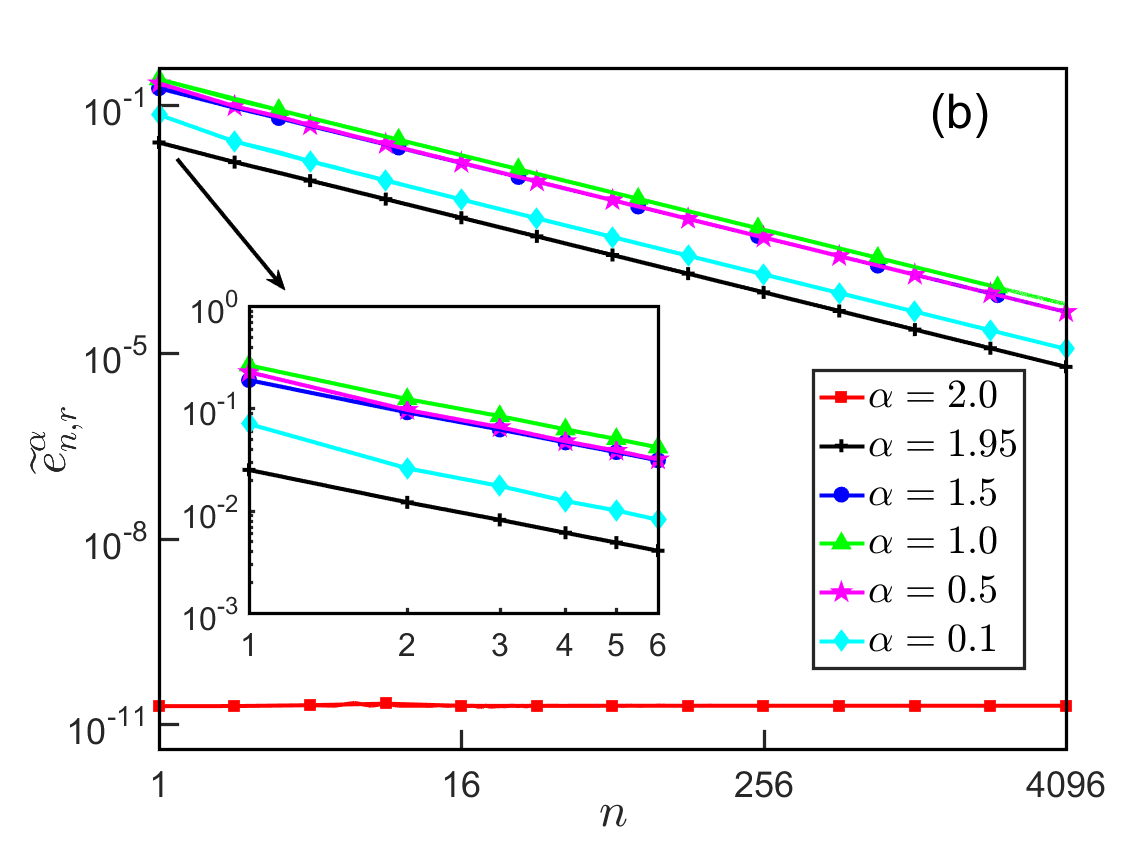,height=6cm,width=7cm,angle=0}}
\centerline{\psfig{figure=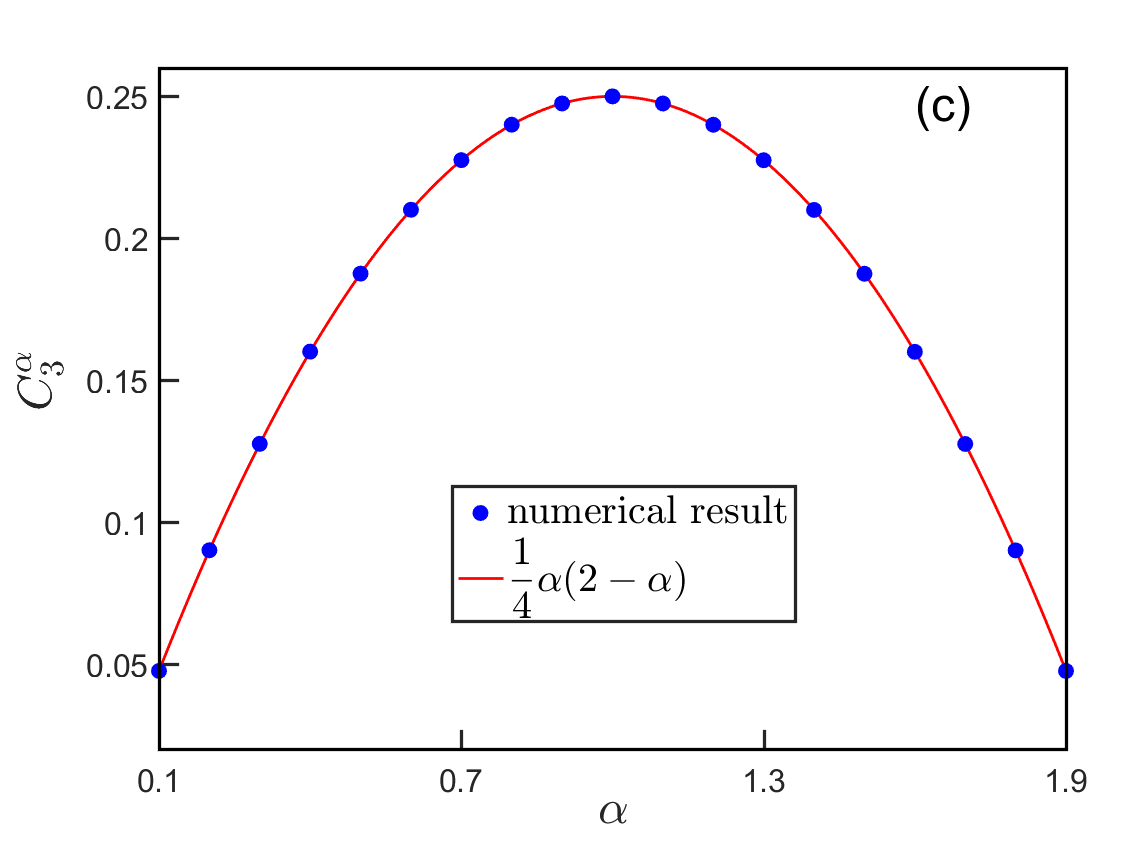,height=6cm,width=7cm,angle=0}
\psfig{figure=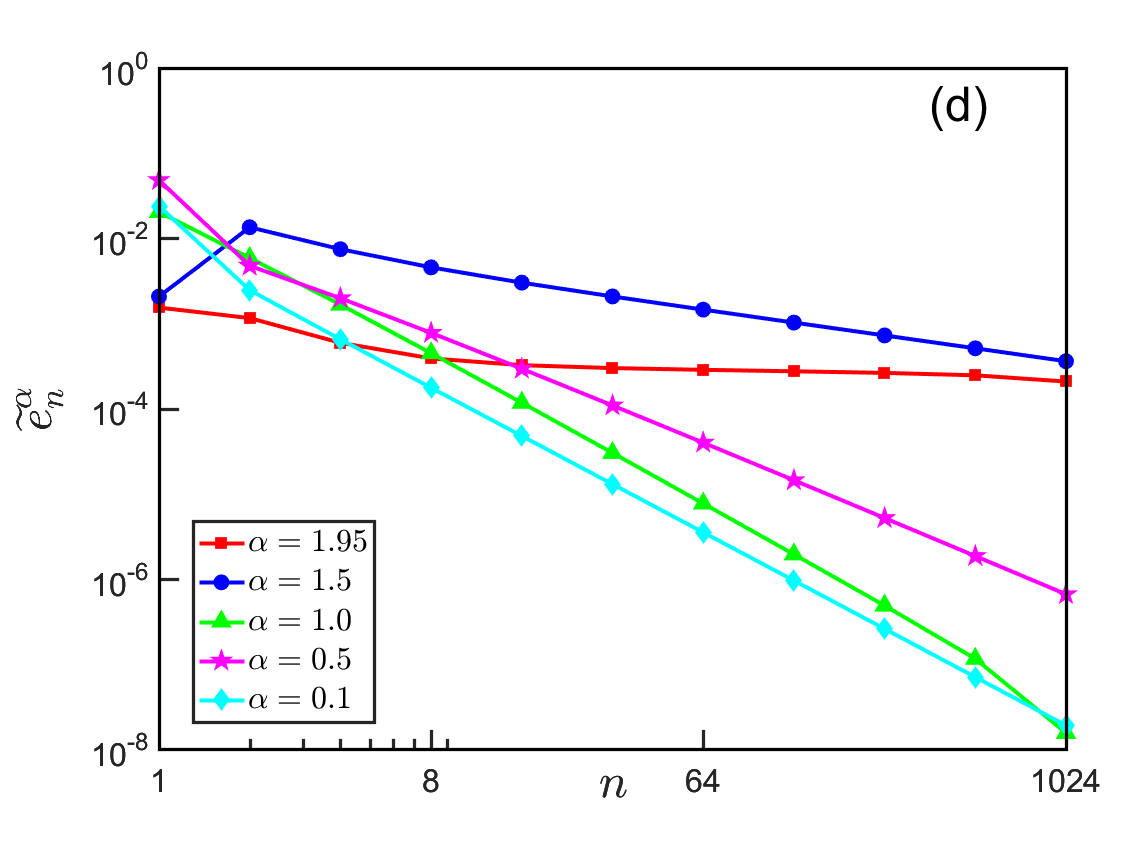,height=6cm,width=7cm,angle=0}}
\caption{(a) Eigenvalues $\lambda_n^\alpha$ ($n=1,2,\ldots,4096$) of \eqref{fproblem} with $\Omega=(-1,1)$ and $V(x)\equiv0$ for different
$\alpha$ (symbols denote numerical results and solids lines are
from the leading order approximation $\tilde \lambda_n^{\alpha }=\left(\frac{n\pi}{2}\right)^\alpha$);
(b) Relative errors $\tilde e_{n,r}^\alpha=\left(\tilde \lambda_n^{\alpha}-\lambda_n^\alpha\right)
/\tilde \lambda_n^{\alpha}$ (symbols denote numerical results and solids lines are from fitting formula $C_3^\alpha n^{-1}$ when $n\gg1$); (c) Fitting results for $C_3^\alpha$; and (d) absolute errors $\tilde e_n^\alpha=\left|\lambda_n^\alpha-
\hat \lambda_n^{\alpha}\right|$ with $\hat \lambda_n^{\alpha}=\tilde \lambda_n^{\alpha}
\left(1-C_3^\alpha\,n^{-1}\right)$. }
\label{fig:eig1dnp}
\end{figure}

From Fig. \ref{fig:eig1dnp}, we can obtain numerically the following
approximations of the eigenvalues of \eqref{fproblem} with $\Omega=(-1,1)$ and $V(x)\equiv0$ as
\be \label{asy-eig-1Dd11}
\lambda_n^\alpha= \hat \lambda_n^{\alpha}+O(n^{\alpha-2})=\tilde \lambda_n^{\alpha} \left[1-\frac{\alpha(2-\alpha)}{4n}+O(n^{-2})\right], \quad n=1,2,\ldots\;,
\ee
where
\be \label{eigapprx}
\tilde \lambda_n^{\alpha}=\left(\frac{n\pi}{2}\right)^\alpha, \quad
\hat \lambda_n^{\alpha}=\left(\frac{n\pi}{2}\right)^\alpha-
\left(\frac{\pi}{2}\right)^\alpha\frac{\alpha(2-\alpha)}{4} n^{\alpha-1}=\tilde \lambda_n^{\alpha}
\left[1-\frac{\alpha(2-\alpha)}{4n}\right],
\quad n\ge1, \quad 0<\alpha\le 2.\quad
\ee
Combining \eqref{asy-eig-1Dd11} and Lemma 2.1, we can immediately obtain the conclusion \eqref{asy-eig-1D}.

  To demonstrate high accuracy of our numerical method, Table \ref{eigs01}
lists eigenvalues of \eqref{fproblem} with $\Omega=(-1,1)$ and $V(x)\equiv0$  for different $\alpha$.

  \begin{table}[h!]
\centering
\begin{tabular}{|c |c |c|c|c|c|c|} \hline
& $\alpha=0.1$&$\alpha=0.5$&$\alpha=1.0$&$\alpha=1.5$& $\alpha=1.95$ & $\alpha=2.0$ \\ \hline
$\lambda_1^\alpha$    &0.9725944  &0.9701654 &1.157773883 &1.5975035456 &2.35198053244 &2.4674011002 \\
$\lambda_2^\alpha$    &1.0921964  &1.6015377 &2.754754742 &5.0597599283 &9.20812426623 &9.8696044010 \\
$\lambda_3^\alpha$    &1.1473224  &2.0288210 &4.316801066 &9.5943057675 &20.3833201062 &22.206609902 \\
$\lambda_4^\alpha$    &1.1868395  &2.3871563 &5.892147470 &15.018786212 &35.7934316323 &39.478417604 \\
$\lambda_5^\alpha$    &1.2165513  &2.6947426 &7.460175739 &21.189425897 &55.3737634238 &61.685027506 \\
$\lambda_6^\alpha$    &1.2412799  &2.9728959 &9.032852690 &28.035207791 &79.0793754673 &88.826439609 \\
$\lambda_7^\alpha$    &1.2619743  &3.2256090 &10.60229309 &35.488011031 &106.871259423 &120.90265391 \\
$\lambda_8^\alpha$    &1.2801923  &3.4610502 &12.17411826 &43.507108689 &138.718756729 &157.91367041 \\
$\lambda_9^\alpha$    &1.2961956  &3.6805940 &13.74410905 &52.051027490 &174.594065184 &199.85948912 \\
$\lambda_{10}^\alpha$ &1.3107082  &3.8884472 &15.31555499 &61.092457389 &214.473975149 &246.74011002 \\
$\lambda_{20}^\alpha$ &1.4082270  &5.5522311 &31.02330309 &174.43784577 &829.684155066 &986.96044010 \\
$\lambda_{40}^\alpha$ &1.5111219  &7.8894197 &62.43917339 &495.95713648 &3207.64320222 &3947.8417604 \\
$\lambda_{60}^\alpha$ &1.5742803  &9.6777480 &93.85508924 &912.11187382 &7073.79138904 &8882.6439609 \\
\hline
\end{tabular}
\caption{Eigenvalues of \eqref{fproblem} with $\Omega=(-1,1)$ and $V(x)\equiv0$  for different $\alpha$.}
\label{eigs01}
\end{table}

\subsection{Asymptotic behaviour of different gaps}

\begin{figure}[h!]
\centerline{
\psfig{figure=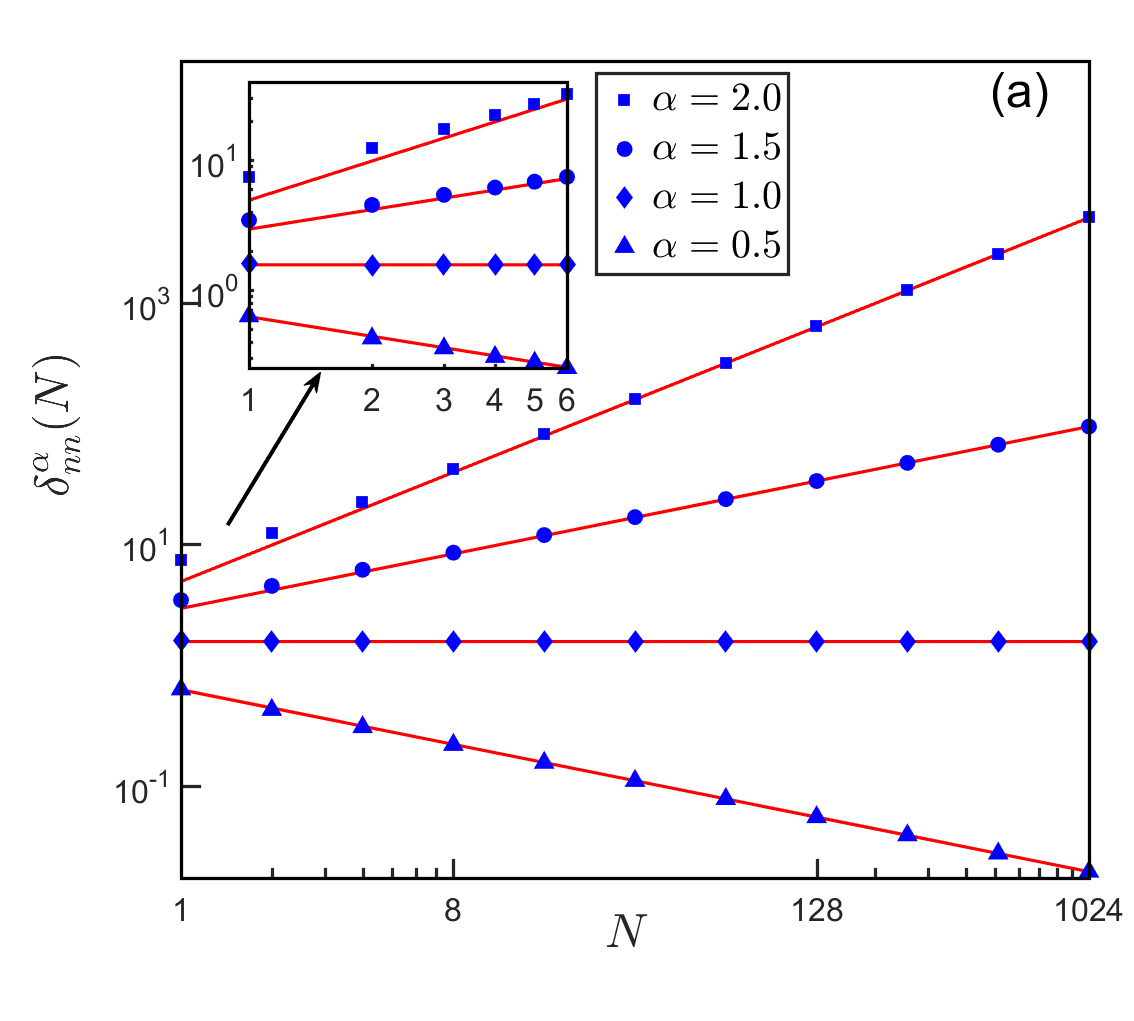,height=5cm,width=7cm,angle=0}
\psfig{figure=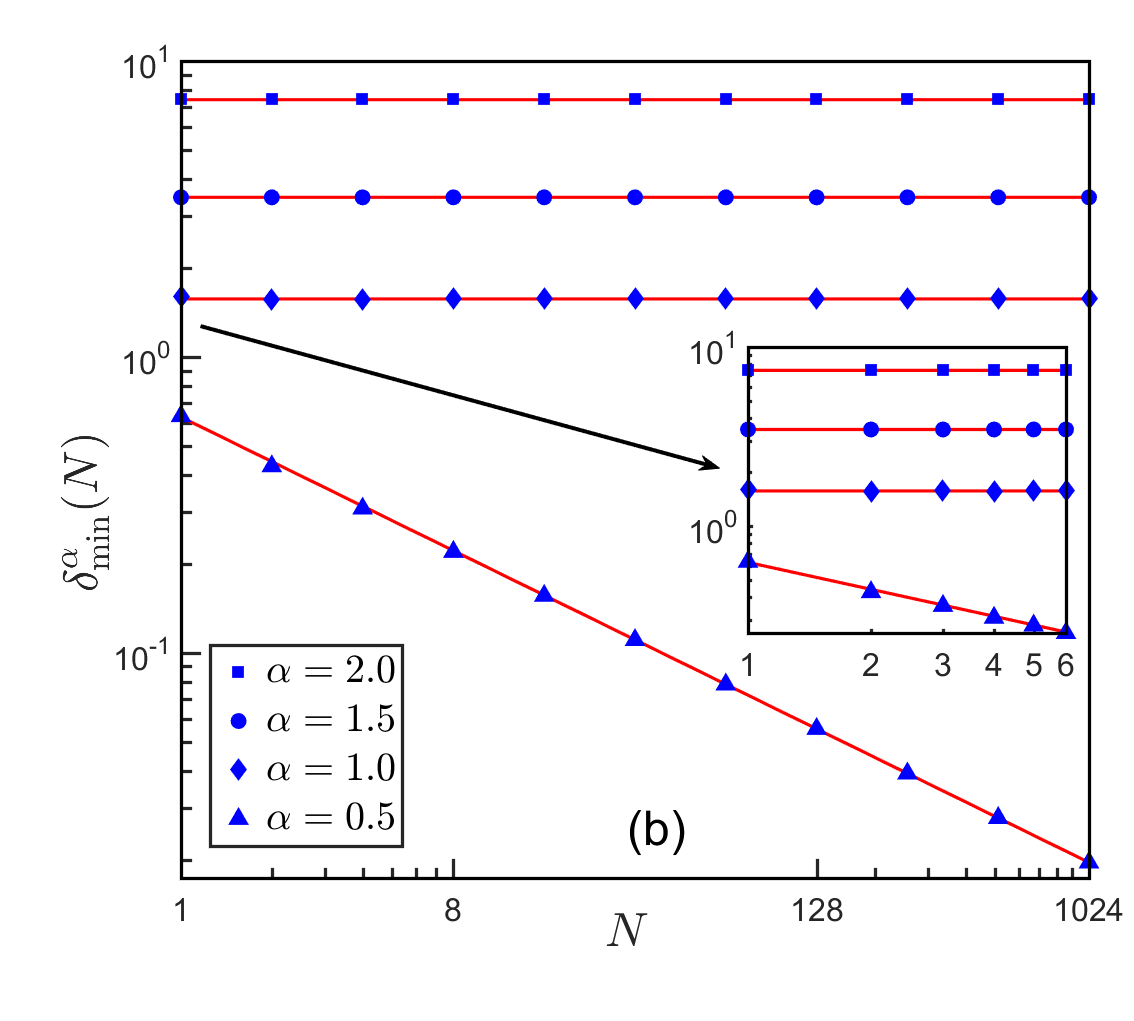,height=5cm,width=7cm,angle=0}}
\centerline{\psfig{figure=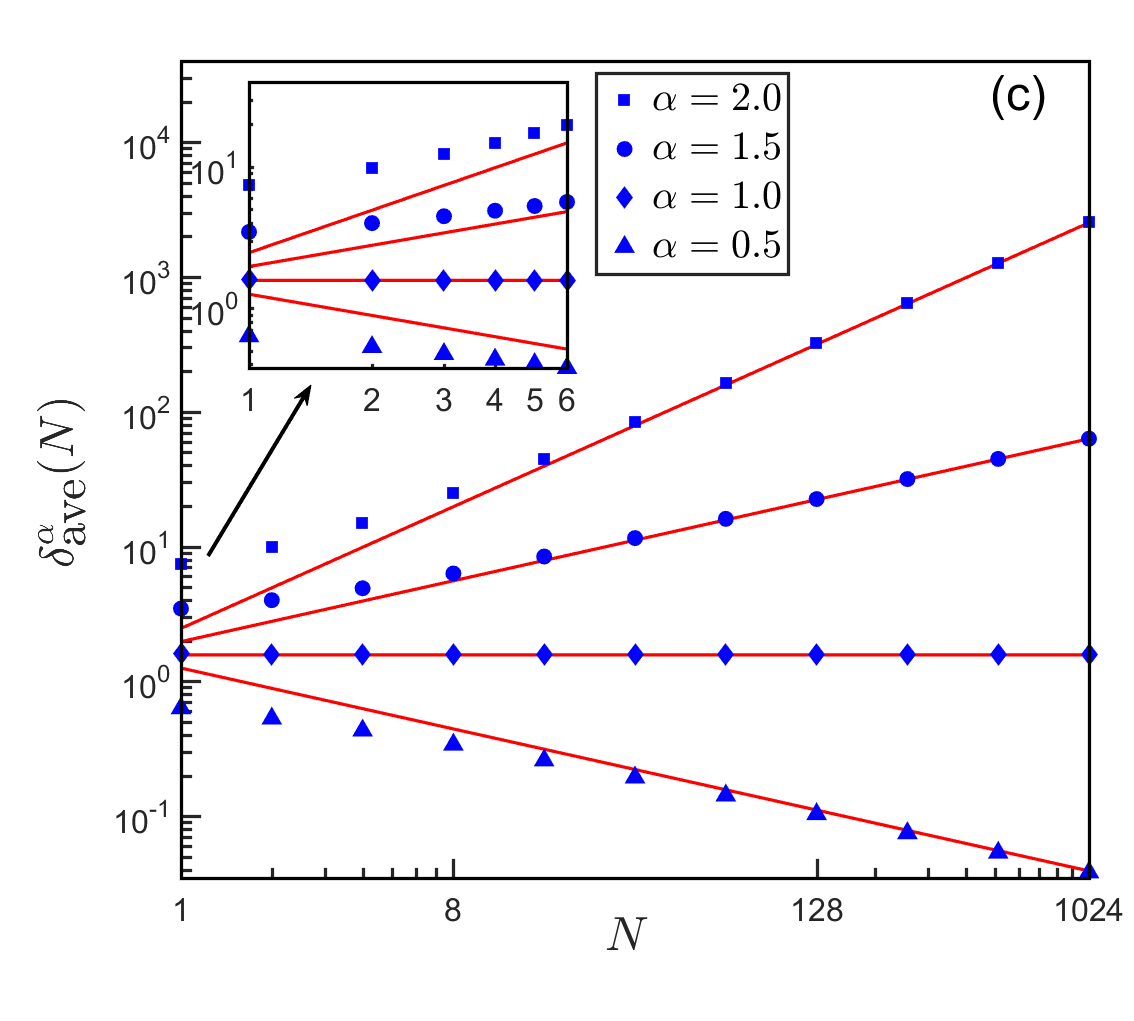,height=5cm,width=7cm,angle=0}
\psfig{figure=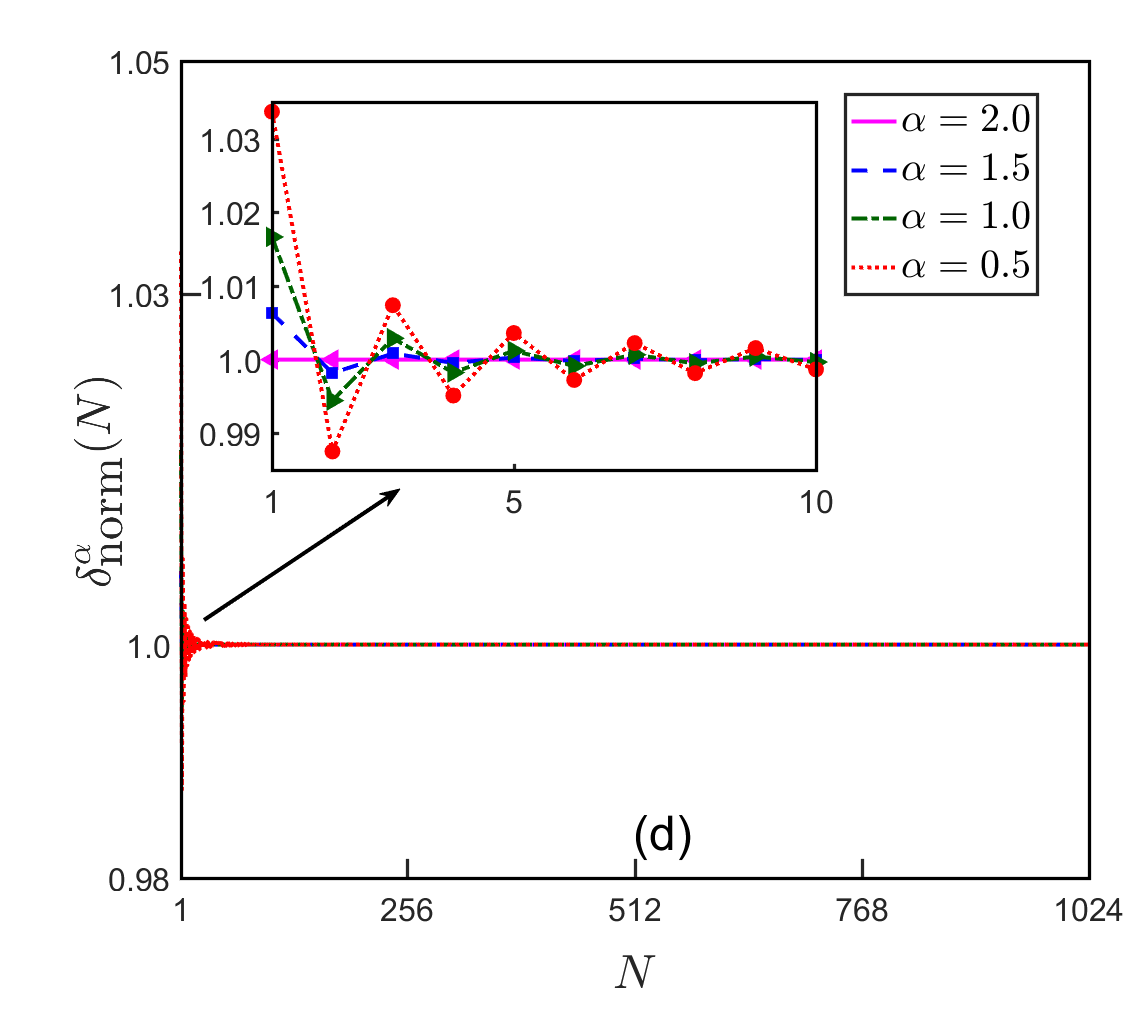,height=5cm,width=7cm,angle=0}}
\caption{Different eigenvalue gaps of \eqref{fproblem}
with $\Omega=(-1,1)$, $V(x)\equiv0$ and different $\alpha$
for (symbols denote numerical results and solids lines are from fitting formula when $N\gg1$): (a) the nearest neighbour gaps $\delta_{\rm nn}^\alpha(N)$, (b)
the minimum gaps $\delta_{\textrm{min}}^\alpha(N)$,
 (c) the average gaps  $\delta_{\textrm{ave}}^\alpha(N)$,
 and (d) the normalized gaps
$\delta_{\rm norm}^\alpha(N)$. }
\label{fig:gaps}
\end{figure}

Figure \ref{fig:gaps} plots different eigenvalue gaps of \eqref{fproblem}
with $\Omega=(-1,1)$, $V(x)\equiv0$ and different $\alpha$.
From Fig. \ref{fig:gaps}, we can draw the following conclusions based on our numerical results: (i) the nearest neighbour gaps $\delta_{\rm nn}^\alpha(N)$
increase and decrease with respect to $N$ when $1<\alpha \le 2$ and $0<\alpha<1$, respectively;
and they are almost constant when $\alpha=1$ (cf. Fig. \ref{fig:gaps}a). (ii) The minimum gaps $\delta_{\textrm{min}}^\alpha(N)$ are almost constants
and decrease with respect to $N$ when $1\le \alpha\le 2$ and $0<\alpha<1$, respectively (cf. Fig. \ref{fig:gaps}b). (iii)
The average gaps  $\delta_{\textrm{ave}}^\alpha(N)$ increase and decrease
with respect to $N$ when $1<\alpha \le 2$ and $0<\alpha<1$, respectively;
and they are almost constant when $\alpha=1$ (cf. Fig. \ref{fig:gaps}c).
(iv) The normalized  gaps
$\delta_{\rm norm}^\alpha(N)\approx 1$
when $N\gg1$  (cf. Fig. \ref{fig:gaps}d).

In fact, based on the numerical asymptotic approximation \eqref{asy-eig-1Dd11},
we can formally obtain the following approximation of the nearest neighbour
gaps as
\bea \label{nngapar}
\delta_{\rm nn}^\alpha(N)&=&\lambda_{N+1}^\alpha-\lambda_{N}^\alpha
\approx \hat \lambda_{N+1}^{\alpha}-\hat \lambda_{N}^{\alpha}\nn\\
&=&\left(\frac{(N+1)\pi}{2}\right)^\alpha-
\left(\frac{\pi}{2}\right)^\alpha\frac{\alpha(2-\alpha)}{4} (N+1)^{\alpha-1}
-\left(\frac{N\pi}{2}\right)^\alpha+
\left(\frac{\pi}{2}\right)^\alpha\frac{\alpha(2-\alpha)}{4} N^{\alpha-1}
\nn\\
&=&\left(\frac{\pi}{2}\right)^\alpha\left[(N+1)^\alpha-N^\alpha-
\frac{\alpha(2-\alpha)}{4} \left((N+1)^{\alpha-1}-N^{\alpha-1}\right) \right]\nn\\
&=&\left(\frac{\pi}{2}\right)^\alpha\left[N^\alpha\left(\left(1
+\frac{1}{N}\right)^\alpha-1\right)-
\frac{\alpha(2-\alpha)}{4} N^{\alpha-1}\left(\left(1+\frac{1}{N}\right)^{\alpha-1}-1\right)\right]\nn\\
&=&\left(\frac{\pi}{2}\right)^\alpha\left[N^\alpha\left(\frac{\alpha}{N}+
\frac{\alpha(\alpha-1)}{N^2}+O(N^{-3})\right)-
\frac{\alpha(2-\alpha)}{4} N^{\alpha-1}\left(\frac{\alpha-1}{N}+O(N^{-2})\right) \right]\nn\\
&=&\left(\frac{\pi}{2}\right)^\alpha\left[\alpha N^{\alpha-1}+ \frac{\alpha(\alpha-1)(2+\alpha)}{4}N^{\alpha-2}+O(N^{\alpha-3})\right],
\qquad N=1,2,\ldots\;.
\eea
Again, this asymptotic results also confirm that the nearest neighbour gaps $\delta_{\rm nn}^\alpha(N)$
increase and decrease with respect to $N$ when $1<\alpha \le 2$ and $0<\alpha<1$, respectively;
and they are almost constant when $\alpha=1$.

Based on the asymptotic results \eqref{nngapar} and the numerical results in
Fig. \ref{fig:gaps}b, we can conclude that
\be\label{mingapar}
\delta_{\rm min}^\alpha(N)=\left\{\ba{ll}
\delta_{\rm nn}^\alpha(1)=\lambda_2^\alpha-\lambda_1^\alpha, &1< \alpha<2,\\
\approx \delta_{\rm nn}^\alpha(1)=\lambda_2^\alpha-\lambda_1^\alpha, & \alpha=1,\\
\delta_{\rm nn}^\alpha(N)\approx \alpha \left(\frac{\pi}{2}\right)^\alpha N^{\alpha-1},  &0<\alpha<1.\\
\ea  \right.
\qquad N=1,2,\ldots\;.
\ee
Again, these asymptotic results suggest that the minimum gaps $\delta_{\textrm{min}}^\alpha(N)$ are almost constants
and decrease with respect to $N$ when $1\le \alpha\le 2$ and $0<\alpha<1$, respectively.

Similarly, we have the asymptotic results for the average gaps as
\bea \label{avegapar}
\delta_{\rm ave}^\alpha(N)&=&\frac{\lambda_{N+1}^\alpha-\lambda_{1}^\alpha}{N}
\approx \frac{\hat \lambda_{N+1}^{\alpha}-\lambda_{1}^{\alpha}}{N}\nn\\
&=&\frac{1}{N}\left[\left(\frac{(N+1)\pi}{2}\right)^\alpha-
\left(\frac{\pi}{2}\right)^\alpha\frac{\alpha(2-\alpha)}{4} (N+1)^{\alpha-1}
-\lambda_{1}^{\alpha}\right]\nn\\
&=&\left(\frac{\pi}{2}\right)^\alpha\left[N^{\alpha-1}\left(1+\frac{1}{N}
\right)^\alpha-\frac{\alpha(2-\alpha)}{4}N^{\alpha-2}\left(1+\frac{1}{N}
\right)^{\alpha-1}-\lambda_{1}^{\alpha}
\left(\frac{2}{\pi}\right)^{\alpha}N^{-1}\right]\nn\\
&=&\left(\frac{\pi}{2}\right)^\alpha\left[N^{\alpha-1}+\alpha N^{\alpha-2}-
\frac{\alpha(2-\alpha)}{4}N^{\alpha-2}-\lambda_{1}^{\alpha}
\left(\frac{2}{\pi}\right)^{\alpha}N^{-1}+O(N^{\alpha-3})\right]\nn\\
&=&\left(\frac{\pi}{2}\right)^\alpha\left[N^{\alpha-1}+
\frac{\alpha(2+\alpha)}{4}N^{\alpha-2}-\lambda_{1}^{\alpha}
\left(\frac{2}{\pi}\right)^{\alpha}N^{-1}+O(N^{\alpha-3})\right],
\qquad N=1,2,\ldots\;.
\eea
Thus when $1<\alpha<2$, we have
\be \label{avegapar1}
\delta_{\rm ave}^\alpha(N)=\left(\frac{\pi}{2}\right)^\alpha\left[N^{\alpha-1}+
\frac{\alpha(2+\alpha)}{4}N^{\alpha-2}+O(N^{-1})\right], \qquad N=1,2,\ldots\;,
\ee
and when $0<\alpha<1$, we have
\be \label{avegapar2}
\delta_{\rm ave}^\alpha(N)=\left(\frac{\pi}{2}\right)^\alpha\left[N^{\alpha-1}
-\lambda_{1}^{\alpha}
\left(\frac{2}{\pi}\right)^{\alpha}N^{-1}+O(N^{\alpha-2})\right],\qquad N=1,2,\ldots\;,
\ee
and when $\alpha=1$, we get
\be \label{avegapar3}
\delta_{\rm ave}^\alpha(N)=\frac{\pi}{2}\left[1+
\left(\frac{3}{4}-\frac{2}{\pi}\lambda_{1}^{\alpha=1}
\right)N^{-1}+O(N^{-2})\right],\qquad N=1,2,\ldots\;.
\ee
Again, these asymptotic results suggest that the average gaps  $\delta_{\textrm{ave}}^\alpha(N)$ increase and decrease with respect to $N$ when $1<\alpha \le 2$ and $0<\alpha<1$, respectively;
and they are almost constants when $\alpha=1$ (cf. Fig. \ref{fig:gaps}c).

Based on the asymptotic results of the eigenvalue $\lambda_n^\alpha$ in \eqref{asy-eig-1Dd11}, noticing \eqref{nngap1}-\eqref{nngap2}, we can get the asymptotic results for the normalized gaps as
\bea \label{normgapar}
\delta_{\rm norm}^\alpha(N)&=&\frac{2}{\pi}\left[\left(\lambda_{N+1}^\alpha\right)^{1/\alpha}-
\left(\lambda_{N}^\alpha\right)^{1/\alpha}\right]\nonumber\\
&=&(N+1)\left(1-\frac{\alpha(2-\alpha)}{4(N+1)}+O((N+1)^{-2})
\right)^{1/\alpha}-N\left(1-\frac{\alpha(2-\alpha)}{4N}+O(N^{-2})
\right)^{1/\alpha}\nonumber\\
&=&N+1-\frac{2-\alpha}{4}-\frac{\tilde C}{N+1}+O((N+1)^{-2})
-N+\frac{2-\alpha}{4}+\frac{\tilde C}{N}-O(N^{-2})\nonumber\\
&=&1+\frac{\tilde C}{N(N+1)}+O(N^{-3}),\qquad N=1,2,\ldots\;,
\eea
where $\tilde C$ is a constant.
Again, this asymptotic result suggests that the normalized  gaps  $\delta_{\textrm{norm}}^\alpha(N)\approx 1$
when $N\gg1$ (cf. Fig. \ref{fig:gaps}d).

Finally, combining \eqref{nngapar}, \eqref{mingapar}, \eqref{avegapar1},
 \eqref{avegapar2}, \eqref{avegapar3}, \eqref{normgapar}
 and \eqref{gapscl}, we can get the conjecture \eqref{asy-gaps-1D}
 stated in Section 1.

\subsection{The gap distribution statistics}

Figure \ref{fig:gaps1da} displays the histogram of the normalized  gaps $\{\delta_{\rm norm}^\alpha(n)\ |\ 1\le n\le N=4096\}$ defined in \eqref{normgapar} for  \eqref{fproblem}
with $\Omega=(-1,1)$, $V(x)\equiv0$ and different $\alpha$.

\begin{figure}[h!]
\centerline{
\psfig{figure=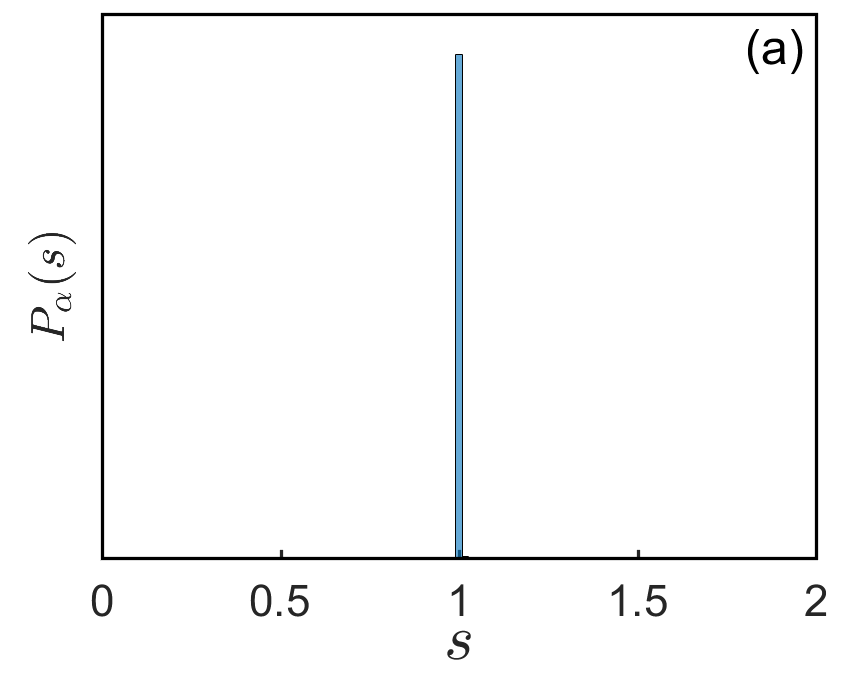,height=4cm,width=5cm,angle=0}
\psfig{figure=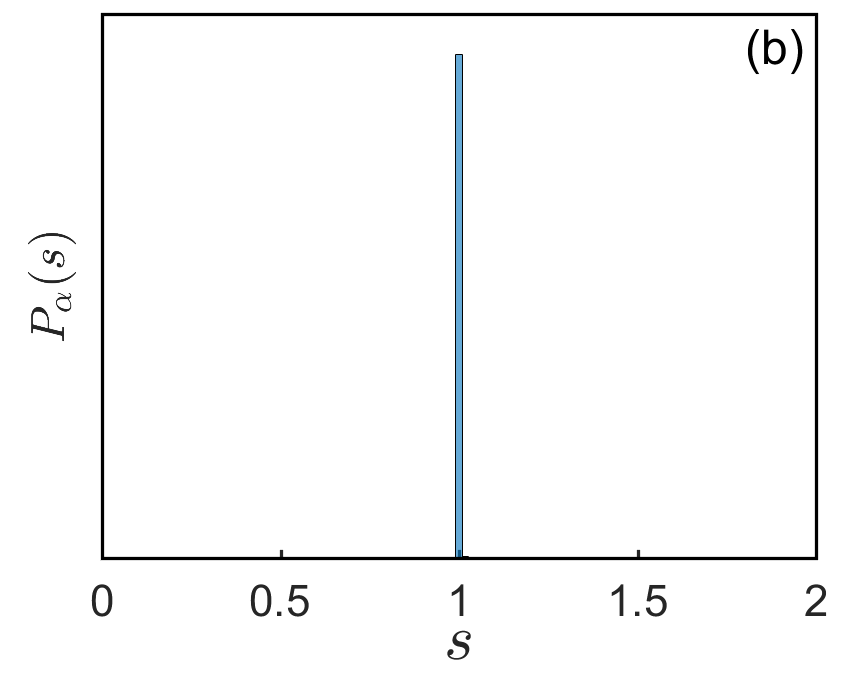,height=4cm,width=5cm,angle=0}
\psfig{figure=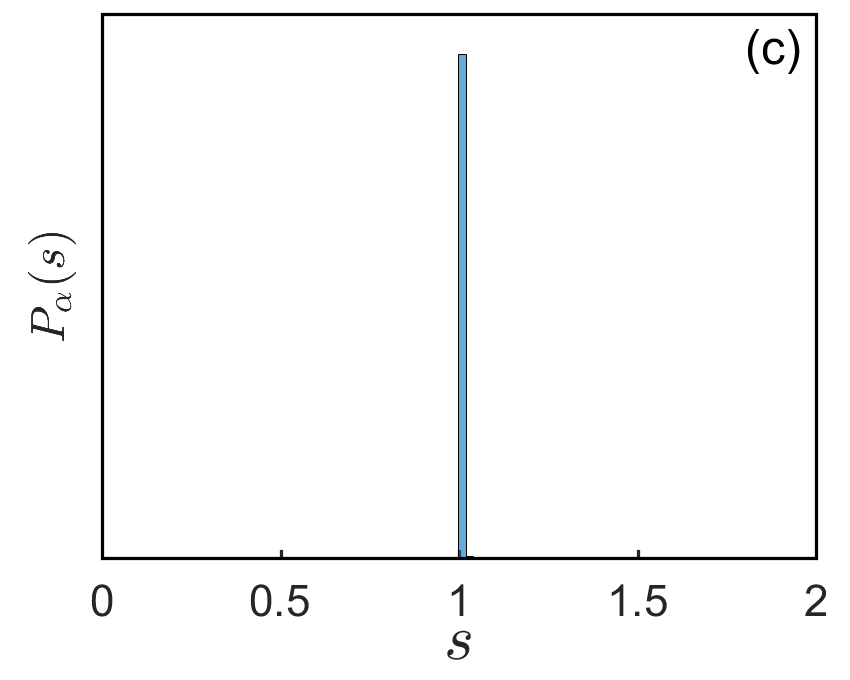,height=4cm,width=5cm,angle=0}}
\centerline{
\psfig{figure=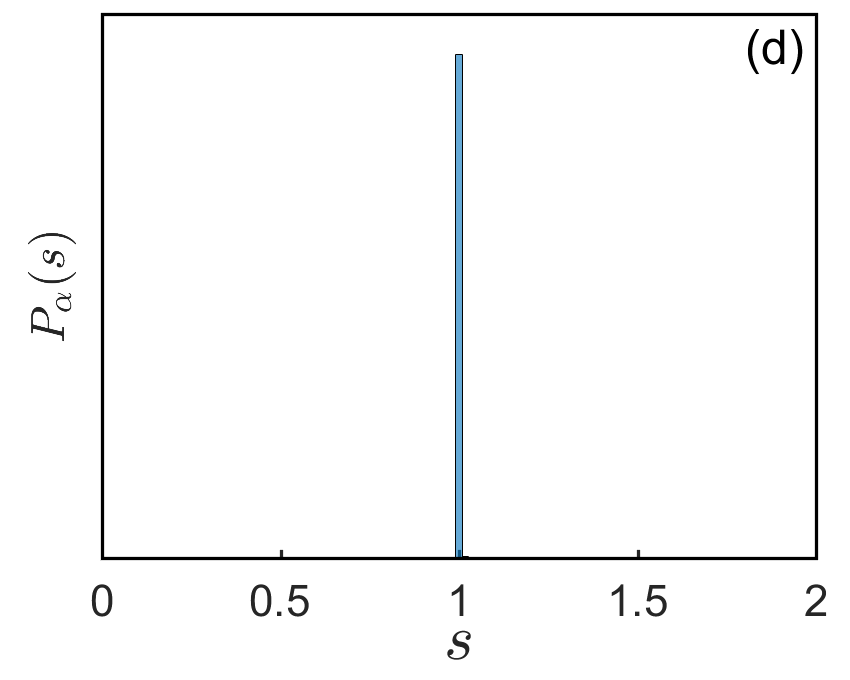,height=4cm,width=5cm,angle=0}
\psfig{figure=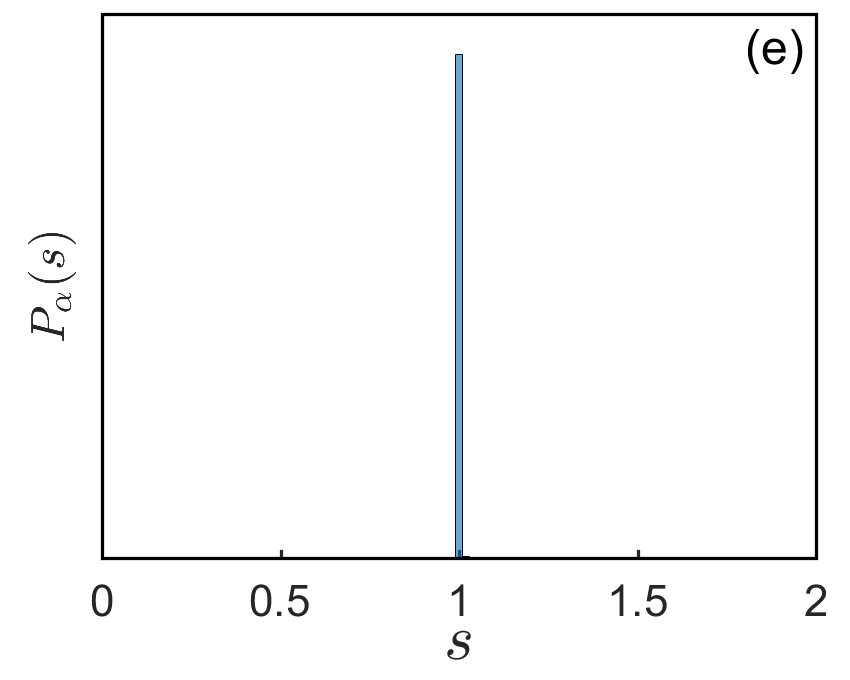,height=4cm,width=5cm,angle=0}
\psfig{figure=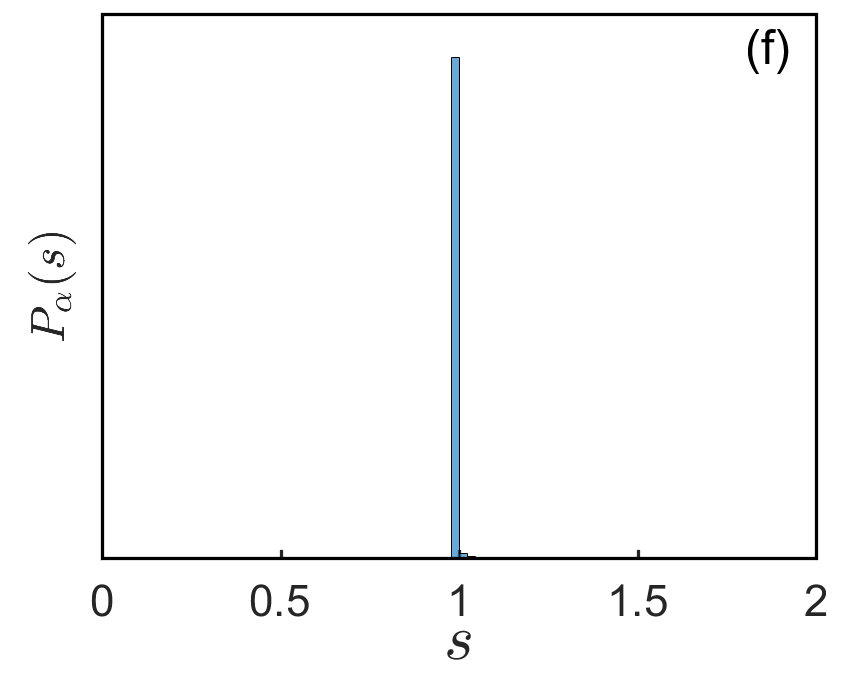,height=4cm,width=5cm,angle=0}}

%\centerline{
%\psfig{figure=gaps/gaps0p8.png,height=4cm,width=5cm,angle=0}
%\psfig{figure=gaps/gaps0p5.png,height=4cm,width=5cm,angle=0}
%\psfig{figure=gaps/gaps0p1.png,height=4cm,width=5cm,angle=0}}
\caption{The histogram of the normalized  gaps  $\{\delta_{\rm norm}^\alpha(n)\ |\ 1\le n\le N=4096\}$ of \eqref{fproblem}
with $\Omega=(-1,1)$ and $V(x)\equiv0$ for different $\alpha$:
(a) $\alpha=2.0$, (b) $\alpha=1.9$, (c) $\alpha=\sqrt{3}$, (d) $\alpha=1.5$, (e) $\alpha=1.0$,
and (f) $\alpha=0.5$.}
\label{fig:gaps1da}
\end{figure}

From Fig. \ref{fig:gaps1da}, we can conclude that the gaps distribution statistics of \eqref{fproblem}
with  $V(x)\equiv0$ is $P_\alpha(s) = \delta(s-1)$ for
$0<\alpha\le 2$.

\subsection{Eigenfunctions and their singularity characteristics}

\begin{figure}[h!]
\centerline{
\psfig{figure=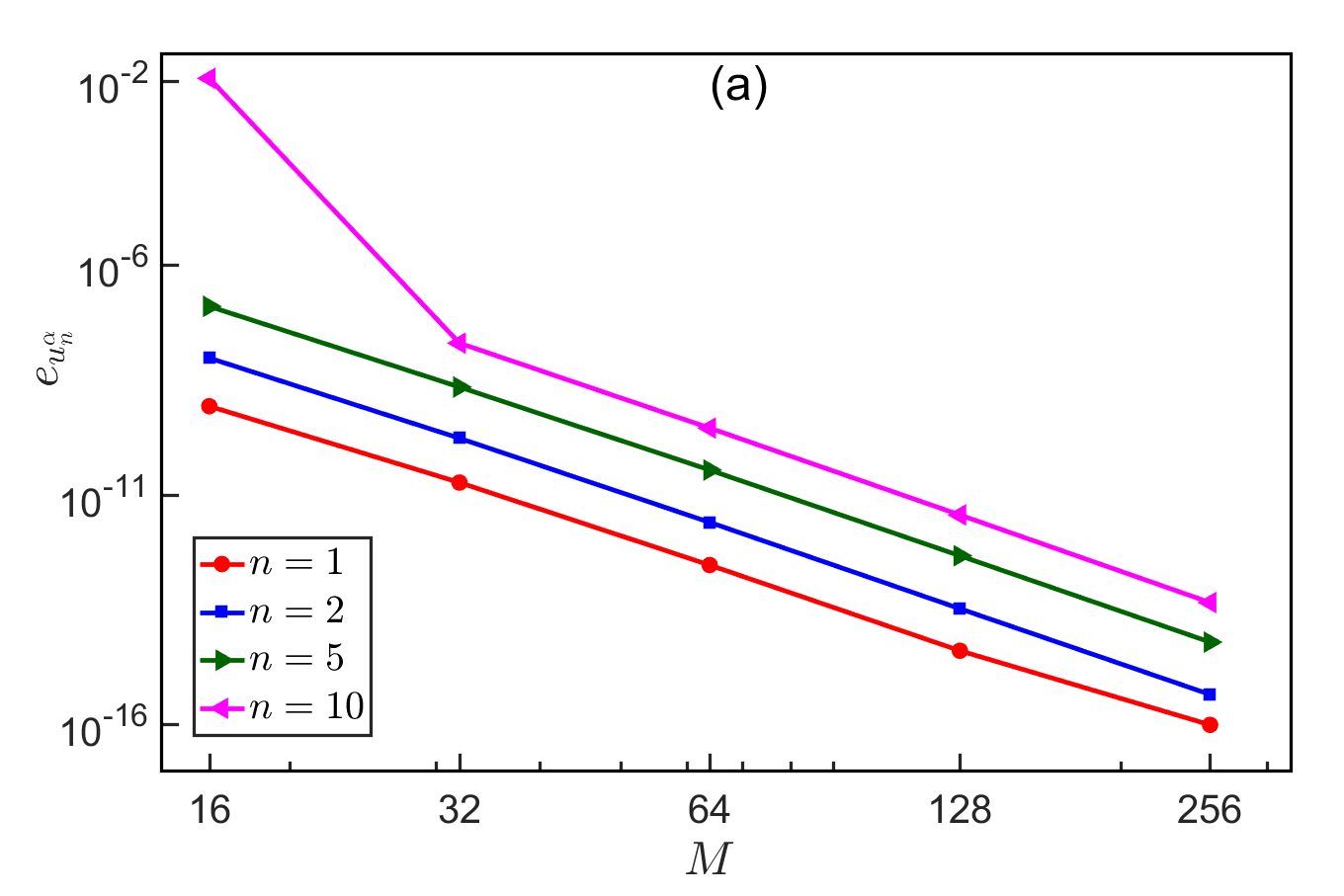,height=5cm,width=7cm,angle=0}\qquad
\psfig{figure=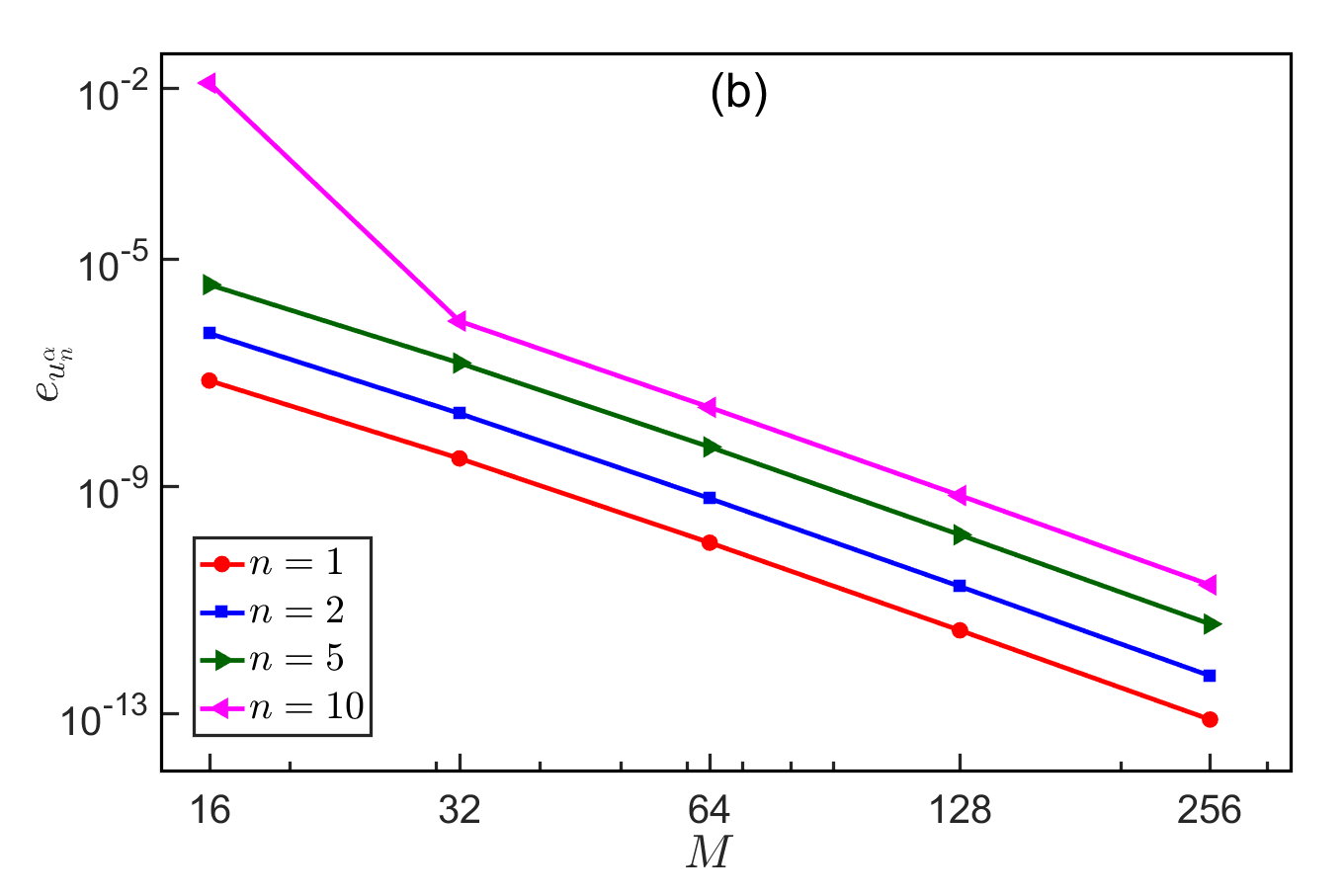,height=5cm,width=7cm,angle=0}}
\centerline{
\psfig{figure=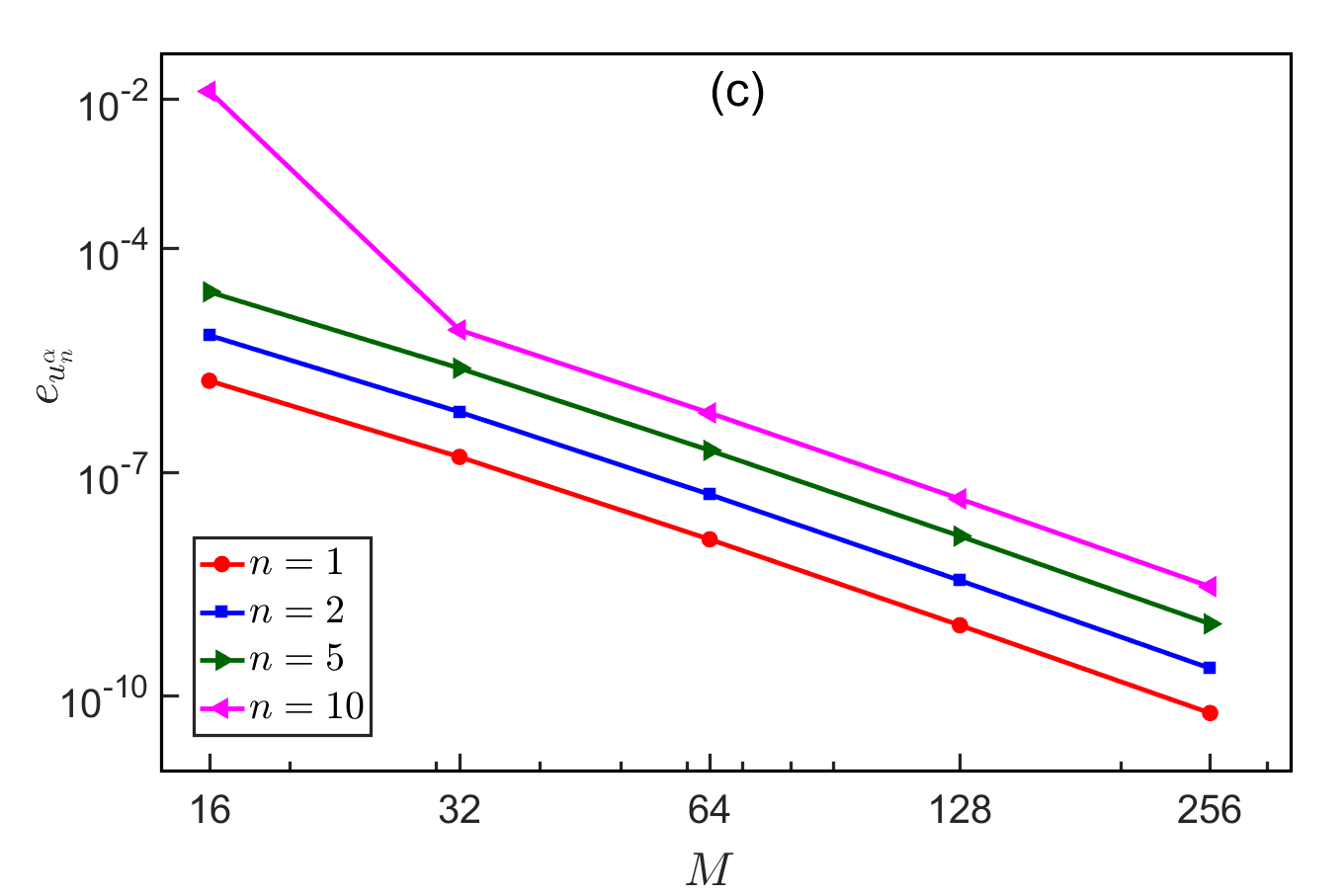,height=5cm,width=7cm,angle=0}\qquad
\psfig{figure=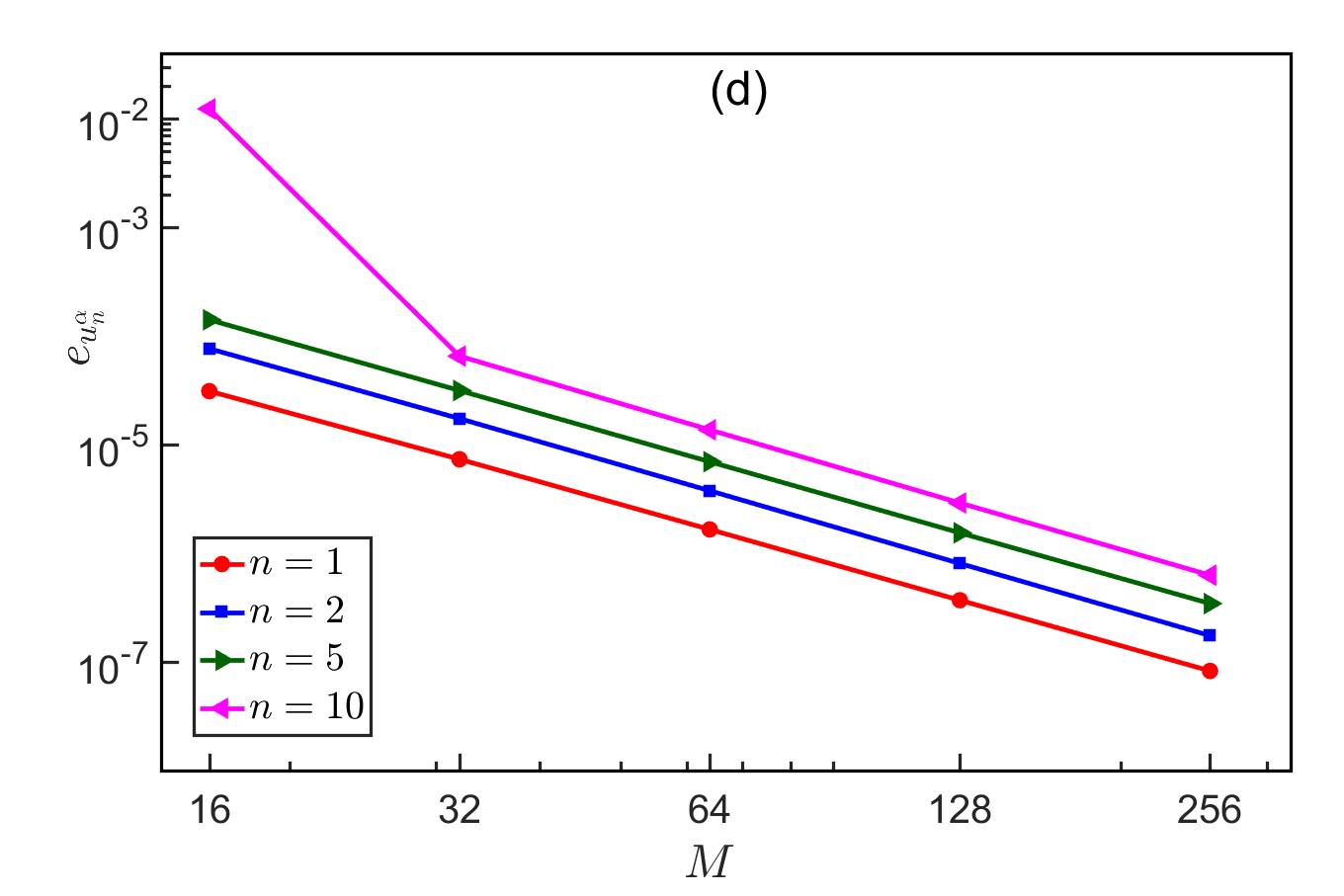,height=5cm,width=7cm,angle=0}}
\caption{Convergence rates of computing different eigenfunctions of \eqref{fproblem} with $\Omega=(-1,1)$, $V(x)\equiv0$ and  different $\alpha$
by using our JSM \eqref{weakn} for:
(a) the first eigenfunction $u_1^\alpha$, (b) the second eigenfunction $u_2^\alpha$, (c) the fifth eigenfunction $u_5^\alpha$,
and (d) the tenth eigenfunction $u_{10}^\alpha$.}
\label{fig:rateef}
\end{figure}

\begin{figure}[h!]
\centerline{
\psfig{figure=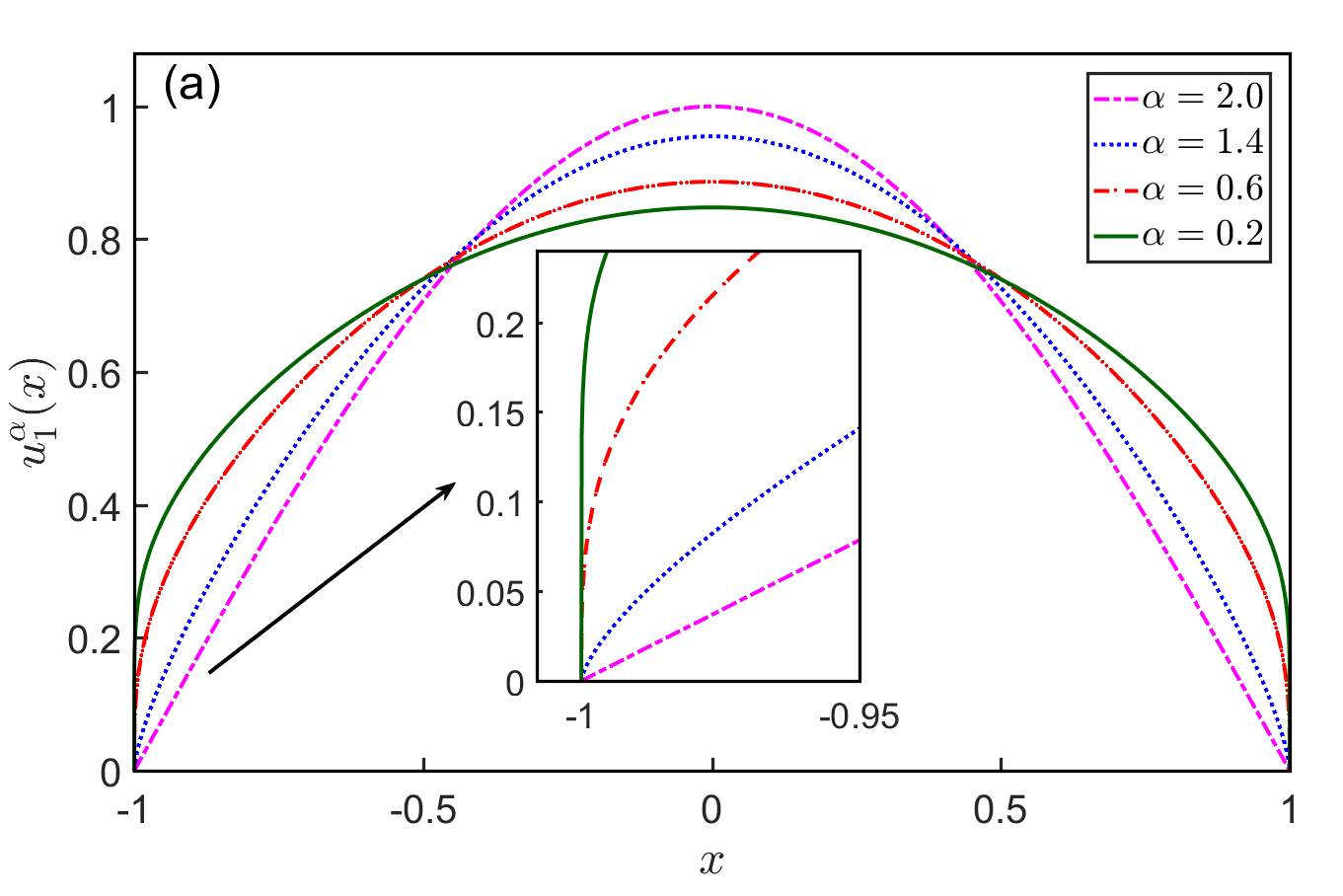,height=5cm,width=7cm,angle=0}
\psfig{figure=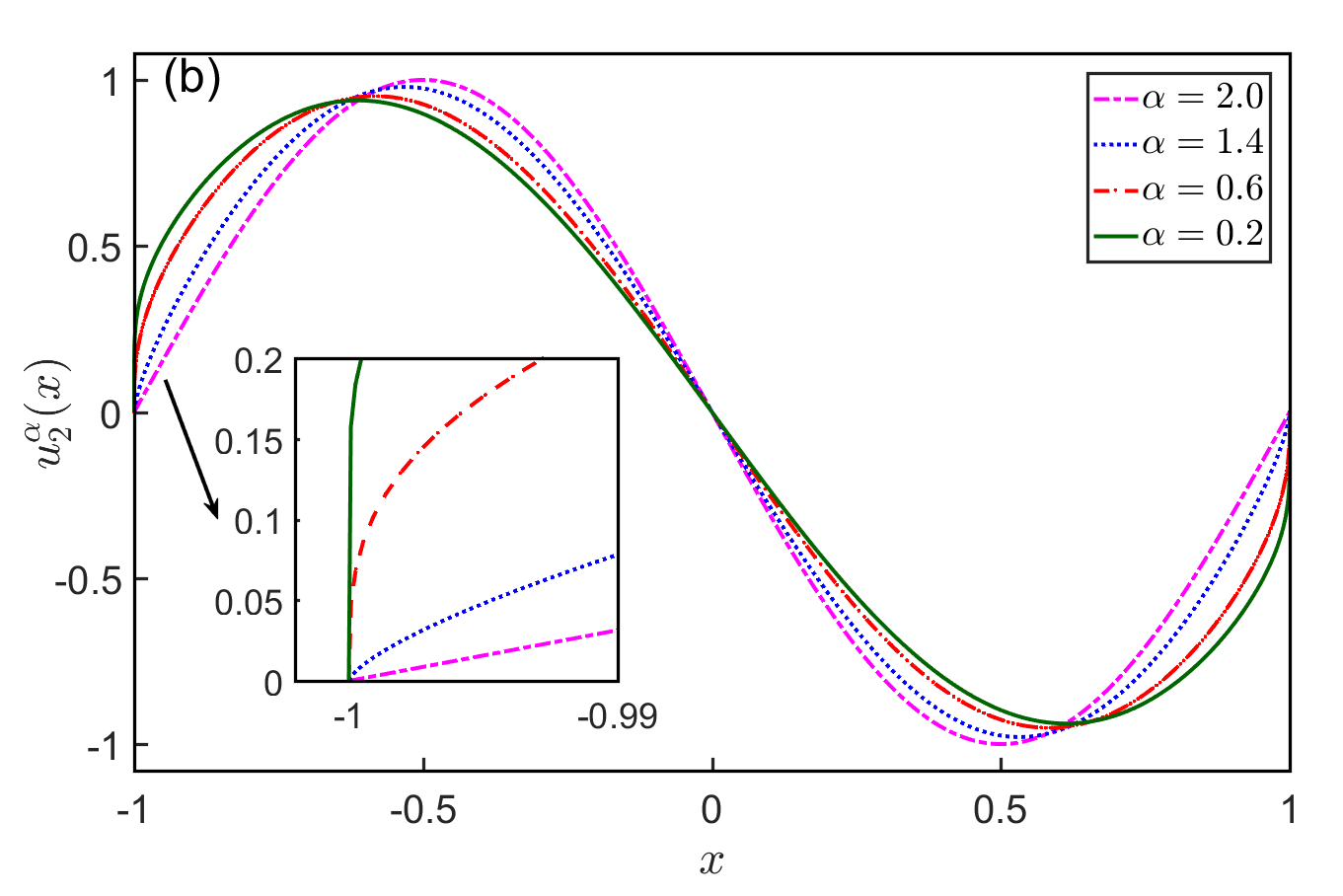,height=5cm,width=7cm,angle=0}}
\centerline{\psfig{figure=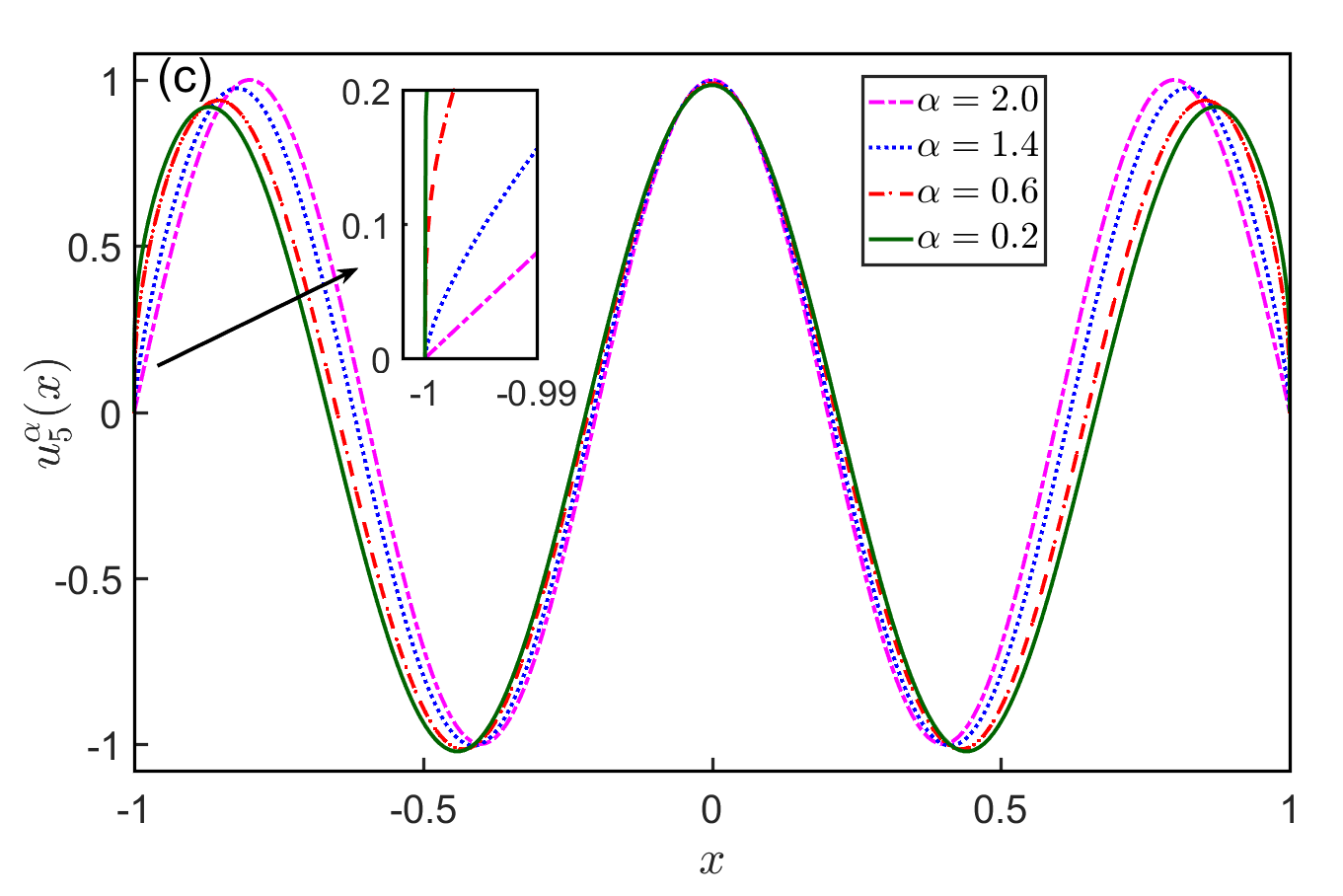,height=5cm,width=7cm,angle=0}
\psfig{figure=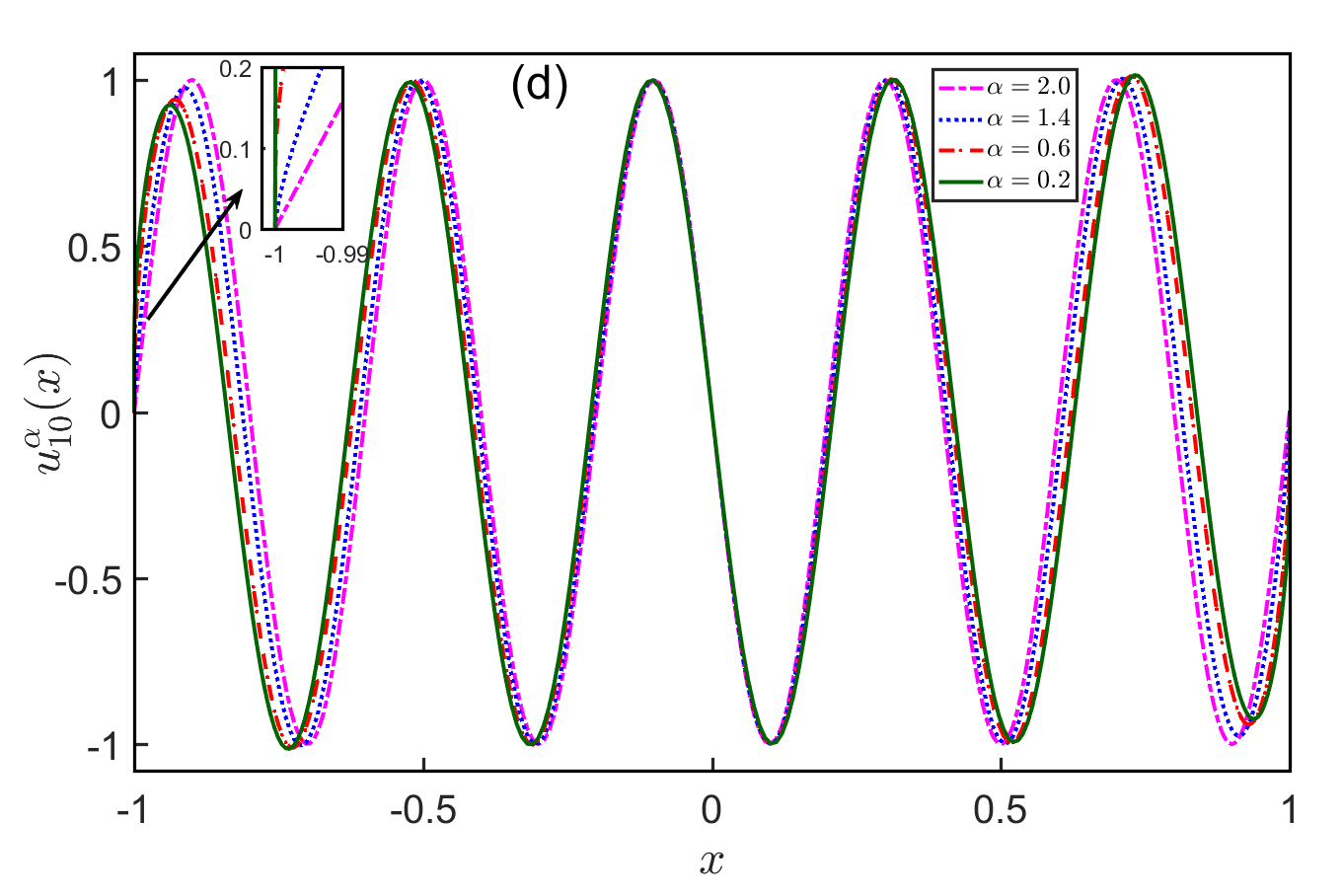,height=5cm,width=7cm,angle=0}}
\caption{Plots of different eigenfunctions of \eqref{fproblem}
with $\Omega=(-1,1)$, $V(x)\equiv0$ and different $\alpha$
for: (a) the first eigenfunction $u_1^\alpha(x)$,
(b) the second eigenfunction $u_2^\alpha(x)$,
(c) the fifth eigenfunction $u_5^\alpha(x)$,
and (d) the tenth eigenfunction $u_{10}^\alpha(x)$.}
\label{fig:eigf}
\end{figure}

Denote $u_n^\alpha(x)$ be the eigenfunction  satisfying $\|u_n^\alpha\|_{L^2(\Omega)}=1$ and $\left.\frac{d u_n^\alpha(x)}{dx}\right|_{x=-1}>0$,
which corresponds to the eigenvalue $\lambda_n^\alpha$ ($n=1,2,\ldots$) of
\eqref{fproblem} with $\Omega=(-1,1)$ and $V(x)\equiv0$.
The `exact' eigenfunctions $u_n^\alpha(x)$ ($n=1,2,\ldots$) are
obtained numerically by using the JSM \eqref{weakn} under a very large DOF $M=M_0$, e.g. $M_0=512$. Let $u_{n,M}^\alpha$ be the numerical approximation of $u_n^\alpha$ ($n=1,2,\ldots, M$) obtained by a numerical method  with the DOF chosen as $M$.
Define the absolute errors of $u_n^\alpha$ as
\be\label{errors23}
e_{u_n^\alpha}:=\|u_n^\alpha-u_{n,M}^\alpha\|_{l^2}, \qquad n=1,2,\ldots\; .
\ee
Figure \ref{fig:rateef} shows convergence rates
of our JSM \eqref{weakn}  for computing the first, second, fifth and tenth eigenfunctions of \eqref{fproblem} with $\Omega=(-1,1)$,
$V(x)\equiv0$ and  different
$\alpha$. Figure \ref{fig:eigf} plots different eigenfunctions of \eqref{fproblem} with $\Omega=(-1,1)$, $V(x)\equiv0$ and  different
$\alpha$. Finally Figure \ref{fig:eigfsig} displays
different eigenfunctions of \eqref{fproblem} with $\Omega=(-1,1)$, $V(x)\equiv0$ and  different $\alpha$ near the boundary layer $0<\xi:=x+1\ll 1$ to show the singularity characteristics of the eigenfunctions $u_n^\alpha$ at the boundary $x=-1$.

\begin{figure}[h!]
\centerline{
\psfig{figure=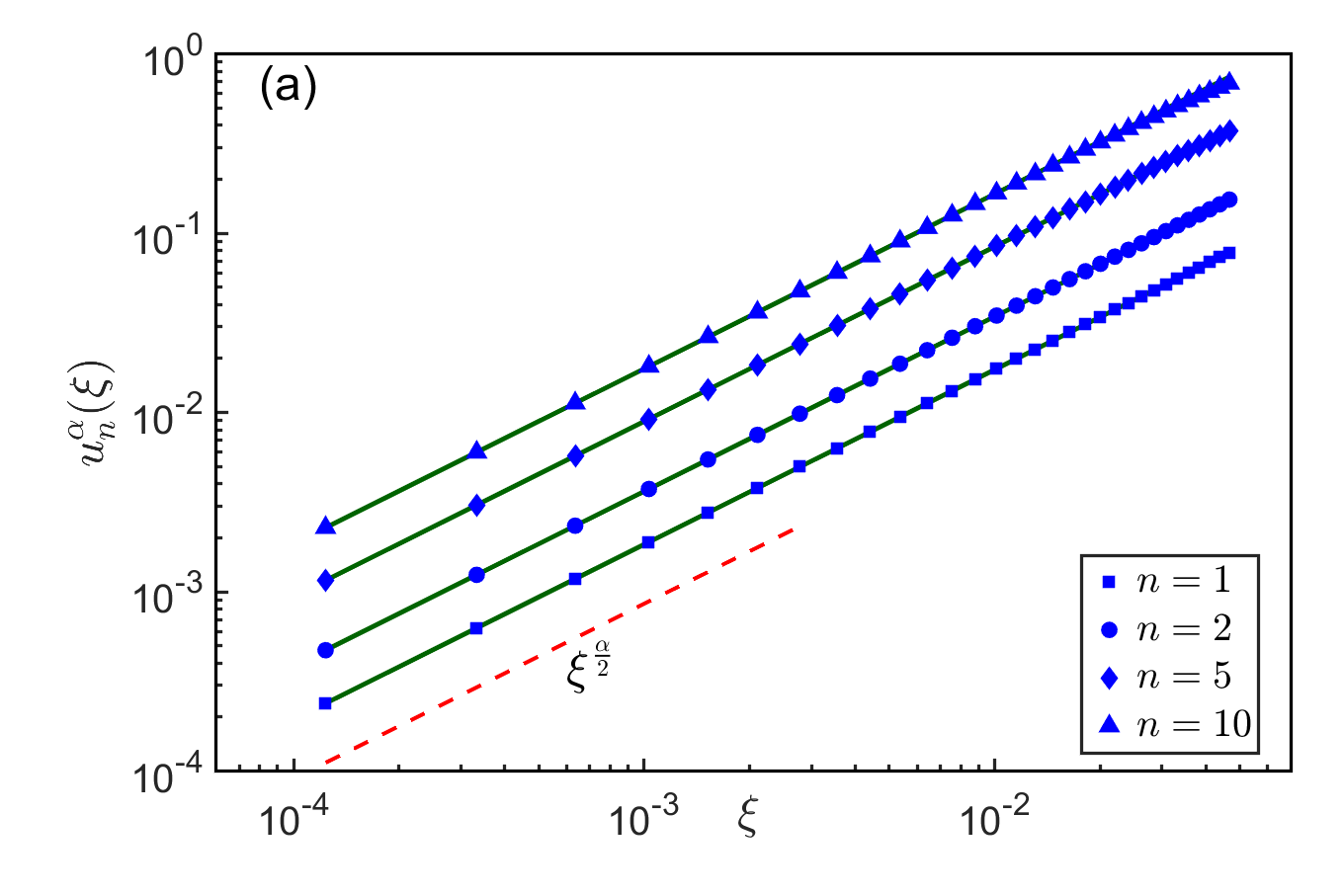,height=5cm,width=7cm,angle=0}\qquad
\psfig{figure=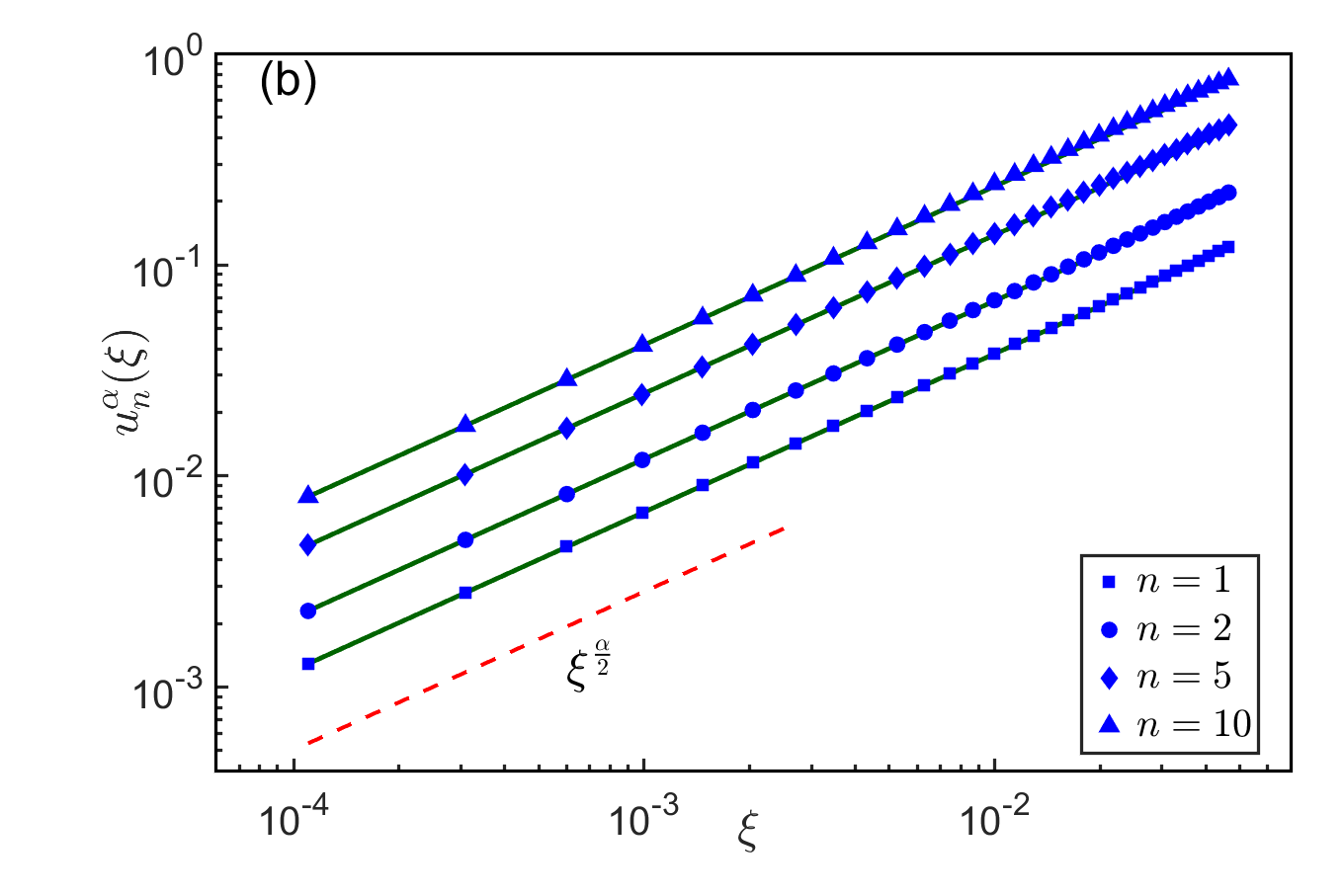,height=5cm,width=7cm,angle=0}}
\centerline{
\psfig{figure=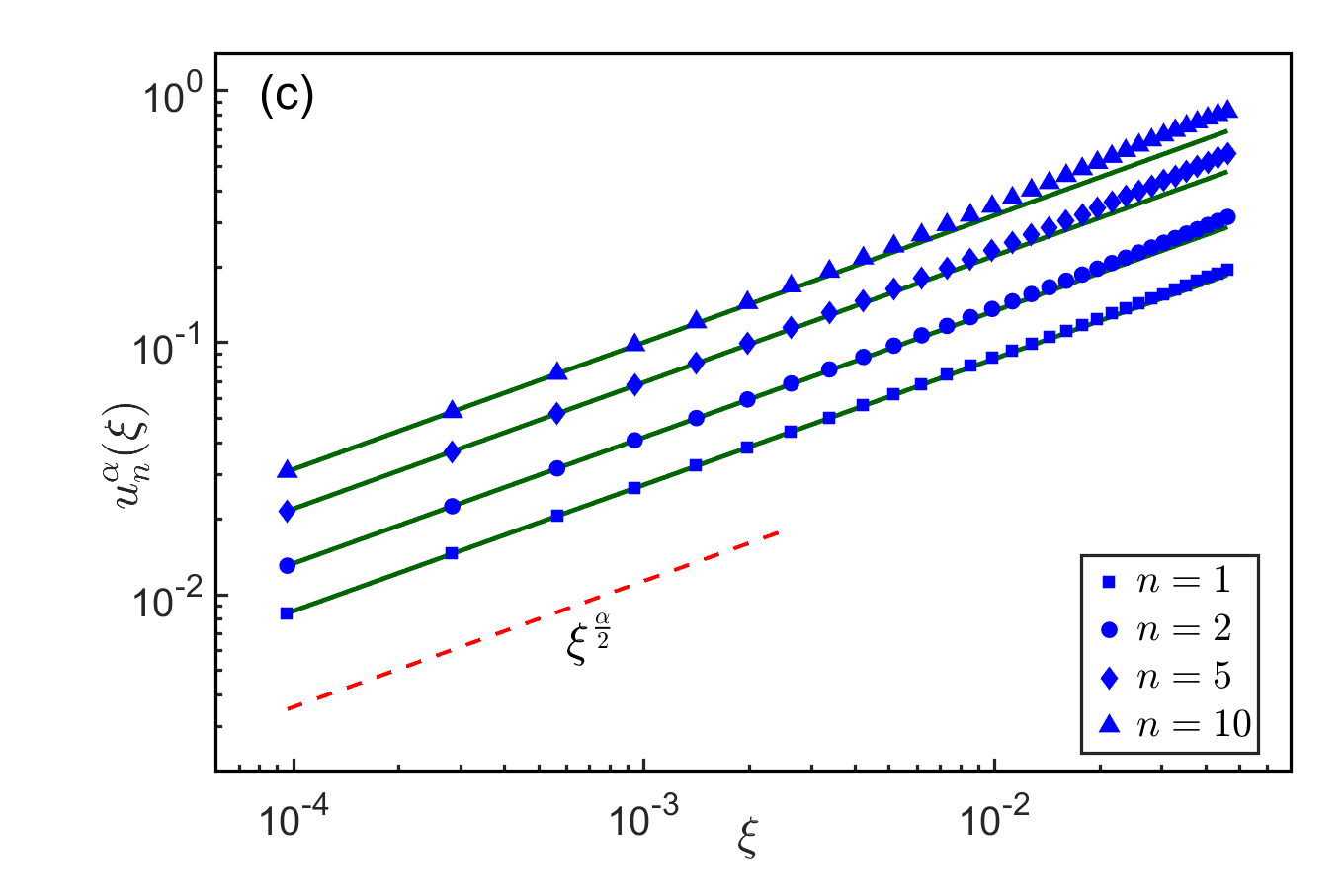,height=5cm,width=7cm,angle=0}\qquad
\psfig{figure=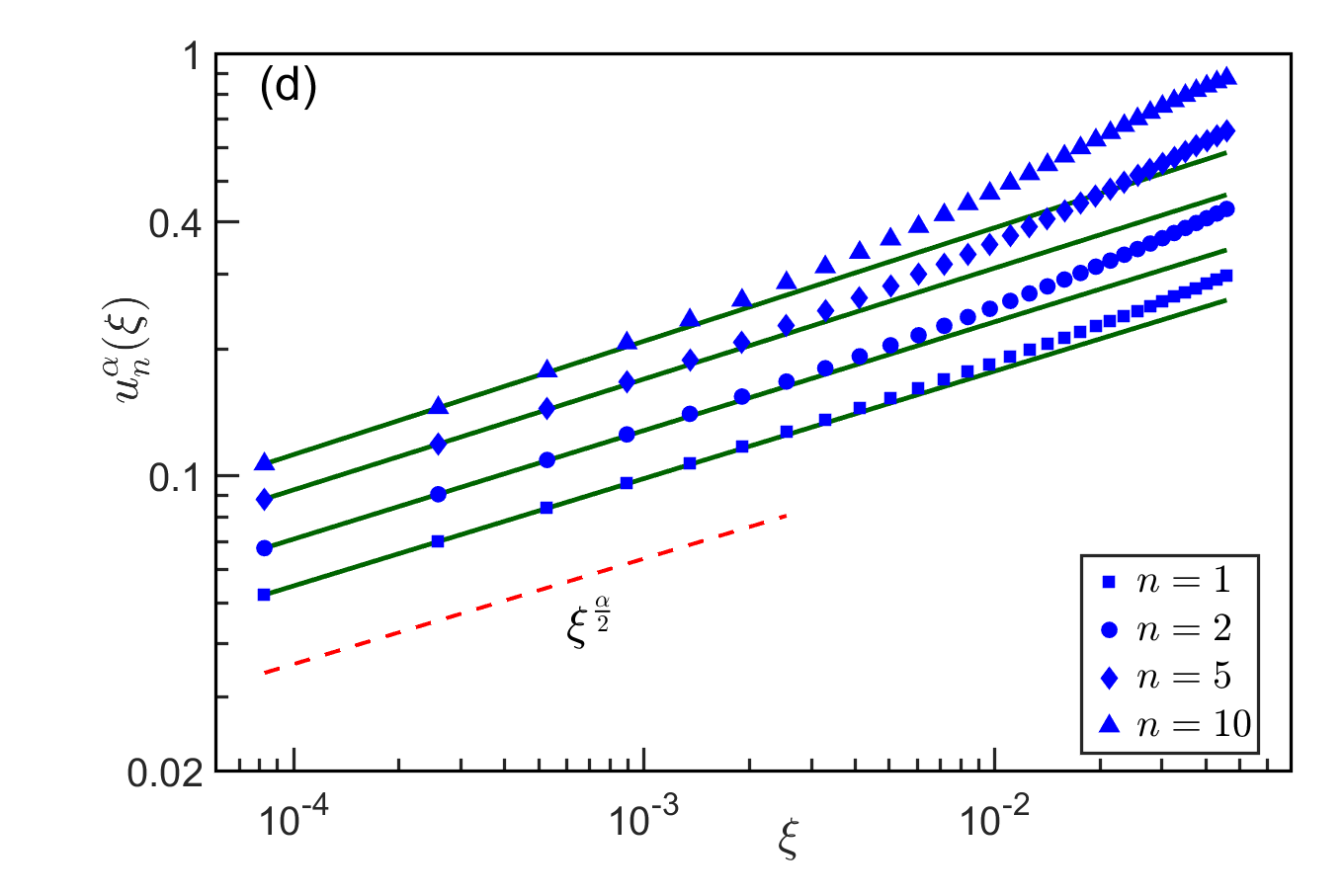,height=5cm,width=7cm,angle=0}}
\caption{Singularity characteristics of different eigenfunctions of \eqref{fproblem} with $\Omega=(-1,1)$, $V(x)\equiv0$ and  different $\alpha$
 for (symbols denote numerical results and solids lines are from fitting formula $C\xi^{\alpha/2}$ when $0<\xi=x+1\ll1$): (a) the first eigenfunction $u_1^\alpha$, (b) the second eigenfunction $u_2^\alpha$, (c) the fifth eigenfunction $u_5^\alpha$,
and (d) the tenth eigenfunction $u_{10}^\alpha$.}
\label{fig:eigfsig}
\end{figure}

From Figs. \ref{fig:rateef}-\ref{fig:eigfsig}, we can draw the following conclusions: (i) Our JSM method \eqref{weakn} converges super-linearly
 with respect to  the DOF $M$ for computing the eigenfunctions $u_n^\alpha$ (cf. Fig. \ref{fig:rateef}). (ii) For fixed $0<\alpha<2$,
the eigenfunctions $u_n^\alpha$ ($n=1,2,\ldots$) can be
characterised as
\be\label{uvx1d}
u_n^\alpha(x) =(1-x^2)^{\alpha/2}\, v_n^\alpha(x), \qquad -1\le x\le 1,
\ee
where $v_n^\alpha$ ($n=1,2,\ldots$) are smooth functions over the
interval $\bar \Omega=[-1,1]$
(cf. Fig. \ref{fig:eigfsig}). In addition, our numerical results
indicate that, when $n\to\infty$ (cf. Fig. \ref{fig:eigf}d),
the eigenfunctions $u_n^\alpha$ ($0<\alpha<2$) of
\eqref{fproblem} with $\Omega=(-1,1)$ and $V(x)\equiv0$
converge to the eigenfunction $u_n^{\alpha=2}=\sin(n\pi x)$ of
\eqref{fproblem} with $\alpha=2$, $\Omega=(-1,1)$ and $V(x)\equiv0$, i.e.
\be\label{eigenfulmt}
u_n^\alpha(x)\to \sin\left(\frac{n\pi (x+1)}{2}\right)=u_n^{\alpha=2}(x), \qquad x\in\bar{\Omega}, \qquad n\to\infty.
\ee

Based on the above results, for the eigenvalue problem of the FSO in
high dimensions:

Find $\lambda \in \mathbb{R}$ and a nonzero real-valued function $u(\bx)\ne 0$ such that
\begin{equation}\label{fproblemhdww}
\begin{split}
L_{\rm FSO}\; u(\bx)&:=\left[(-\Delta)^{\alpha/2}+V(\bx)\right]u(\bx)
=\lambda\; u(\bx), \qquad \bx\in \Omega\subset {\mathbb R}^d, \\
u(\bx)&=0, \qquad  \bx \in \Omega^c:=\mathbb{R}^d \backslash \Omega,
\end{split}
\end{equation}
where $d\ge2$, $0<\alpha<2$, $\Omega$ is a bounded domain and the fractional Laplacian $(-\Delta)^{\alpha/2}$ is defined via the Fourier transform \cite{CS07,NPV12}, we conjecture here that the eigenfunction $u(\bx)$ can be written as
\begin{equation}\label{uvxhd}
u(\bx)= v(\bx)\left({\rm dist}(\bx,\partial \Omega)\right)^{\alpha/2},
\qquad \bx\in\bar \Omega,
\end{equation}
where $v(\bx)$ is a smooth function over $\bar \Omega$ and
${\rm dist}(\bx,\partial \Omega)$ represents the distance from $\bx\in\Omega$ to $\partial \Omega$.

We remark here that the singularity
characteristics of the eigenfunctions in  \eqref{uvx1d} (or \eqref{uvxhd})
is quite different
with the singularity characteristics given in \cite{BBROA18}
for fractional PDEs as
\begin{equation}\label{uvxhhd}
u(\bx)\approx \left({\rm dist}(\bx,\partial \Omega)\right)^{\alpha/2}+v(\bx),
\qquad \bx\in\bar \Omega,
\end{equation}
where $v(\bx)$ is a smooth function over $\bar \Omega$.
From our numerical results, we speculate that  the correct singularity characteristics of the solution of fractional PDEs should be
\eqref{uvxhd} instead of \eqref{uvxhhd}!

\section{Numerical results of FSO in 1D with potential}\label{sec:pot}
\setcounter{equation}{0}

In this section, we report numerical results on eigenvalues of \eqref{fproblem} with $\Omega=(-1,1)$ and $V(x)\ne0$ by using our JSM \eqref{weakn} under the DOF $M=8192$. All results are based on the first
$4096$ eigenvalues, i.e. we use half of the eigenvalues obtained numerically
to present the results and to calculate distribution statistics.
Here we consider four different external potentials given as:

Case I. ~~$V(x) = \frac{x^2}{2}$;

Case II. ~$V(x)=4 x^2$; % \frac{(x-0.1)^2}{2}$;

Case III. $V(x)=4x^2+\sin(\frac{\pi}{2} x)$;   %x^2+\sin(x)$;

Case IV. $V(x)=50x^2+\sin(2\pi x)$.  % x^2+ \cos(x)$.

\subsection{Eigenvalues and their asymptotics}

\begin{table}[h!]
\centering
\begin{tabular}{ |c |c|c|c|c|c|} \hline
    & $\alpha=0.5$          &  $\alpha=1.0$          &  $\alpha=1.5$            & $\alpha=1.9$ & $\alpha=2.0$ \\ \hline
$\lambda_1^\alpha$    &1.0599238    &  1.240244372 &  1.6707307180 &2.31063679348 & 2.53245197432  \\
$\lambda_2^\alpha$    &1.7684725    &  2.918074603 &  5.2120578091 &8.73899699079 & 10.0106621605 \\
$\lambda_3^\alpha$    &2.1903345    &  4.481368142 &  9.7550085449 &18.8734566366 & 22.3620761310\\
$\lambda_4^\alpha$    &2.5518267    &  6.058660406 &  15.182580104 &32.6230979973 & 39.6388288214\\
$\lambda_5^\alpha$    &2.8580498    &  7.626501974 &  21.354271585 &49.8832020720 & 61.8477048695\\
$\lambda_6^\alpha$    &3.1370031    &  9.199495156 &  28.200700106 &70.5802261928 & 88.9903414346\\
$\lambda_7^\alpha$    &3.3893161    &  10.76885112 &  35.653816621 &94.6494682651 & 121.067291745\\
$\lambda_8^\alpha$    &3.6251388    &  12.34077821 &  43.673146060 &122.040857583 & 158.078785000\\
$\lambda_9^\alpha$    &3.8445549    &  13.91072820 &  52.217197374 &152.708819987 & 200.024930128\\
$\lambda_{10}^\alpha$ &4.0526430    &  15.48221913 &  61.258734930 &186.615849002 & 246.905784303\\
$\lambda_{20}^\alpha$ &5.5522311    &  31.02330310 &  174.43784577 &697.513597025 & 986.960440109\\
$\lambda_{40}^\alpha$ &7.8894197    &  62.43917340 &  495.71364899 &2606.30876720 & 3947.84176043\\
$\lambda_{60}^\alpha$ &9.6777480    &  93.85508927 &  912.11187382 &5633.40862247 & 8882.64396098\\
\hline
\end{tabular}
\caption{Different eigenvalues of \eqref{fproblem} with $\Omega=(-1,1)$, $V(x)=\frac{x^2}{2}$ and different $\alpha$ obtained numerically by our JSM \eqref{weakn}.}
\label{eigs}
\end{table}

\begin{figure}[h!]
\centerline{
\psfig{figure=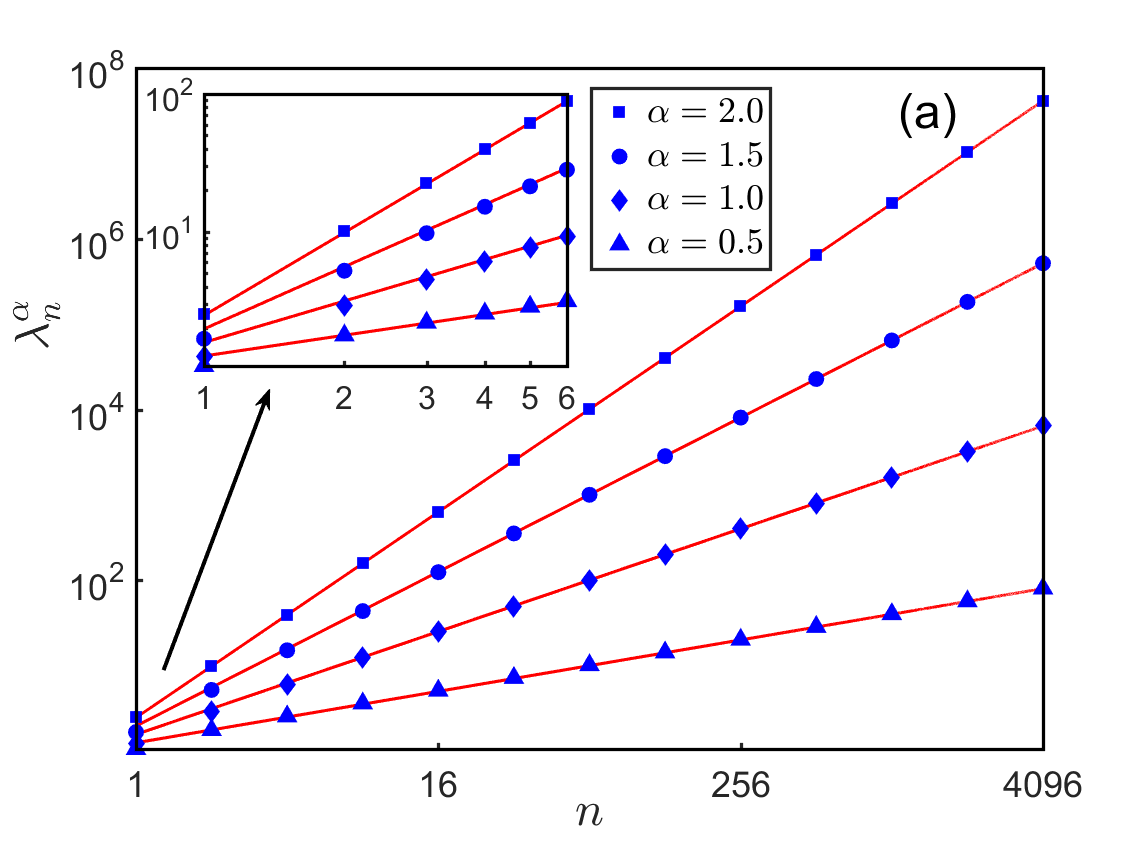,height=5cm,width=7cm,angle=0}
\psfig{figure=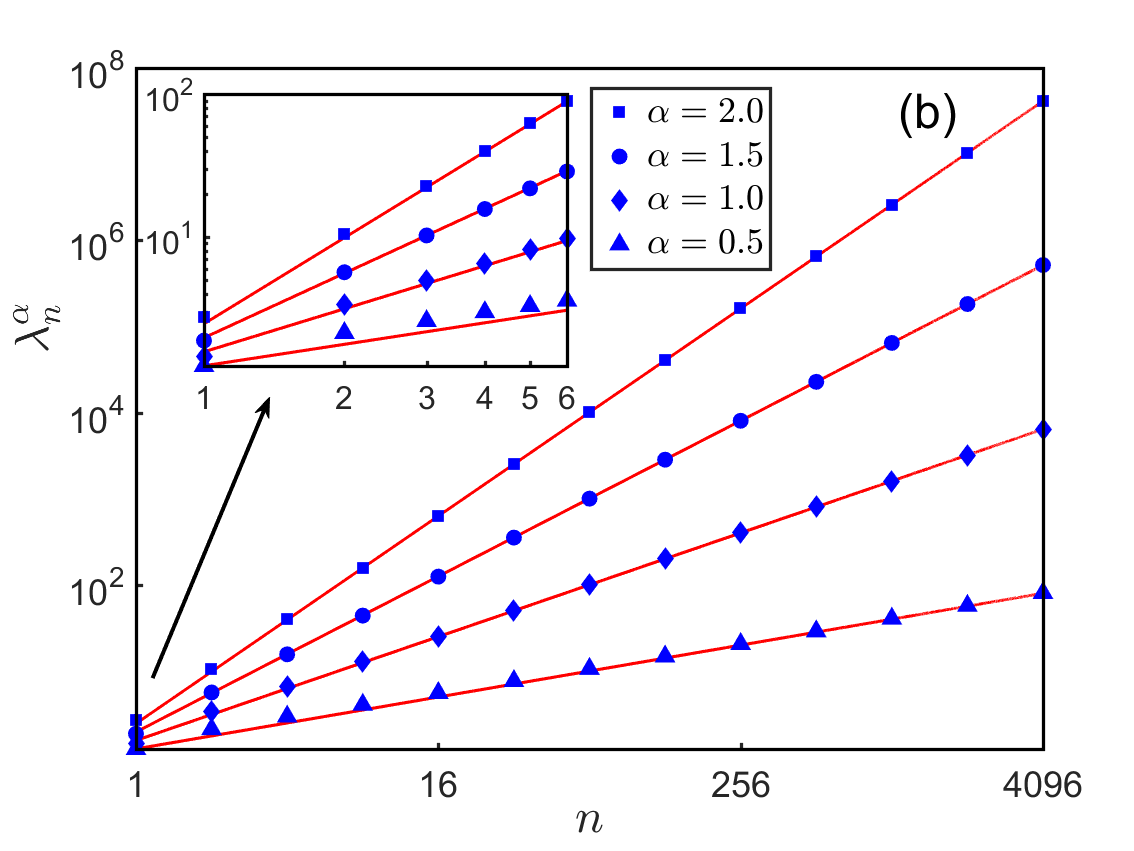,height=5cm,width=7cm,angle=0}}
\centerline{\psfig{figure=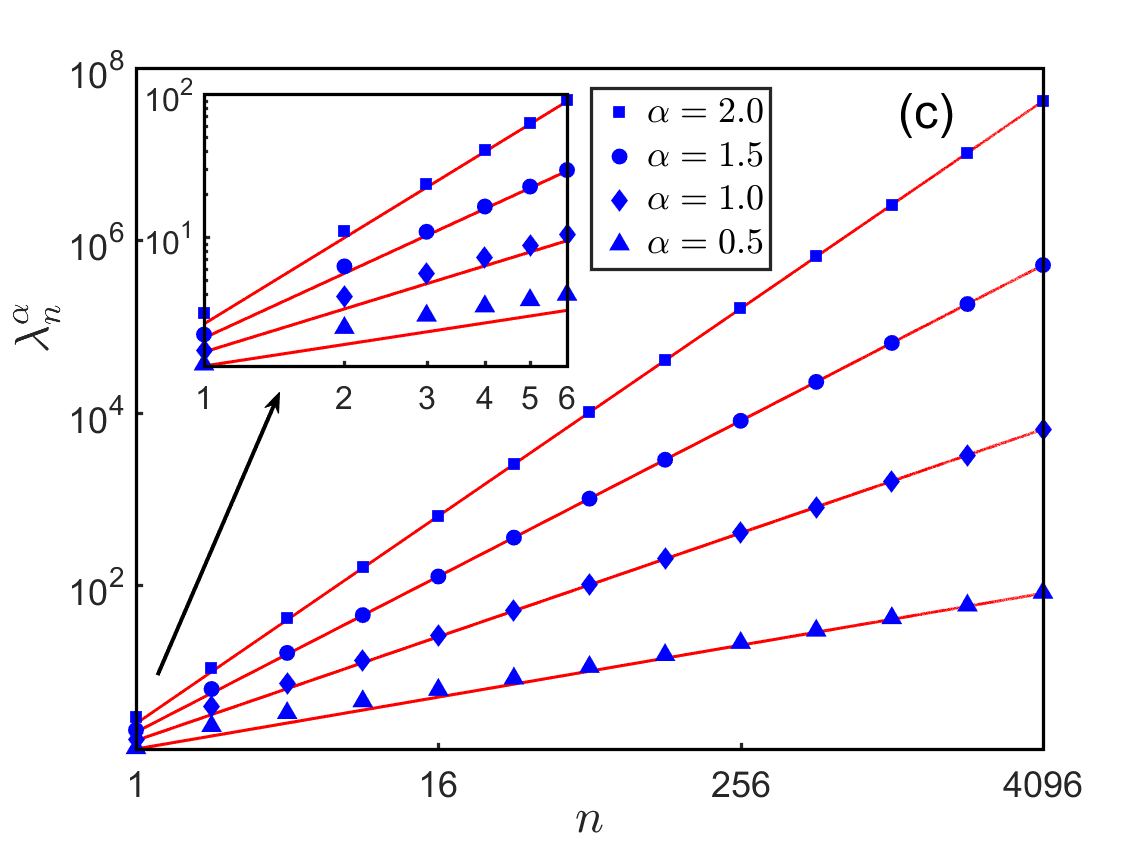,height=5cm,width=7cm,angle=0}
\psfig{figure=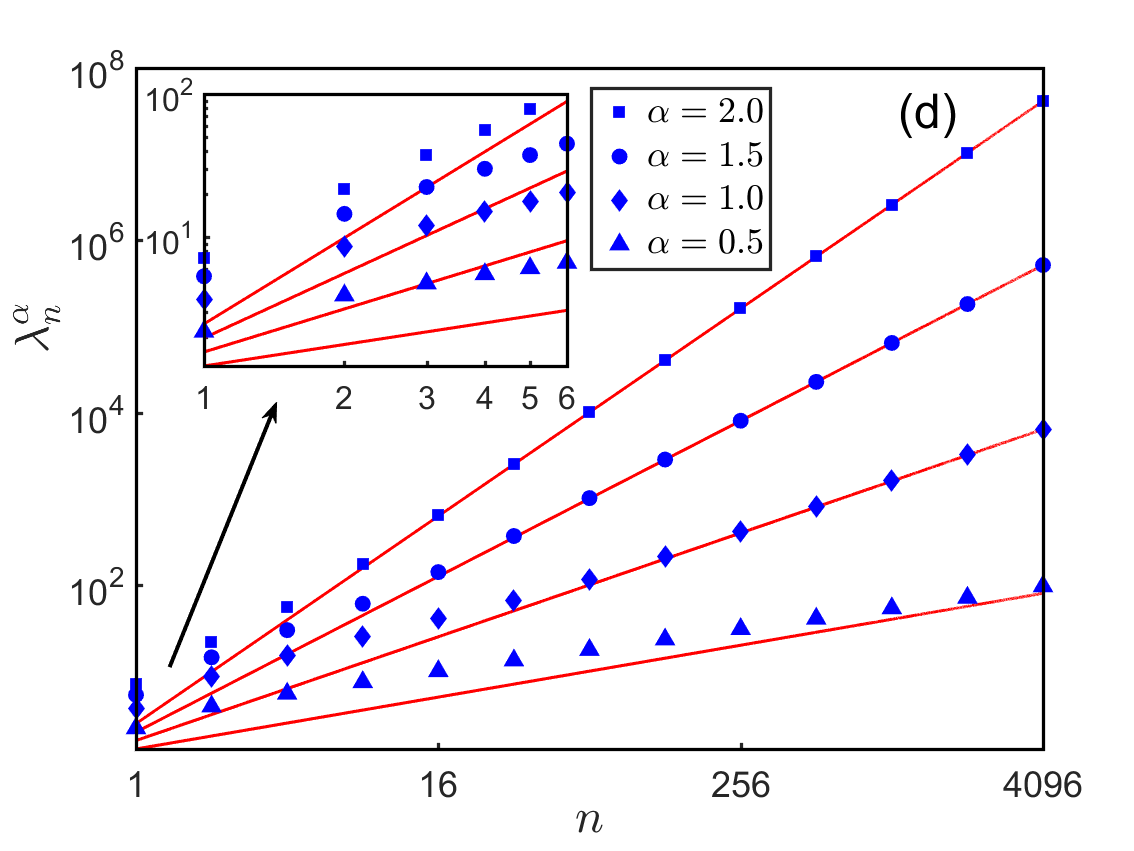,height=5cm,width=7cm,angle=0}}
\caption{Eigenvalues $\lambda_n^\alpha$ ($n=1,2,\ldots,4096$) of \eqref{fproblem} with $\Omega=(-1,1)$ and different $\alpha$ for
differential external potentials (symbols denote numerical results and solids lines are from fitting formula when $n\gg1$):
(a) Case I, (b) Case II, (c) Case III, and (d) Case IV. }
\label{fig:eig1dnpepwp}
\end{figure}

\begin{figure}[h!]
\centerline{
\psfig{figure=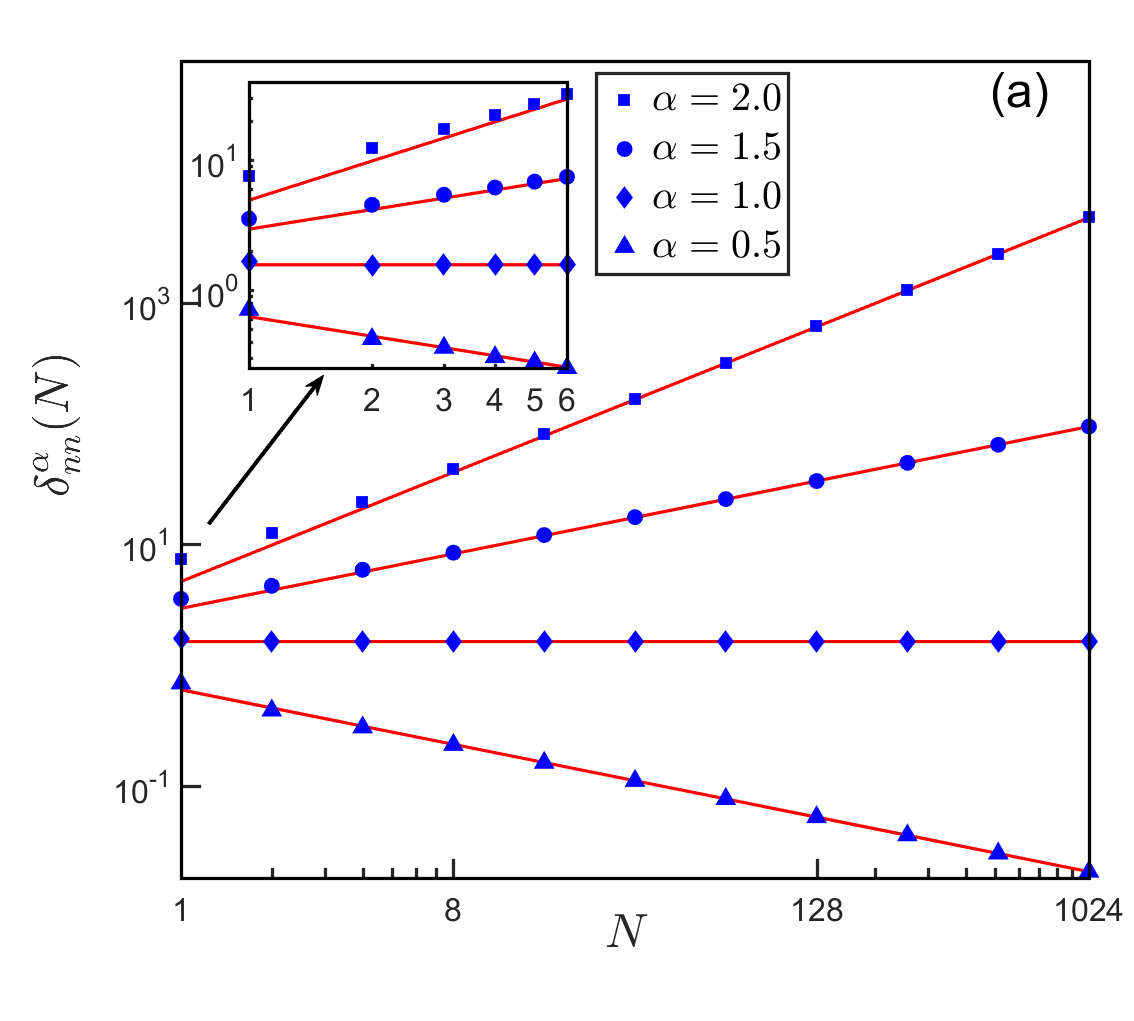,height=5cm,width=7cm,angle=0}
\psfig{figure=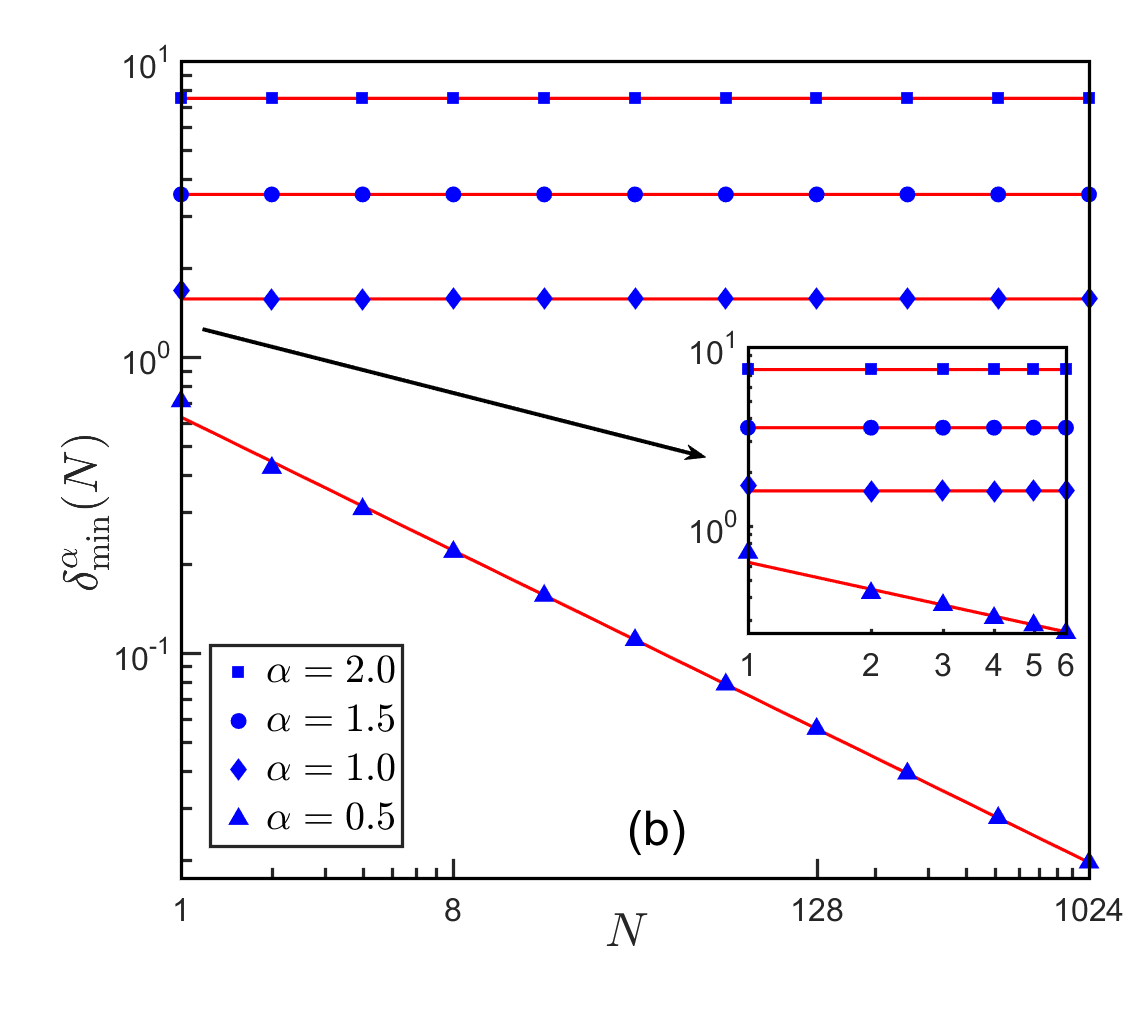,height=5cm,width=7cm,angle=0}}
\centerline{\psfig{figure=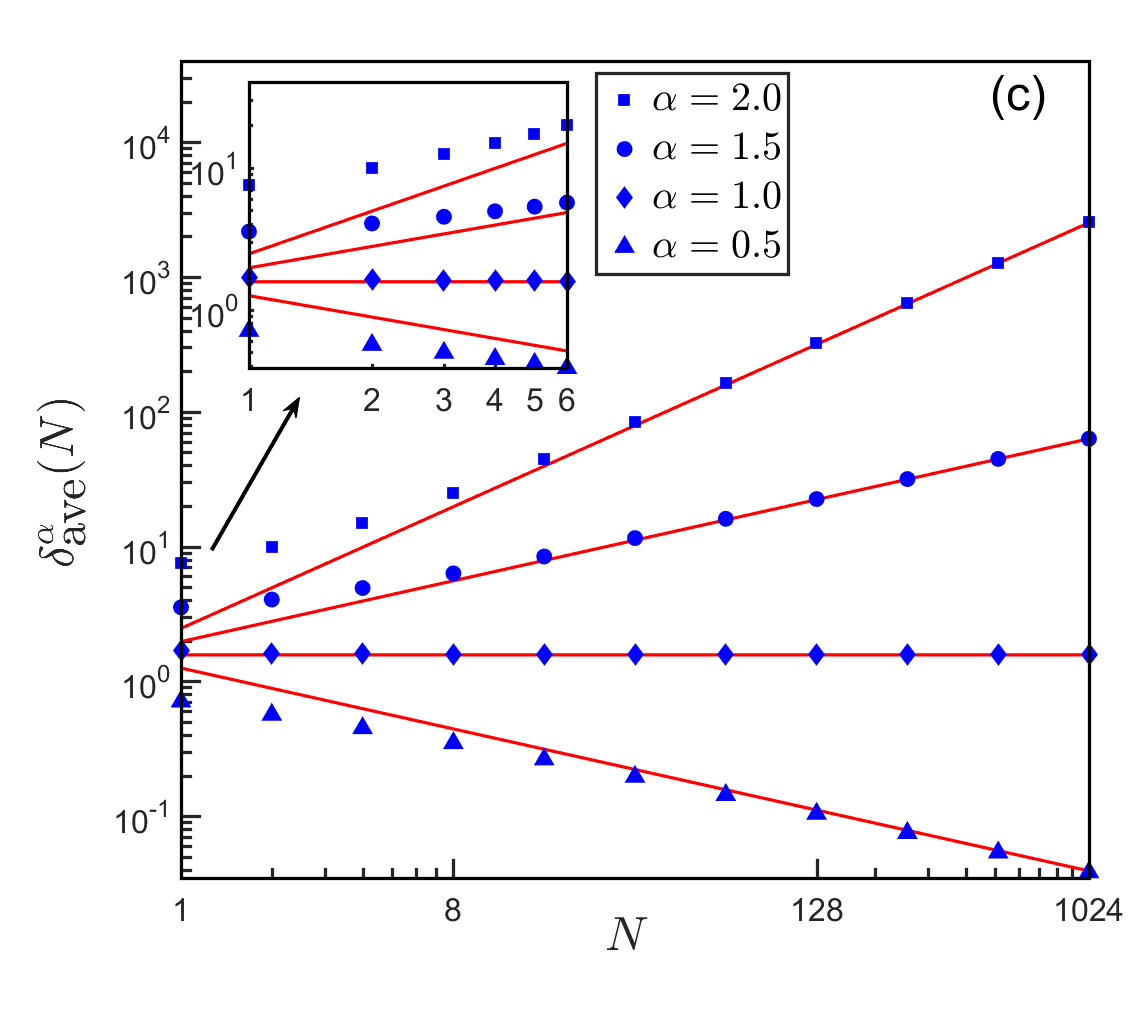,height=5cm,width=7cm,angle=0}
\psfig{figure=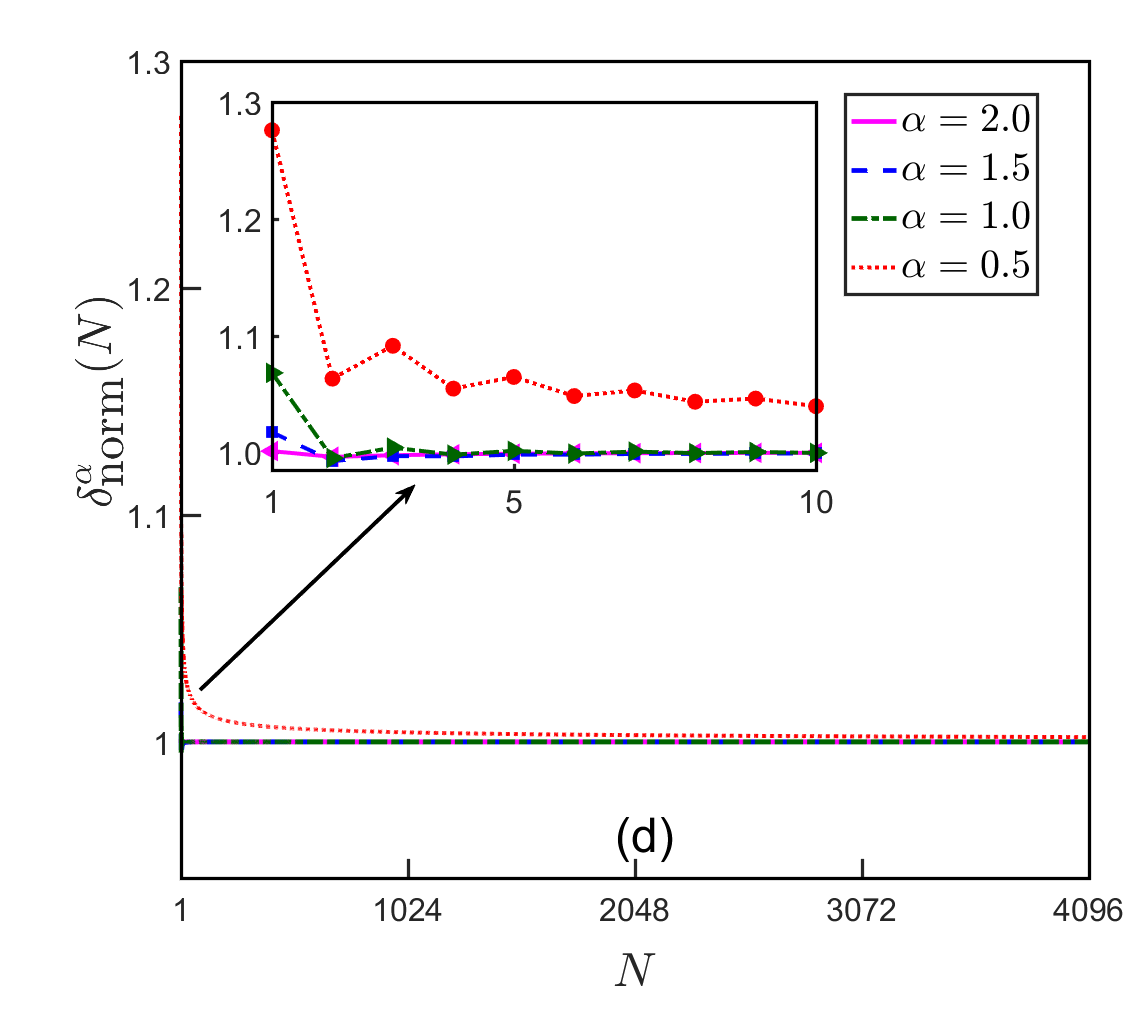,height=5cm,width=7cm,angle=0}}
\caption{Different eigenvalue gaps of \eqref{fproblem}
with $\Omega=(-1,1)$, $V(x)=\frac{x^2}{2}$ and different $\alpha$
for (symbols denote numerical results and solids lines are from fitting formula  when $N\gg1$ in a-c): (a) the nearest neighbour gaps $\delta_{\rm nn}^\alpha(N)$, (b)
the minimum gaps $\delta_{\textrm{min}}^\alpha(N)$,
 (c) the average gaps  $\delta_{\textrm{ave}}^\alpha(N)$,
 and (d) the normalized gaps  $\delta_{\textrm{norm}}^\alpha(N)$. }
\label{fig:gapsharwp}
\end{figure}

 Table \ref{eigs}
lists the eigenvalues of \eqref{fproblem} with $\Omega=(-1,1)$ and $V(x)=\frac{x^2}{2}$  for different $\alpha$.
Figure \ref{fig:eig1dnpepwp} plots the eigenvalues of \eqref{fproblem} with $\Omega=(-1,1)$, different external potentials $V(x)$  and different $\alpha$.

From Fig. \ref{fig:eig1dnpepwp}, we can conclude that,
when $n\gg1$,  the leading order asymptotics of the eigenvalues $\lambda_n^\alpha$
in \eqref{asy-eig-1Dd11} is still valid for the eigenvalue problem of FSO \eqref{fproblem} with potential $V(x)$.

\subsection{Gaps and their distribution statistics}

\begin{figure}[h!]
\centerline{
\psfig{figure=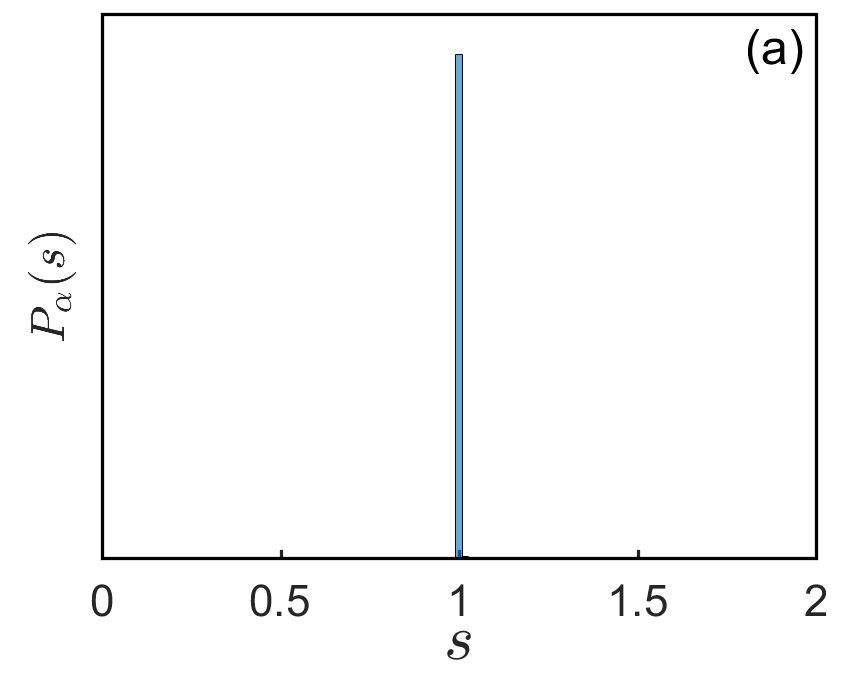,height=4cm,width=5cm,angle=0}
\psfig{figure=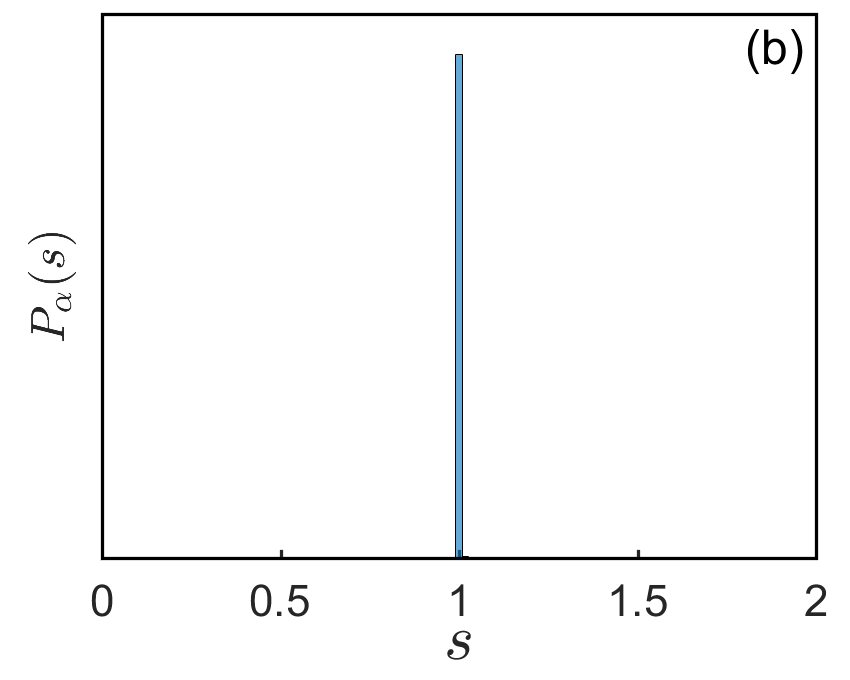,height=4cm,width=5cm,angle=0}
\psfig{figure=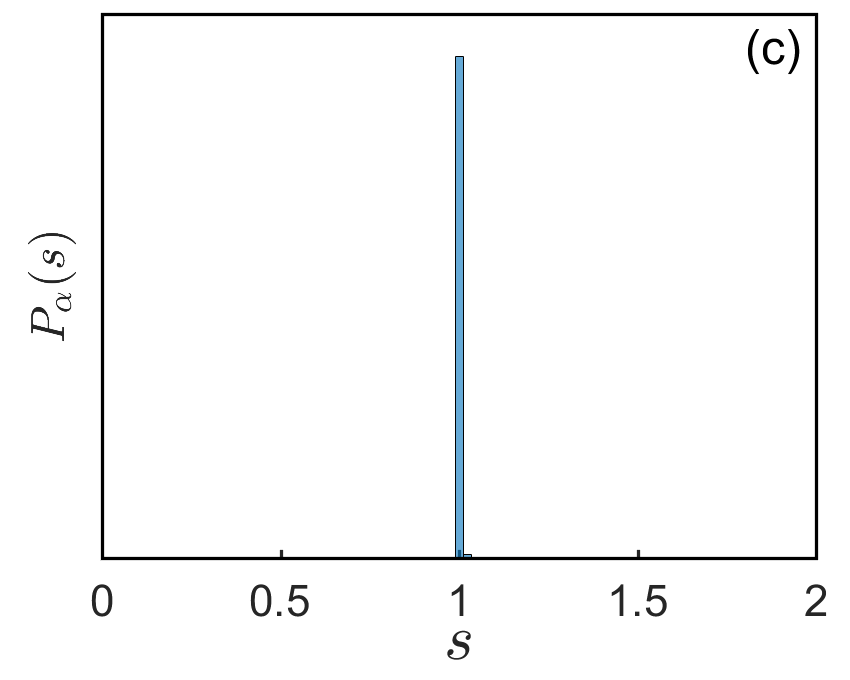,height=4cm,width=5cm,angle=0}}
\centerline{
\psfig{figure=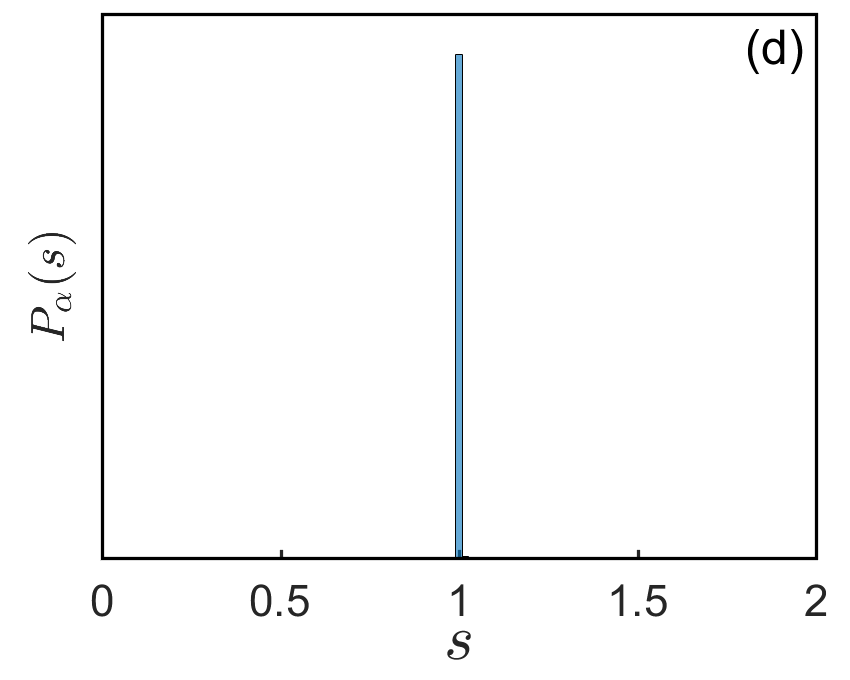,height=4cm,width=5cm,angle=0}
\psfig{figure=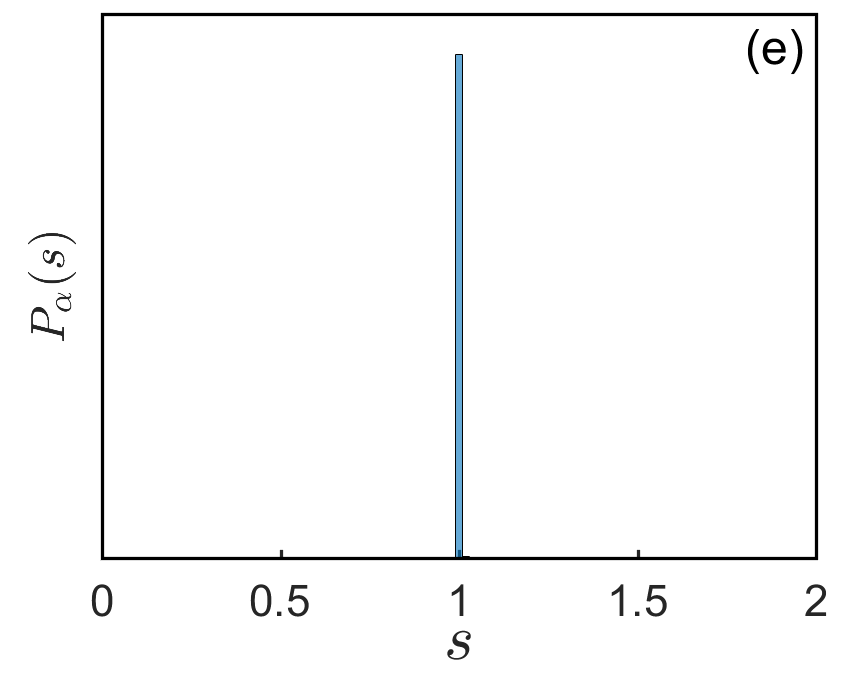,height=4cm,width=5cm,angle=0}
\psfig{figure=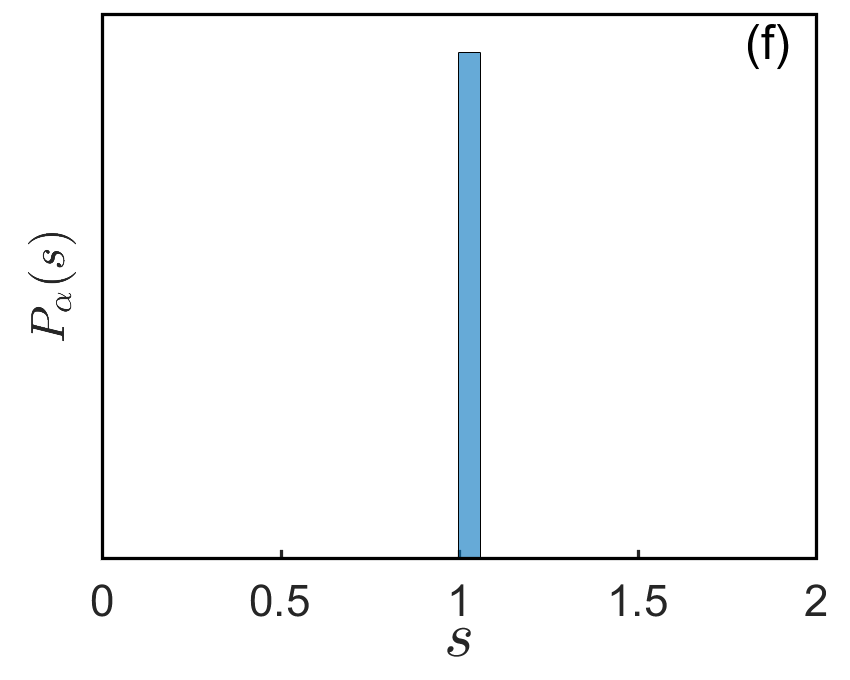,height=4cm,width=5cm,angle=0}}
%\centerline{
%\psfig{figure=gaps/gaps0p8.png,height=4cm,width=5cm,angle=0}
%\psfig{figure=gaps/gaps0p5.png,height=4cm,width=5cm,angle=0}
%\psfig{figure=gaps/gaps0p1.png,height=4cm,width=5cm,angle=0}}
\caption{The histogram of the normalized gaps  $\{\delta_{\rm norm}^\alpha(n)\ | \ 1\le n\le N=4096\}$ of \eqref{fproblem}
with $\Omega=(-1,1)$ and $V(x)=\frac{x^2}{2}$ for different $\alpha$:
(a) $\alpha=2.0$, (b) $\alpha=1.9$, (c) $\alpha=\sqrt{3}$, (d) $\alpha=1.5$, (e) $\alpha=1.0$,
and (f) $\alpha=0.5$. }
\label{fig:gaps1dawp}
\end{figure}

Figure \ref{fig:gapsharwp} plots different eigenvalue gaps of \eqref{fproblem}
with $\Omega=(-1,1)$, $V(x)=\frac{x^2}{2}$ and different $\alpha$.
Figure \ref{fig:gaps1dawp} displays the histogram of the normalized gaps $\{\delta_{\rm norm}^\alpha(n)\ |\ 1\le n\le N=4096\}$ defined in \eqref{normgapar} for \eqref{fproblem}
with $\Omega=(-1,1)$, $V(x)=\frac{x^2}{2}$ and different $\alpha$.
For other potentials, our numerical results show similar behavior
on eigenvalues and their gaps, which are omitted here for brevity.

Again, from Figs. \ref{fig:gapsharwp} and \ref{fig:gaps1dawp}, we can conclude that,
when $n\gg1$,  the asymptotics of the eigenvalue gaps given
in \eqref{nngapar}, \eqref{mingapar}, \eqref{avegapar1},
 \eqref{avegapar2}, \eqref{avegapar3} and \eqref{normgapar}
 are still valid for the eigenvalue problem of FSO \eqref{fproblem} with potential $V(x)$. In addition, the gaps distribution statistics is still $P_\alpha(s) = \delta(s-1)$ for
$0<\alpha\le 2$ in this case.

\subsection{Comparison on eigenvalues of \eqref{fproblem} without and with potential}

Let $0<\lambda_1^{\alpha,0}<\lambda_2^{\alpha,0}<\ldots<\lambda_n^{\alpha,0}<\ldots$ be all eigenvalues of \eqref{fproblem} with $\Omega=(-1,1)$ and $V(x)\equiv 0$, and denote all eigenvalues of \eqref{fproblem} with a potential $V$ as in \eqref{eigenorder}. Figure \ref{fig:eigdiffv} plots differences of the eigenvalues of \eqref{fproblem} with a potential $V$ and without a potential, i.e. $\delta_n^V:=\lambda_n^\alpha-\lambda_n^{\alpha,0}-C_V$ ($1\le n\le N=4096$) for different potentials $V(x)$ and $\alpha$, where $C_V=\frac{1}{2}\int_{-1}^1 V(x)dx$.

\begin{figure}[h!]
\centerline{
\psfig{figure=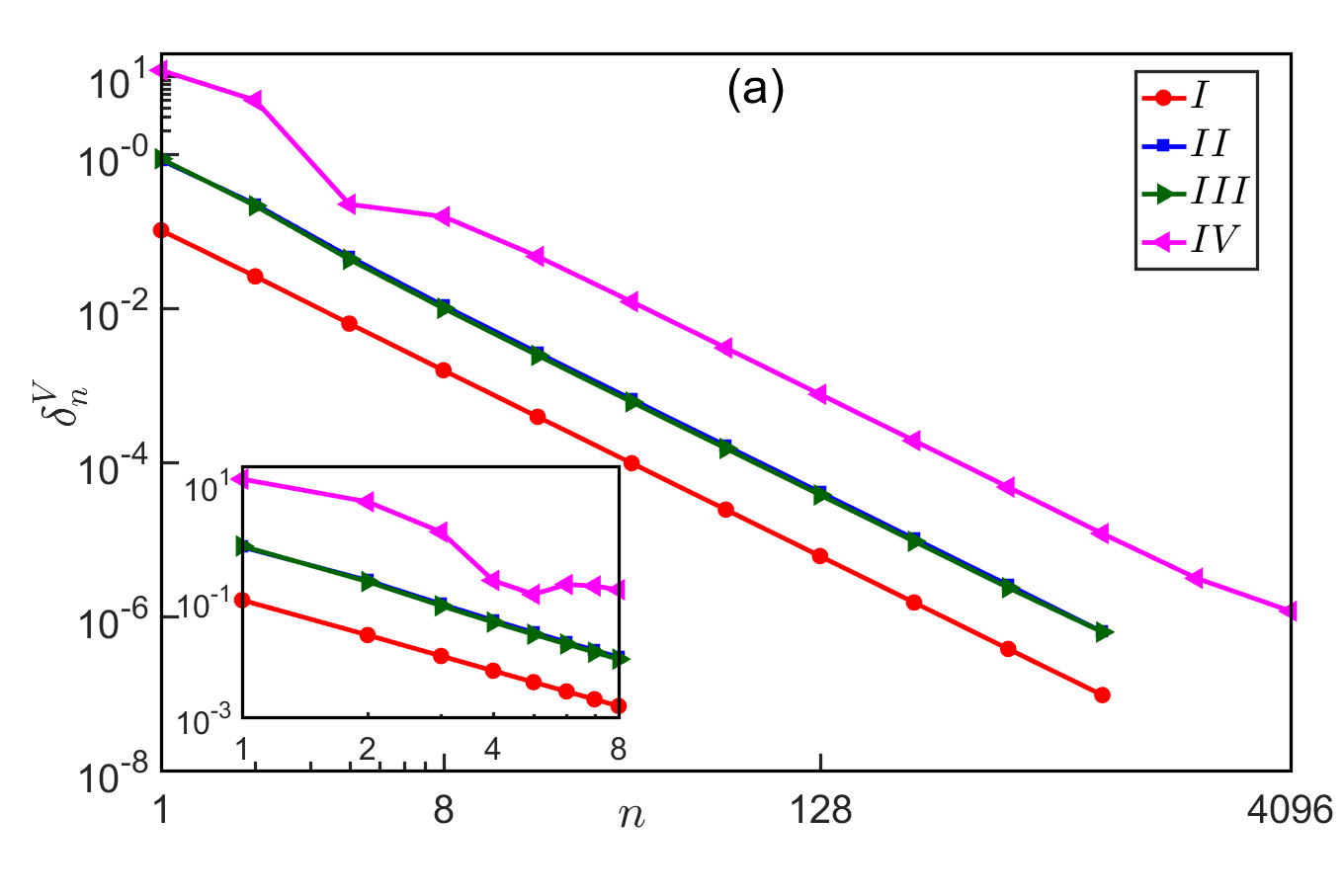,height=5cm,width=7cm,angle=0}
\psfig{figure=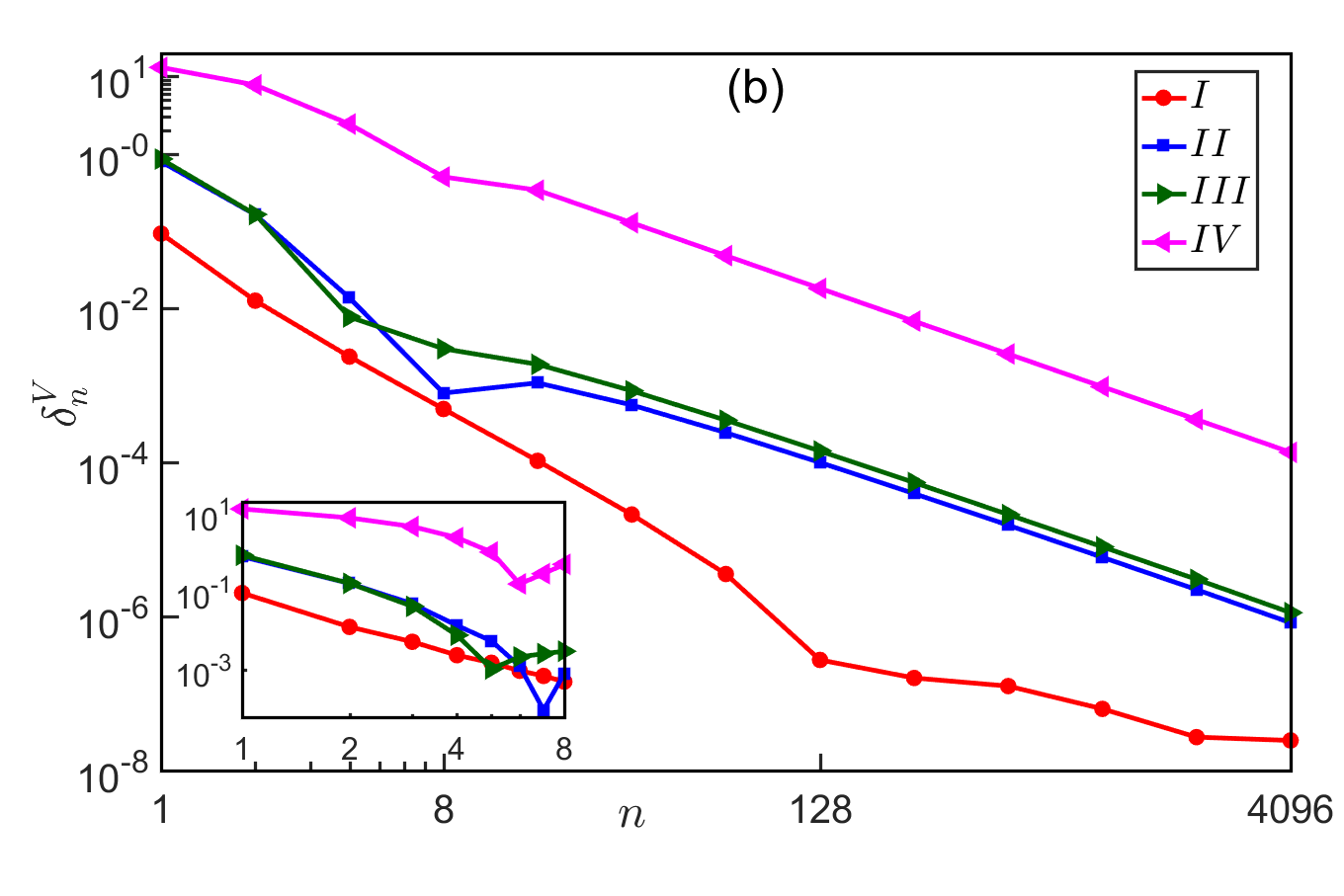,height=5cm,width=7cm,angle=0}}
\centerline{\psfig{figure=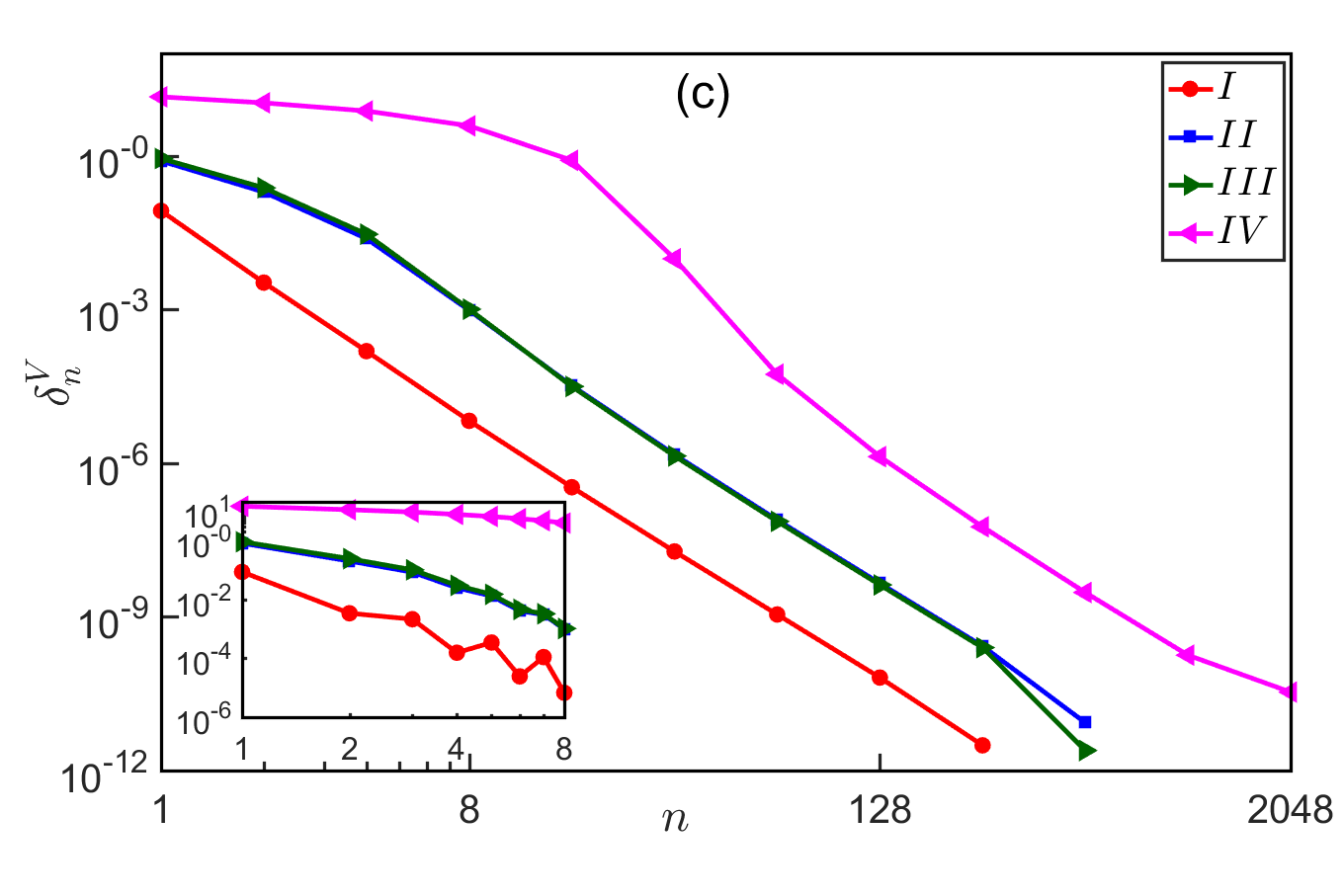,height=5cm,width=7cm,angle=0}
\psfig{figure=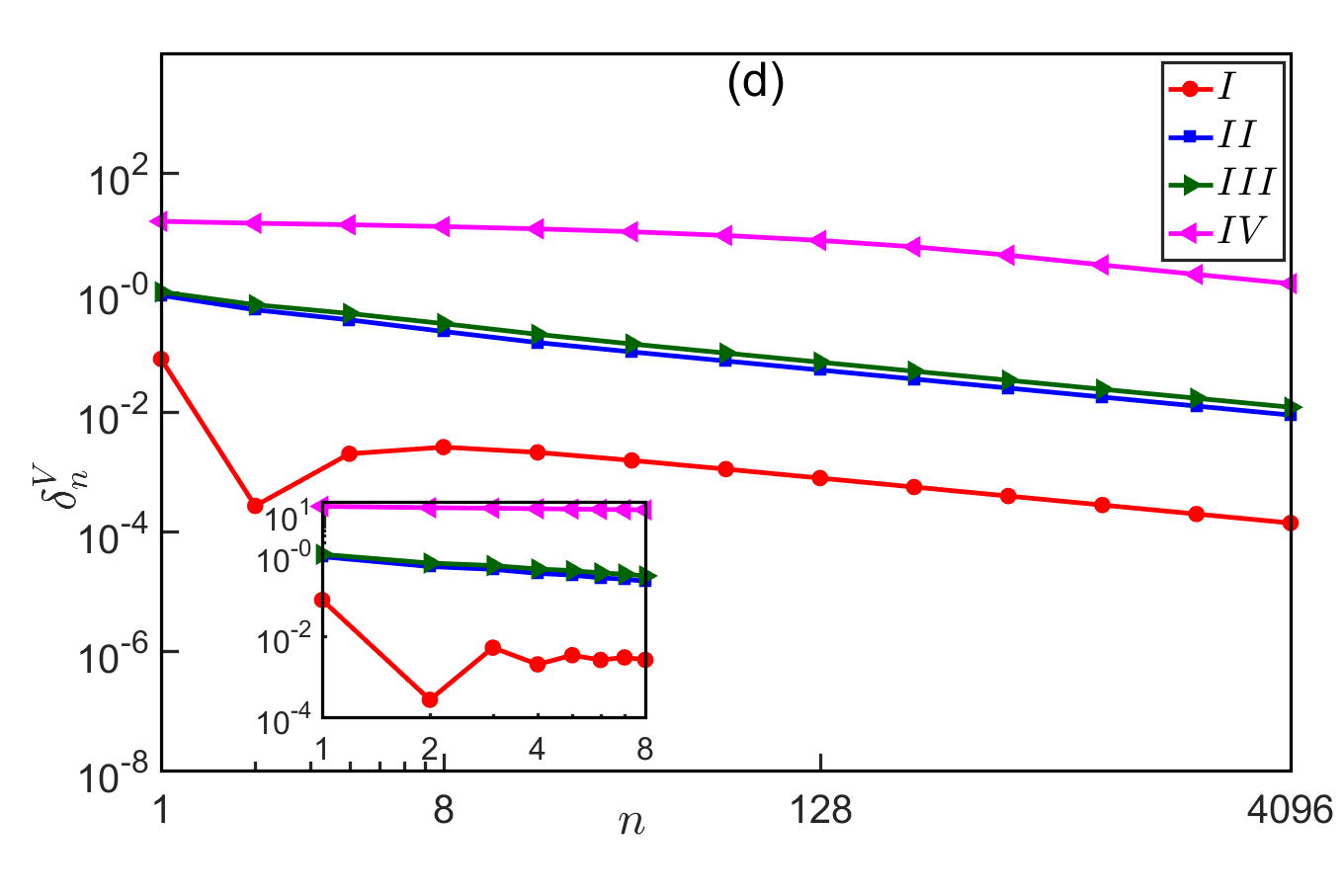,height=5cm,width=7cm,angle=0}}
\caption{Differences of the eigenvalues of \eqref{fproblem} with a potential $V$ and without a potential, i.e. $\delta_n^V:=\lambda_n^\alpha-\lambda_n^{\alpha,0}-C_V$ ($1\le n\le N=4096$) for different potentials $V(x)$ and $\alpha$:
(a) $\alpha=2$, (b) $\alpha=\sqrt{2}$, (c) $\alpha=1$, and (d) $\alpha=0.5$. }
\label{fig:eigdiffv}
\end{figure}

From Fig. \ref{fig:eigdiffv}, we can draw the following conclusion for the eigenvalues of \eqref{fproblem} with a potential $V$:
\be
\lambda_n^\alpha=\lambda_n^{\alpha,0}+C_V+O
\left(n^{-\tau_1(\alpha)}\right), \qquad n\gg1,
\ee
where
\be
\tau_1(\alpha)=\left\{\ba{ll}
\alpha &0<\alpha\le 2 \& \alpha\ne1,\\
\approx 4.5 &\alpha=1,\\
\ea\right.
\ee

\subsection{Eigenfunctions}

\begin{figure}[h!]
\centerline{
\psfig{figure=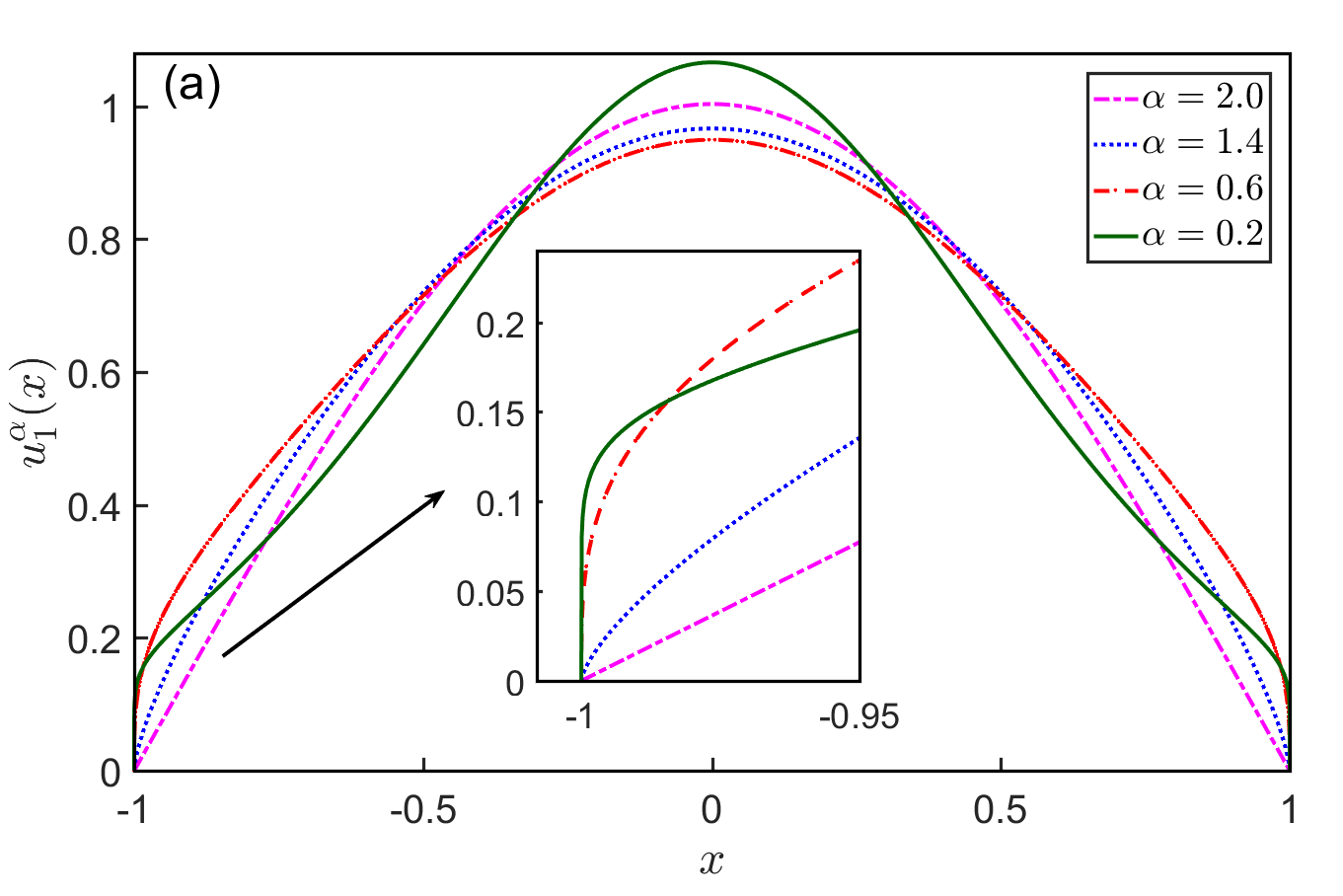,height=5cm,width=7cm,angle=0}
\psfig{figure=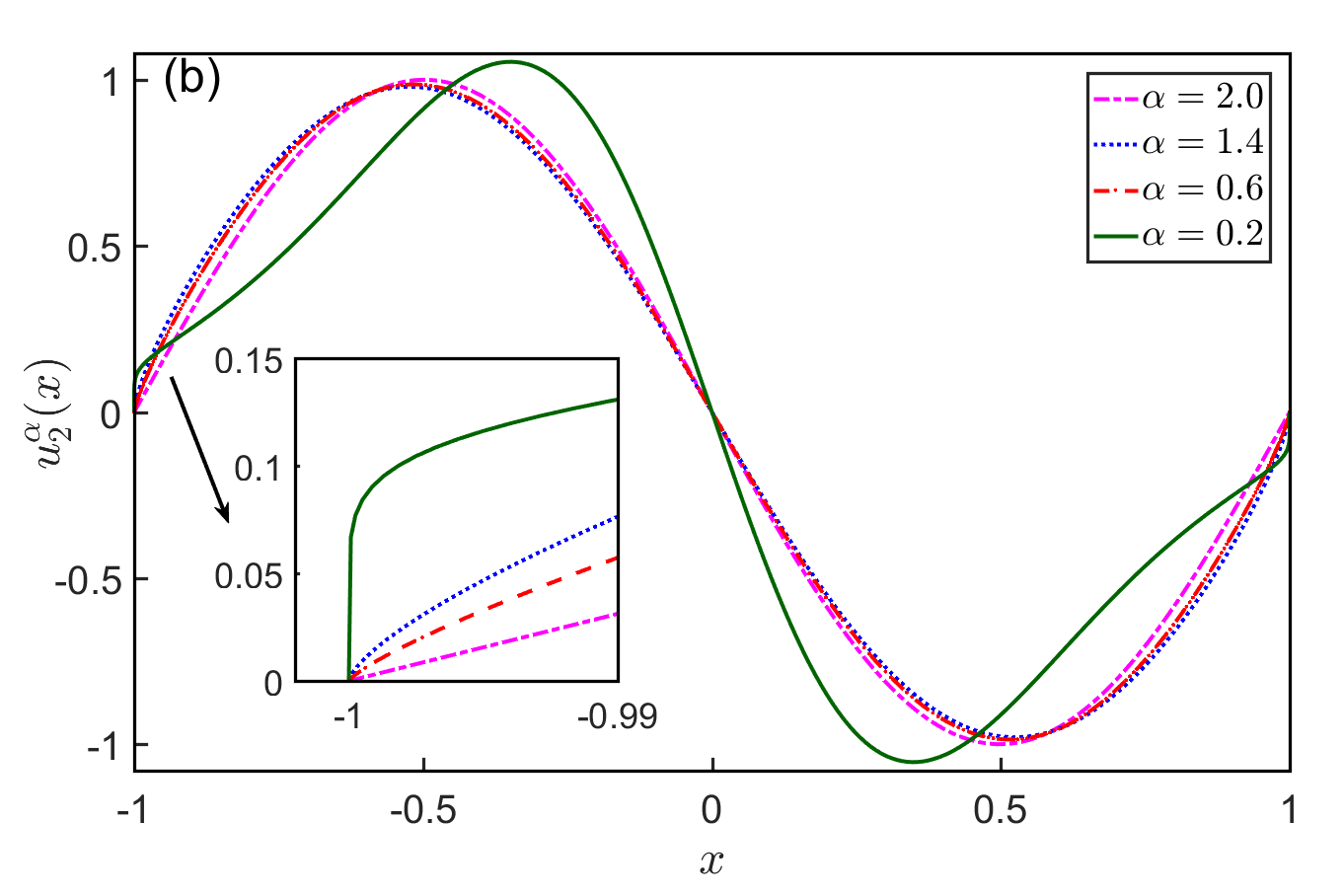,height=5cm,width=7cm,angle=0}}
\centerline{\psfig{figure=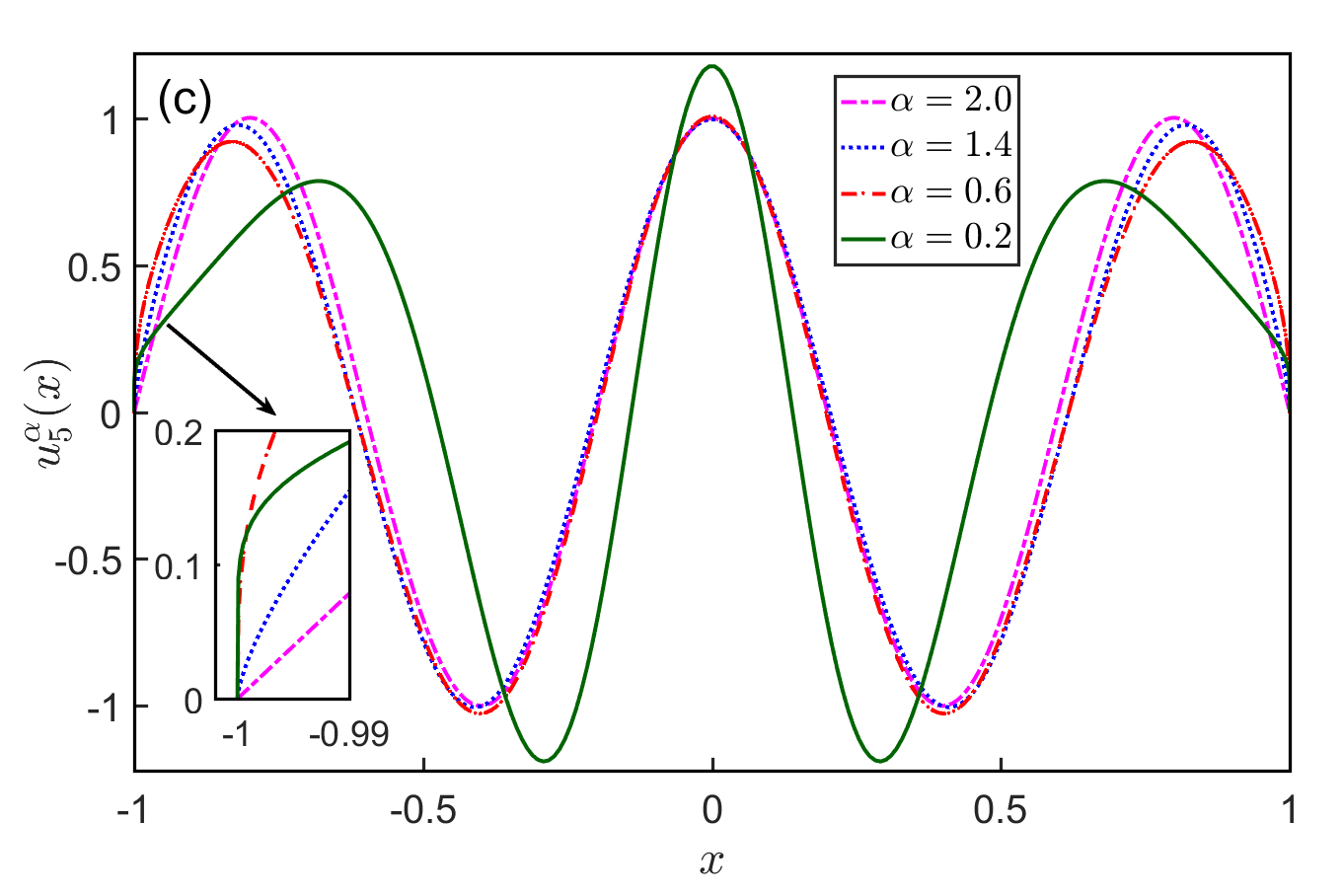,height=5cm,width=7cm,angle=0}
\psfig{figure=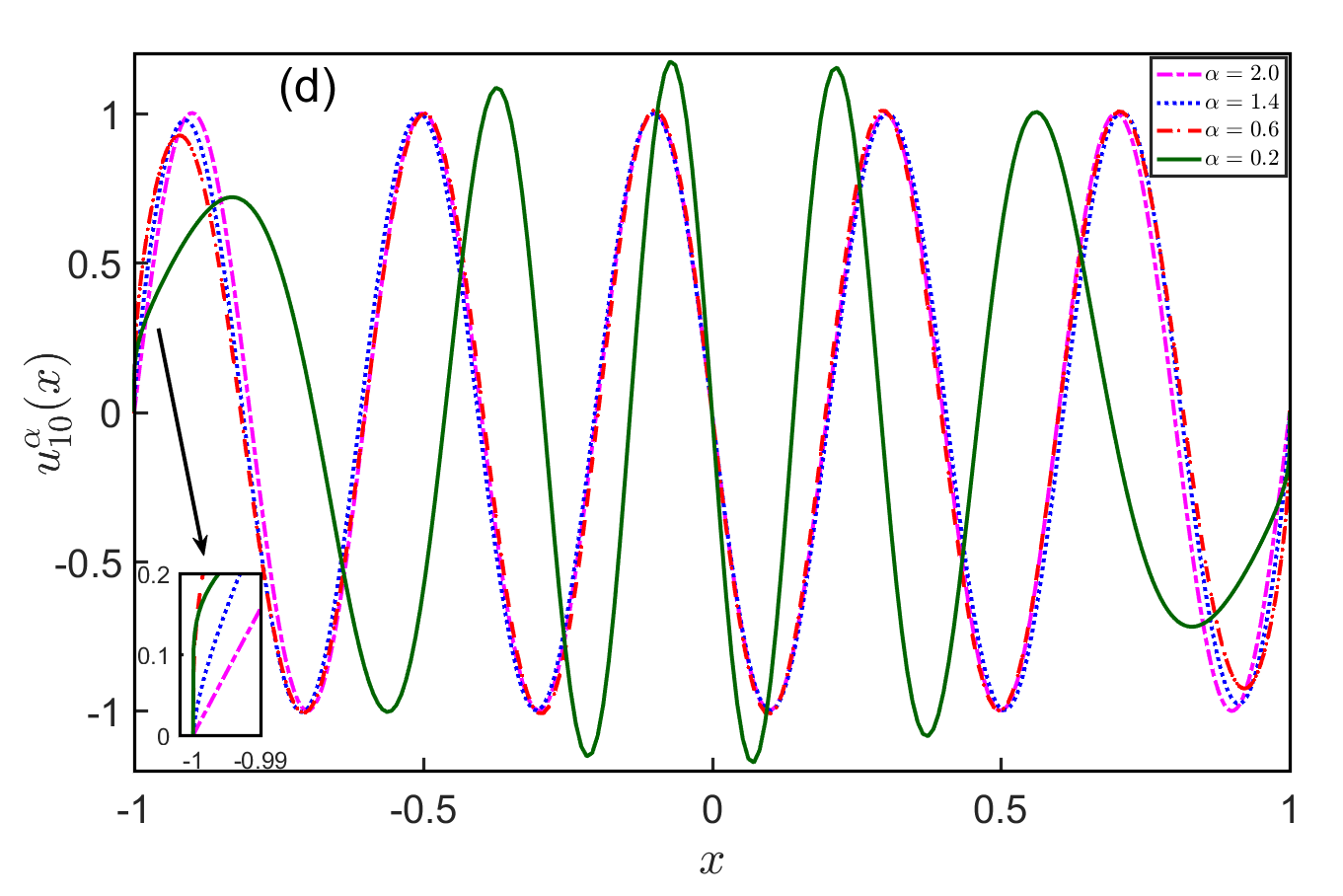,height=5cm,width=7cm,angle=0}}
\caption{Plots of different eigenfunctions of \eqref{fproblem}
with $\Omega=(-1,1)$, $V(x)=\frac{x^2}{2}$ and different $\alpha$: (a) the first eigenfunction $u_1^\alpha(x)$,
(b) the second eigenfunction $u_2^\alpha(x)$,
(c) the fifth eigenfunction $u_5^\alpha(x)$,
and (d) the tenth eigenfunction $u_{10}^\alpha(x)$.}
\label{fig:eigfwp}
\end{figure}

Figure \ref{fig:eigfwp} plots different eigenfunctions $u_n^\alpha$ of \eqref{fproblem} with $\Omega=(-1,1)$ and $V(x)=\frac{x^2}{2}$ for different $\alpha$.

From Fig. \ref{fig:eigfwp}, the singularity characteristics of the eigenfunction given in \eqref{uvx1d} is still valid for the eigenvalue problem of FSO \eqref{fproblem} with potential $V(x)$. In addition, our numerical results indicate that, when $n\to\infty$ (cf. Fig. \ref{fig:eigf}d), the eigenfunctions $u_n^\alpha$ ($0<\alpha<2$) of \eqref{fproblem} with a potential $V$ converge
to the eigenfunction $u_n^{\alpha=2}=\sin\left(\frac{n\pi (x+1)}{2}\right)$  which is the eigenfunction of \eqref{fproblem} with $\alpha=2$ and $V(x)\equiv0$.

\bigskip

Finally, based on our extensive numerical results and observations, we speculate the following observation (or conjecture) for the FSO in \eqref{fproblem} with potential:

\smallskip

{\bf Conjecture II} (Gaps and their distribution statistics of FSO in \eqref{fproblem} with potential) Assume $1<\alpha\le2$ and $V(x)\in C(\bar\Omega)$ in \eqref{fproblem}, then we have the following asymptotics of its eigenvalues:
\be \label{asy-eig-1Dw}
\lambda_n^\alpha=\left\{\ba{ll}
\left(\frac{n\pi}{b-a}\right)^\alpha-
\left(\frac{\pi}{b-a}\right)^\alpha\frac{\alpha(2-\alpha)}{4} n^{\alpha-1}+C_V +O(n^{\alpha-2}), &1<\alpha\le2,\\
\frac{n\pi}{b-a}-
\frac{\pi}{4(b-a)}+C_V +O(n^{-1}), &\alpha=1,\\
\left(\frac{n\pi}{b-a}\right)^\alpha+C_V-
\left(\frac{\pi}{b-a}\right)^\alpha\frac{\alpha(2-\alpha)}{4} n^{\alpha-1} +O(n^{-\alpha}), &0<\alpha<1,\\
\ea\right.
\qquad n\gg1,
\ee
where
\be
C_V=\frac{1}{|\Omega|}\int_\Omega V(x) dx=\frac{1}{b-a}\int_a^b V(x) dx.
\ee
In addition, we have the following asymptotics for different gaps:
\be\label{asy-gaps-1Dw}
\begin{split}
&\delta_{\rm nn}^\alpha(N)\approx\left(\frac{\pi}{b-a}\right)^\alpha\left[\alpha N^{\alpha-1}+\frac{\alpha(\alpha-1)(2+\alpha)}{4}N^{\alpha-2}
+O(N^{\alpha-3})\right],\qquad 0<\alpha\le 2,\\
&\delta_{\rm min}^\alpha(N)=\lambda_{N+1}^\alpha-\lambda_N^\alpha \approx \alpha \left(\frac{\pi}{b-a}\right)^\alpha N^{\alpha-1},  \qquad 0<\alpha<1,\\
&\delta_{\rm ave}^\alpha(N)\approx\left(\frac{\pi}{b-a}\right)^\alpha
\left\{\ba{ll}\vspace{2mm}
\left[N^{\alpha-1}+\frac{\alpha(2+\alpha)}{4}N^{\alpha-2}+O(N^{-1})\right], &1<\alpha\le2,\\ \vspace{2mm}
\left[1+\left(\frac{3}{4}-\frac{b-a}{\pi}\lambda_1^{\alpha=1}\right)N^{-1}+O(N^{-2})\right],
 &\alpha=1,\\
\left[N^{\alpha-1}-\left(\frac{b-a}{\pi}\right)^\alpha \lambda_1^\alpha N^{-1} +O(N^{\alpha-2})\right], &0<\alpha<1,\\
\ea\right.\\
&\delta_{\rm norm}^\alpha(N)\approx 1+O(N^{-2}), \qquad 0<\alpha\le 2,
%\left\{\ba{ll}\vspace{2mm}
%\alpha-\frac{\alpha(2+\alpha)}{4}N^{-1}+O(N^{-2})\le \alpha, &1<\alpha<2,\\ %\vspace{2mm}
%1+\left(\frac{b-a}{\pi}\lambda_1^{\alpha=1}-\frac{3}{4}\right)
%N^{-1}+O(N^{-2}),
% &\alpha=1,\\
%\alpha+\alpha\left(\frac{b-a}{\pi}\right)^\alpha \lambda_1^\alpha %N^{-\alpha} +O(N^{-1})\ge\alpha, &0<\alpha<1.\\
%\ea\right.
\end{split}
\qquad  N\gg1.
\ee
In addition, for the gaps distribution statistics defined
in \eqref{Pas},  we have
\be\label{gap-statw}
P_{\alpha}(s)=\delta(s-1), \qquad s\ge0, \qquad 0<\alpha\le2.
\ee

\section{Extension to directional fractional Schr\"{o}dinger operator in high dimensions}\label{sec:2dgaps}
\setcounter{equation}{0}

In this section, we extend the Jacobi spectral method ({\bf JSM}) presented in Section 2 to directional
fractional Schr\"{o}dinger operator (D-FSO) in high dimensions and
apply it to study numerically its eigenvalues and their gaps as well as gap distribution statistics.

\subsection{D-FSO  in high dimensions}

Consider the eigenvalue problem related to the directional
fractional Schr\"{o}dinger operator (D-FSO) in high dimensions:

Find $\lambda \in \mathbb{R}$ and a nonzero real-valued function $u(\bx)\ne 0$ such that
\begin{equation}\label{fproblemhd}
\begin{split}
L_{\hbox{D-FSO}}\; u(\bx)&:=\left[\mathcal{D}^{\alpha}_{\mathbf{x}}
+V(\bx)\right]u(\bx)
=\lambda\; u(\bx), \qquad \bx\in \Omega:=(-L_1,L_1)\times \ldots (-L_d,L_d)
\subset {\mathbb R}^d, \\
u(\bx)&=0, \qquad  \bx \in \Omega^c:=\mathbb{R}^d \backslash \Omega,
\end{split}
\end{equation}
where $d\ge2$, $\bx=(x_1,x_2,\ldots,x_d)^T$, $0<\alpha\le 2$, $V(\bx)\in L^2(\Omega)$ is a given real-valued function and the directional fractional Laplacian operator $\mathcal{D}^{\alpha}_{\mathbf{x}}:=\sum_{j=1}^d
(-\partial_{x_jx_j})^{\alpha/2}$ is defined via
the Fourier transform (see \cite{CS07,NPV12,LPGSGZMCMA18} and references therein) as
  \begin{equation}\label{foperator2d}
\mathcal{D}^{\alpha}_{\mathbf{x}}\;u(\bx)={\mathcal F}^{-1} \left(\sum_{j=1}^d |\xi_j|^{\alpha}(\mathcal{F}u)(\boldsymbol{\xi})\right)\qquad \bx,\boldsymbol{\xi}\in \mathbb{R}^d,
\end{equation}
with $\boldsymbol{\xi}=(\xi_1,\xi_2,\ldots,\xi_d)^T$, ${\mathcal F}$ and ${\mathcal F}^{-1}$ the Fourier transform
and inverse Fourier transform over ${\mathbb R}^d$ \cite{P99,STW11}, respectively. We remark here that the directional fractional Laplacian operator $\mathcal{D}^{\alpha}_{\mathbf{x}}$ has been widely used in the literature for different fractional PDEs, see \cite{LCTBA13,ZLLBTA14,MS16,FLJT17,LPGSGZMCMA18} and references therein.
Without loss of generality, we assume that $L_1\ge L_2\ge \ldots \ge L_d>0$.

  Again, since all eigenvalues of  \eqref{fproblemhd} are distinct (or all
spectrum are discrete and no continuous spectrum), similar to \eqref{eigenorder} for \eqref{fproblem}, we can also rank (or order) all eigenvalues of \eqref{fproblemhd} as \eqref{eigenorder}, while again
the times that an eigenvalue $\lambda$ of \eqref{fproblemhd} appears in the sequence
\eqref{eigenorder} is the same as its algebraic multiplicity. Under the order of all eigenvalues in \eqref{eigenorder} for \eqref{fproblemhd}, we define
the fraction of repeated eigenvalues of \eqref{fproblemhd} as
\begin{equation}\label{rpratio}
R^\alpha(N):=\frac{\#\{2\le n\le N \ |\ \lambda_n^\alpha=\lambda_{n-1}^\alpha\}}{N}, \qquad N=2,3,\ldots\;.
\ee
In addition, let $0<\lambda_1^{\alpha,0}<\lambda_2^{\alpha,0}<\ldots<\lambda_n^{\alpha,0}<\ldots$ be all eigenvalues of \eqref{fproblem} with $\Omega=(-1,1)$ and $V(x)\equiv 0$, and $u_n^{\alpha,0}(x)$ ($n=1,2,\ldots$) be the corresponding eigenfunctions.
Then when $V(\bx)\equiv 0$ in \eqref{fproblemhd}, all
eigenvalues of the problem \eqref{fproblemhd} can be given as
\be\label{eigen2D}
\lambda_{j_1\ldots j_d}^\alpha = \sum_{l=1}^d L_l^{-\alpha}\lambda_{j_l}^{\alpha,0}, \qquad j_1,\ldots, j_d=1,2,\ldots,
\ee
and their corresponding eigenfunctions can be given as
\be\label{eigenv2D}
u_{j_1\ldots j_d}^\alpha(\bx)=\Pi_{l=1}^d u_{j_l}^{\alpha,0}(x_l/L_l),
\qquad \bx\in \bar \Omega, \qquad j_1,\ldots, j_d=1,2,\ldots\;.
\ee
The above results immediately imply that the fundamental gap of
\eqref{fproblemhd} with $V(\bx)\equiv 0$ can be obtained as
\be
\delta_{\rm fg}(\alpha)=L_1^{-\alpha}\lambda_{2}^{\alpha,0}+
\sum_{l=2}^d L_l^{-\alpha}\lambda_{1}^{\alpha,0}-\sum_{l=1}^d L_l^{-\alpha}\lambda_{1}^{\alpha,0}=
L_1^{-\alpha}\left(\lambda_{2}^{\alpha,0}-\lambda_{1}^{\alpha,0}\right)
\ge \frac{\lambda_{2}^{\alpha,0}-\lambda_{1}^{\alpha,0}}{(D/2)^\alpha},
\ee
where $D$ is the diameter of $\Omega$.

The JSM presented in Section 2 can be easily
extended to solve the eigenvalue problem \eqref{fproblemhd} by tensor product \cite{MS16}. The details are omitted here for brevity.

\begin{figure}[t!]
\centerline{
\psfig{figure=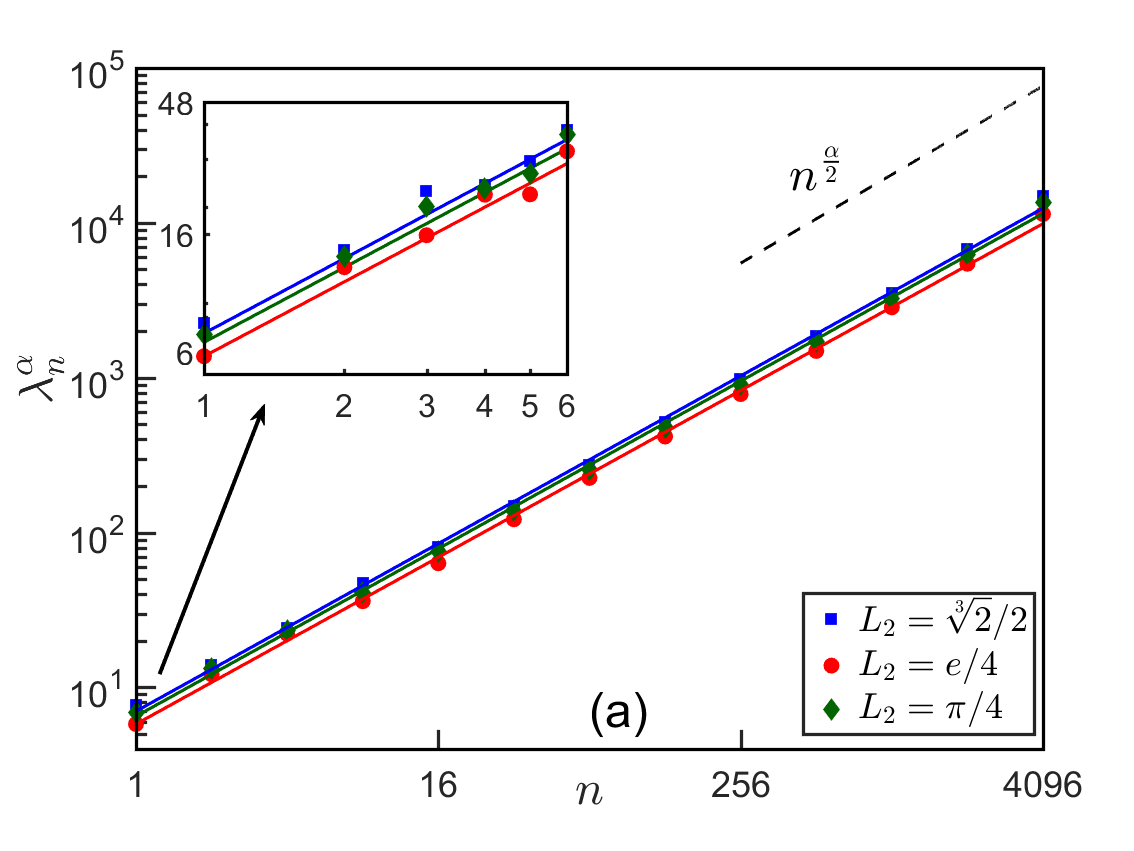,height=5cm,width=7cm,angle=0}
\psfig{figure=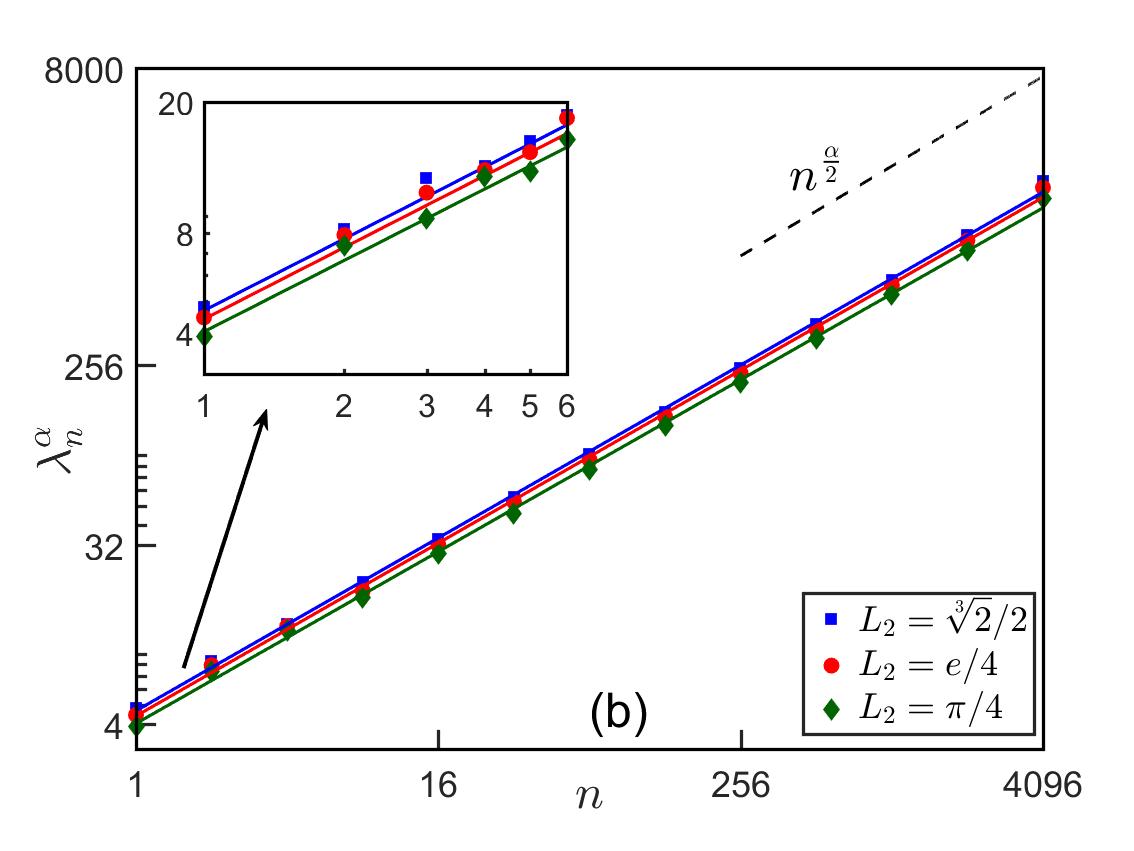,height=5cm,width=7cm,angle=0}}
\centerline{\psfig{figure=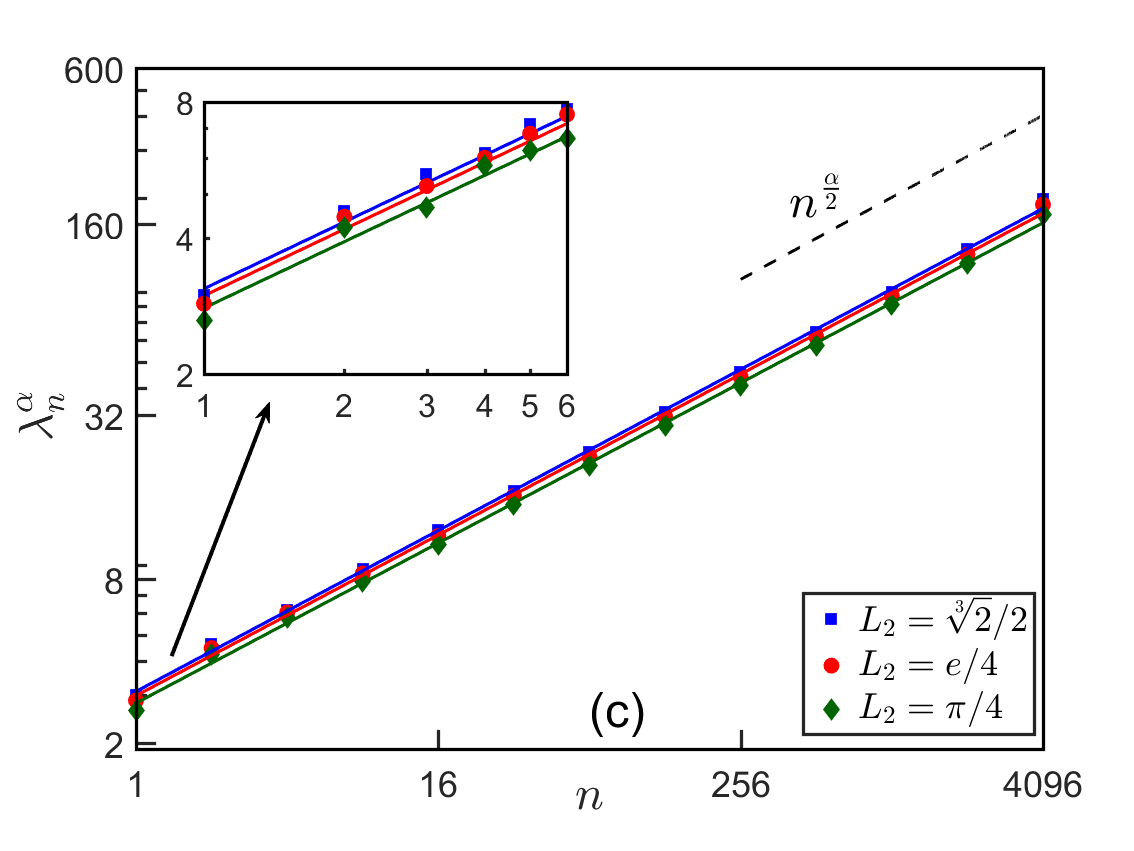,height=5cm,width=7cm,angle=0}
\psfig{figure=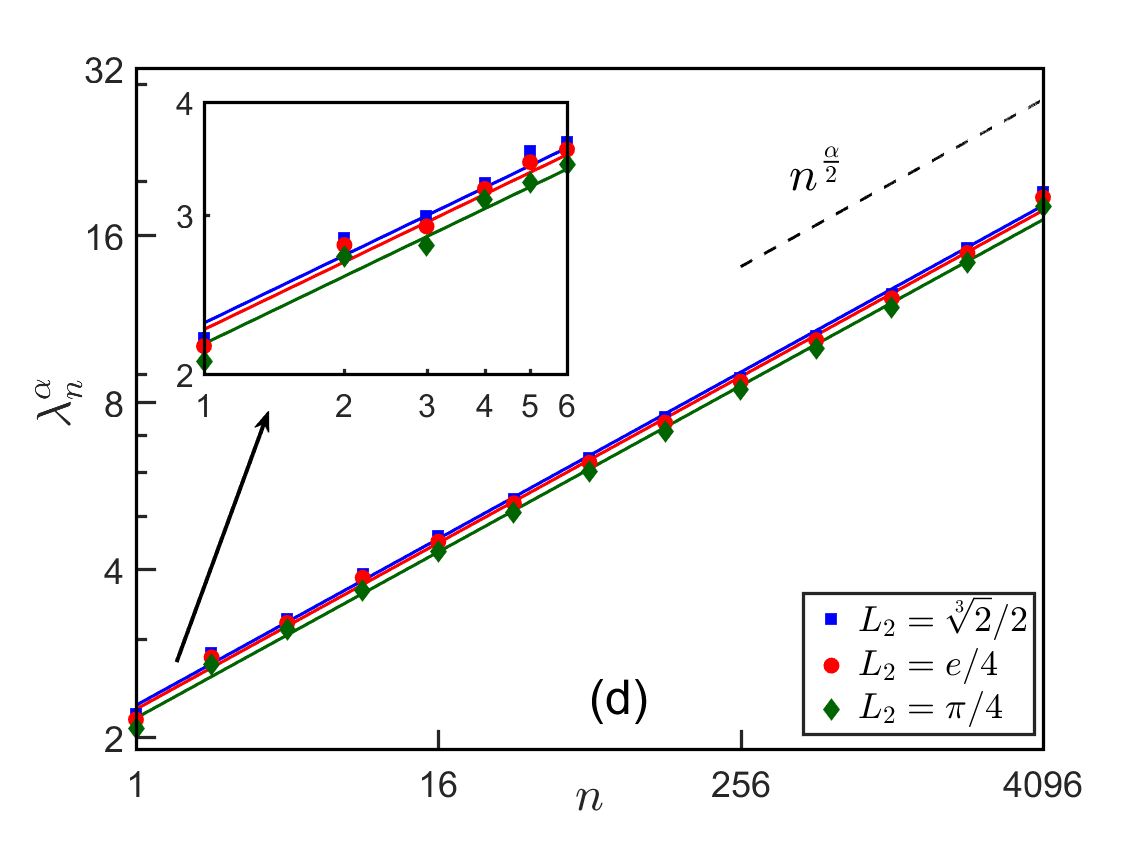,height=5cm,width=7cm,angle=0}}
\caption{Eigenvalues of \eqref{fproblemhd} with
$d=2$, $L_1=1$, $V(\bx)\equiv 0$ and  different $L_2$ and $\alpha$
(symbols denote numerical results and solids lines are from fitting formula $C_2^\alpha n^{\alpha/2}$ when $n\gg1$):
(a) $\alpha=1.9$, (b) $\alpha=1.5$,
(c) $\alpha=1.0$, and (d) $\alpha=0.5$. }
\label{fig:eig1dnphd}
\end{figure}
%%%%%%%%%%%%%%%%
%%%%%%%%%%%
\begin{figure}[t!]
\centerline{
\psfig{figure=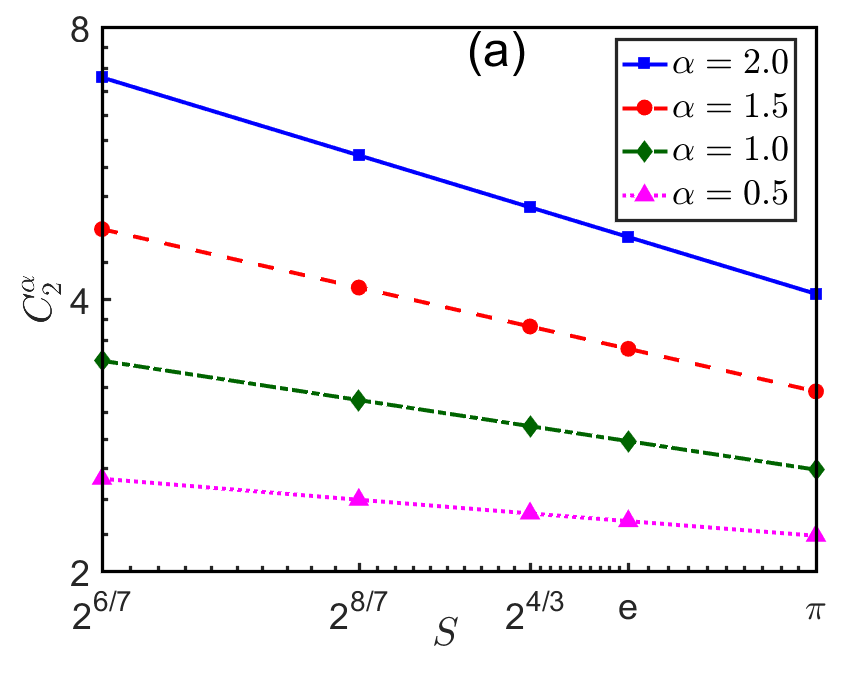,height=5cm,width=7cm,angle=0}
\psfig{figure=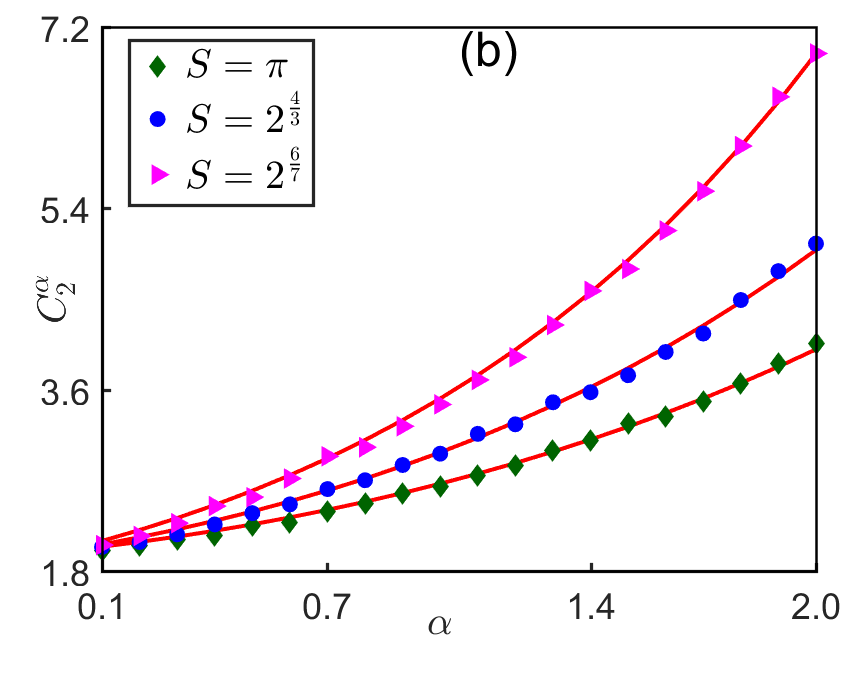,height=5cm,width=7cm,angle=0}}
\caption{Numerical results of $C_2^{\alpha}$ (symbols denote numerical results and solids lines are from fitting formula) for different areas
$S=|\Omega|=4L_2$ and $\alpha$: (a) plots of $C_2^{\alpha}$ as a function
of $S$ for different $\alpha$, and (b) plots of $C_2^{\alpha}$ as a function
of $\alpha$ for different $S$.  }
\label{fig:eig2dfit}
\end{figure}

\begin{figure}[t!]
\centerline{
\psfig{figure=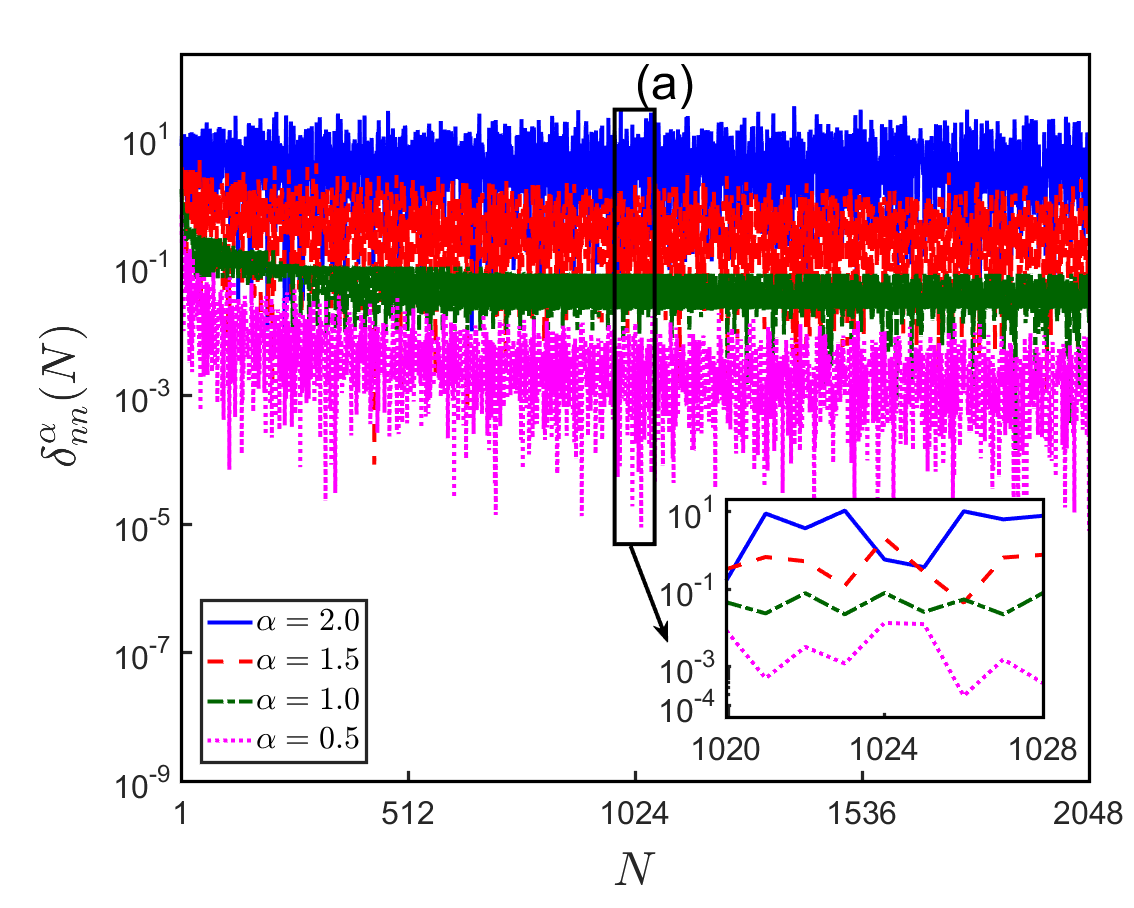,height=5cm,width=7cm,angle=0}
\psfig{figure=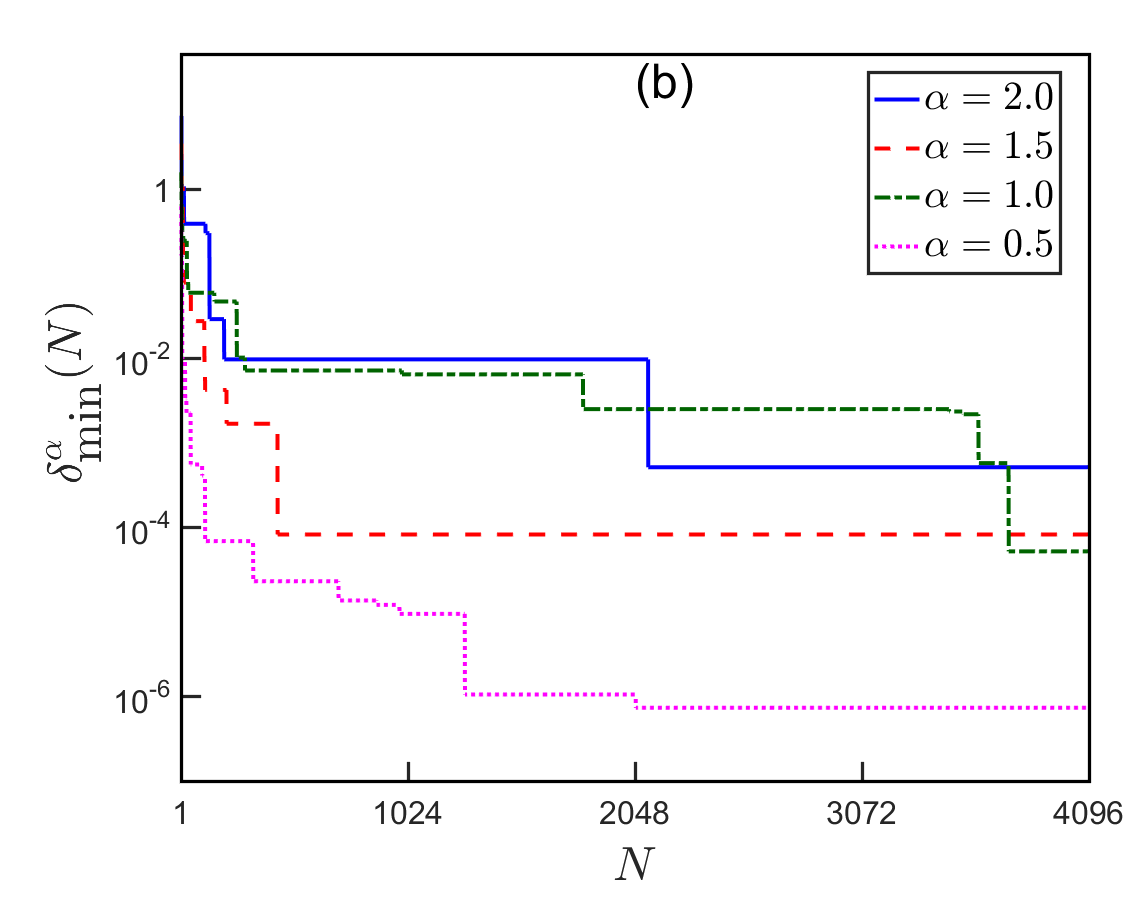,height=5cm,width=7cm,angle=0}}
\centerline{\psfig{figure=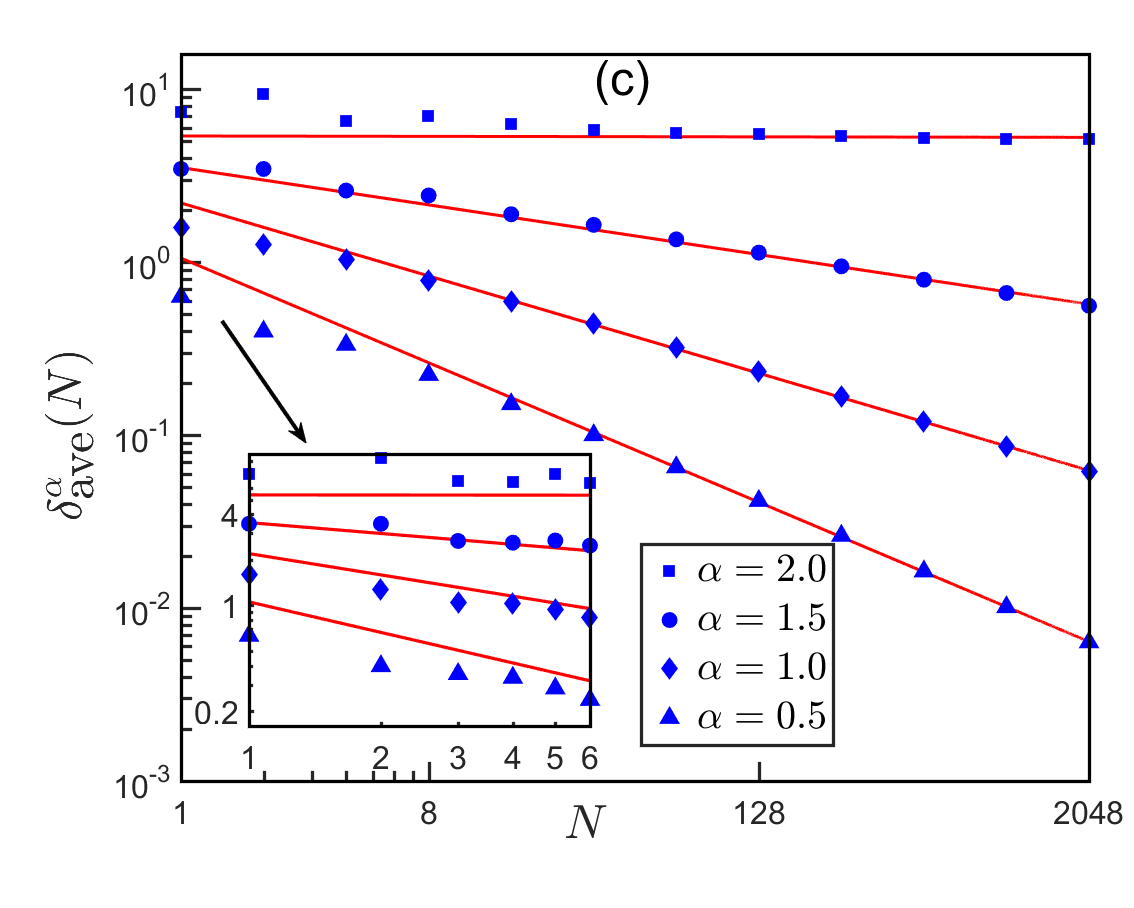,height=5cm,width=7cm,angle=0}
\psfig{figure=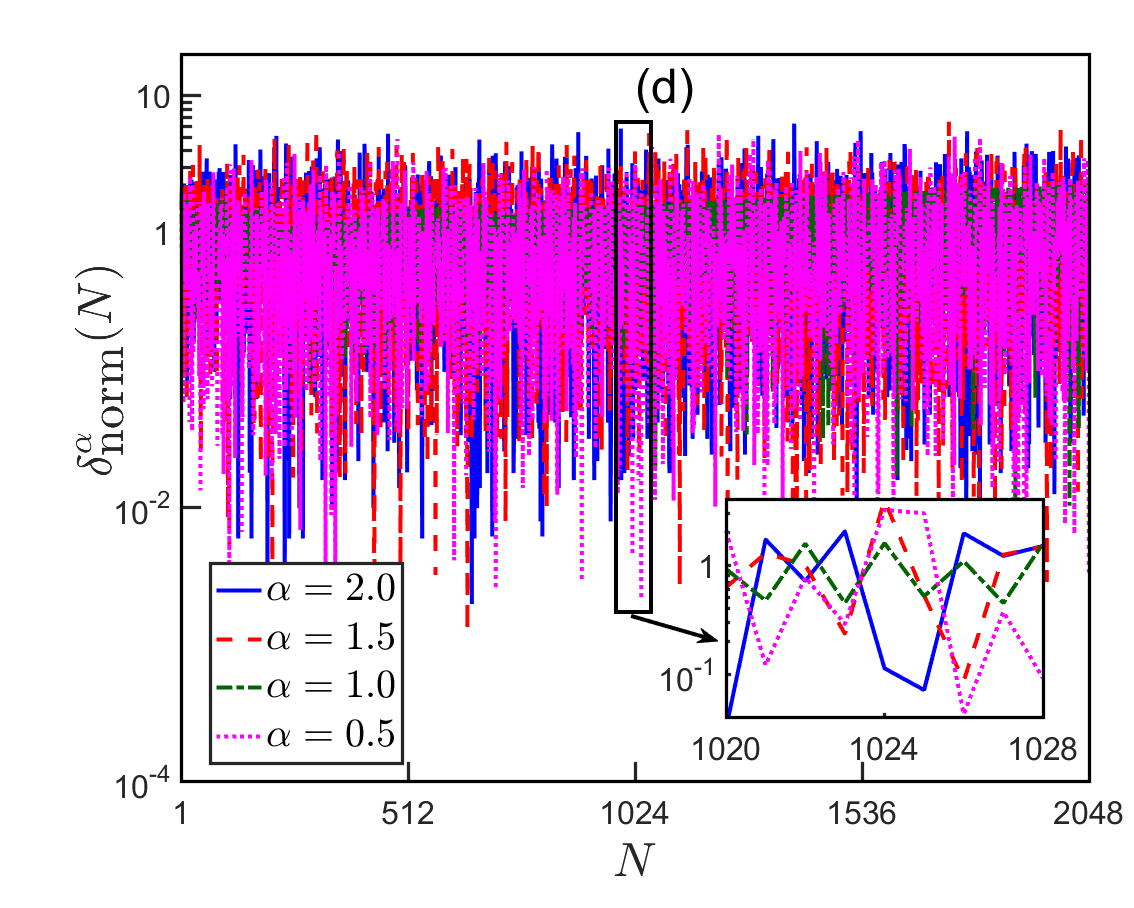,height=5cm,width=7cm,angle=0}}
\caption{Different eigenvalue gaps of \eqref{fproblemhd}
with $d=2$, $L_1=1$, $V(\bx)\equiv 0$, $L_2=\frac{\sqrt[3]{2}}{2}$ and different
$\alpha$ for:
(a) the nearest neighbour gaps $\delta_{\rm nn}^\alpha(N)$, (b)
the minimum gaps $\delta_{\textrm{min}}^\alpha(N)$,
 (c) the average gaps  $\delta_{\textrm{ave}}^\alpha(N)$ (symbols denote numerical results and solids lines are from fitting formula when $N\gg1$),
 and (d) the normalized gaps
$\delta_{\rm norm}^\alpha(N)$. }
\label{fig:gapshd}
\end{figure}

\begin{figure}[t!]
\centerline{
\psfig{figure=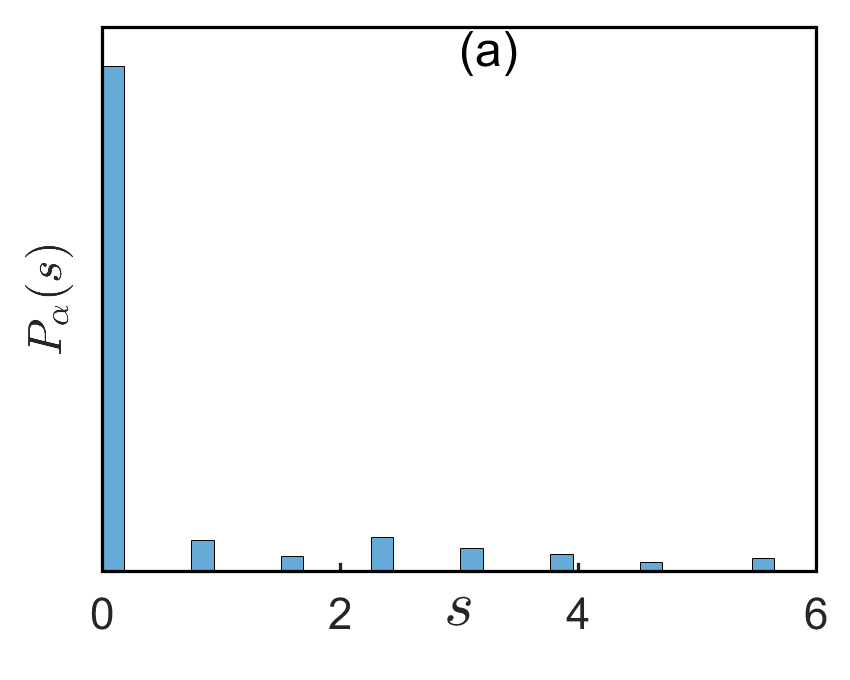,height=4cm,width=5cm,angle=0}
\psfig{figure=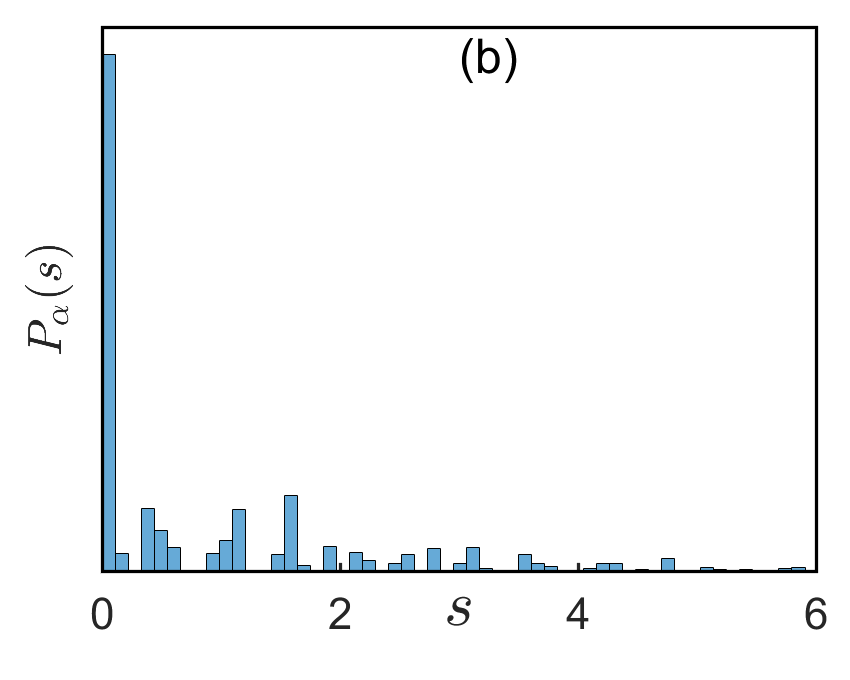,height=4cm,width=5cm,angle=0}
\psfig{figure=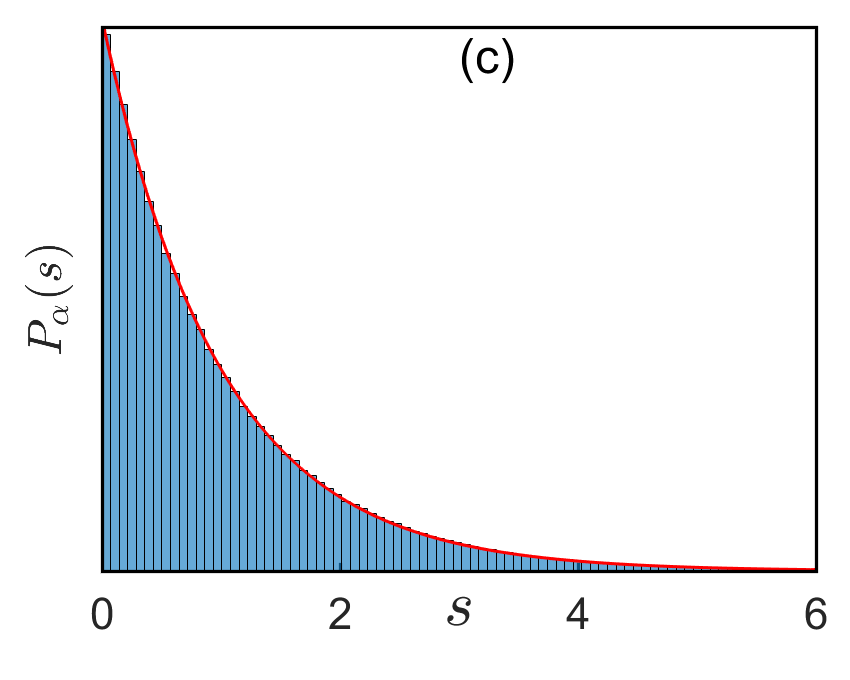,height=4cm,width=5cm,angle=0}}
\centerline{
\psfig{figure=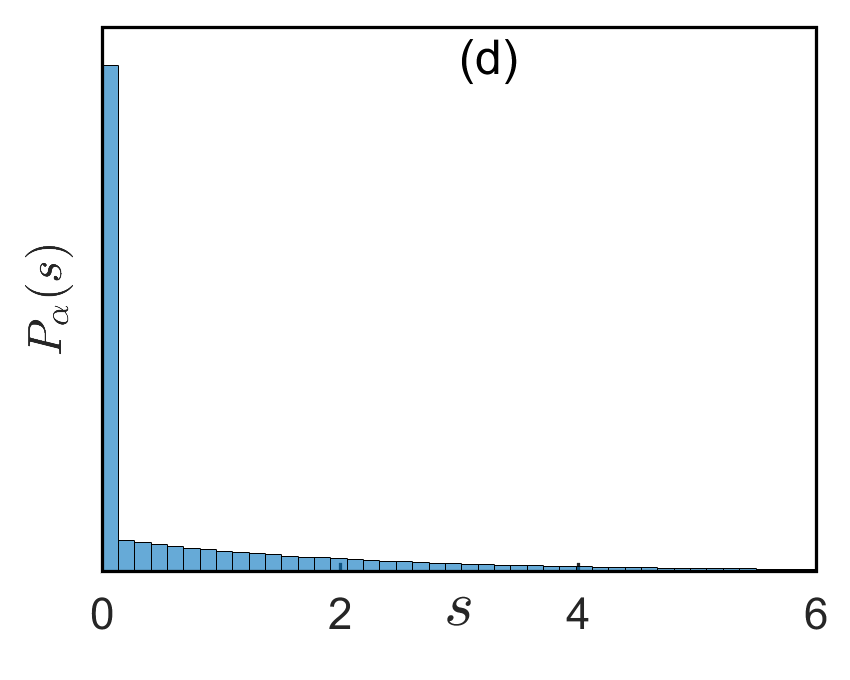,height=4cm,width=5cm,angle=0}
\psfig{figure=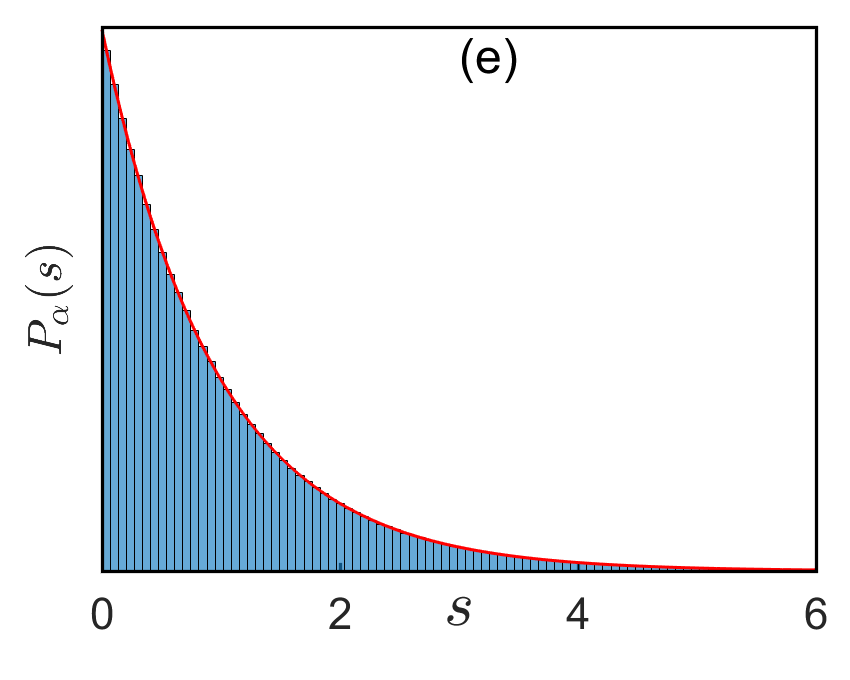,height=4cm,width=5cm,angle=0}
\psfig{figure=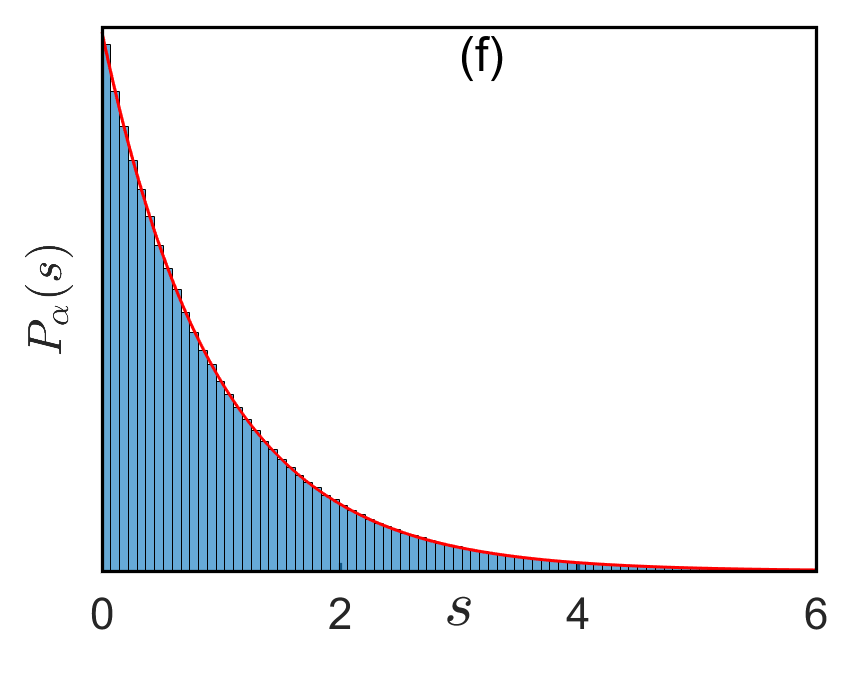,height=4cm,width=5cm,angle=0}}
\centerline{
\psfig{figure=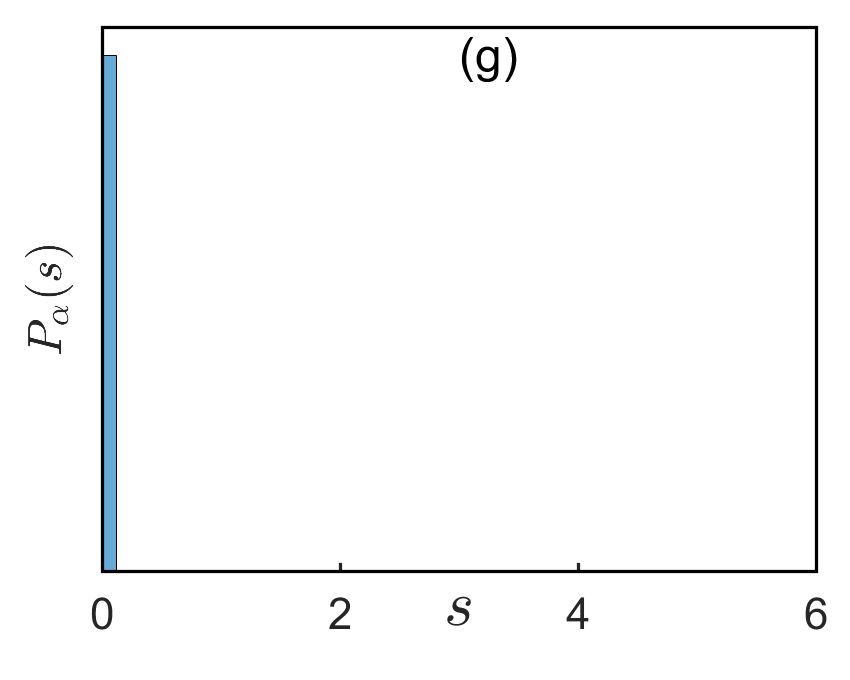,height=4cm,width=5cm,angle=0}
\psfig{figure=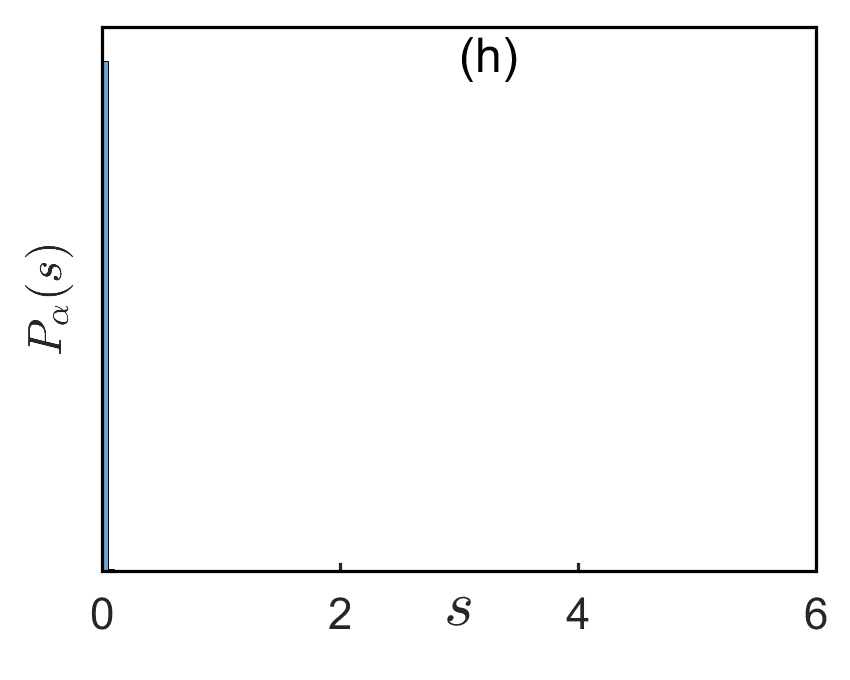,height=4cm,width=5cm,angle=0}
\psfig{figure=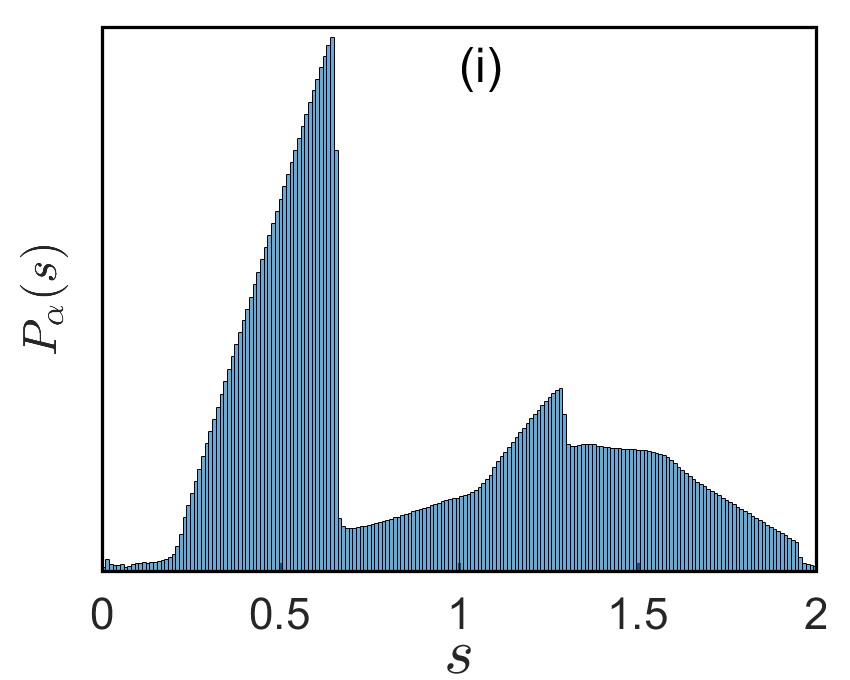,height=4cm,width=5cm,angle=0}}
\centerline{
\psfig{figure=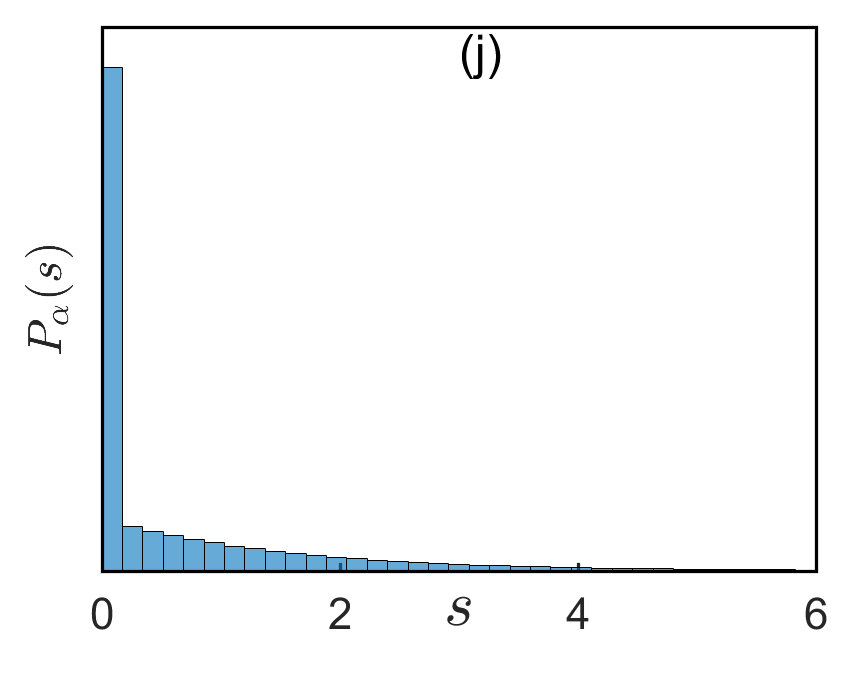,height=4cm,width=5cm,angle=0}
\psfig{figure=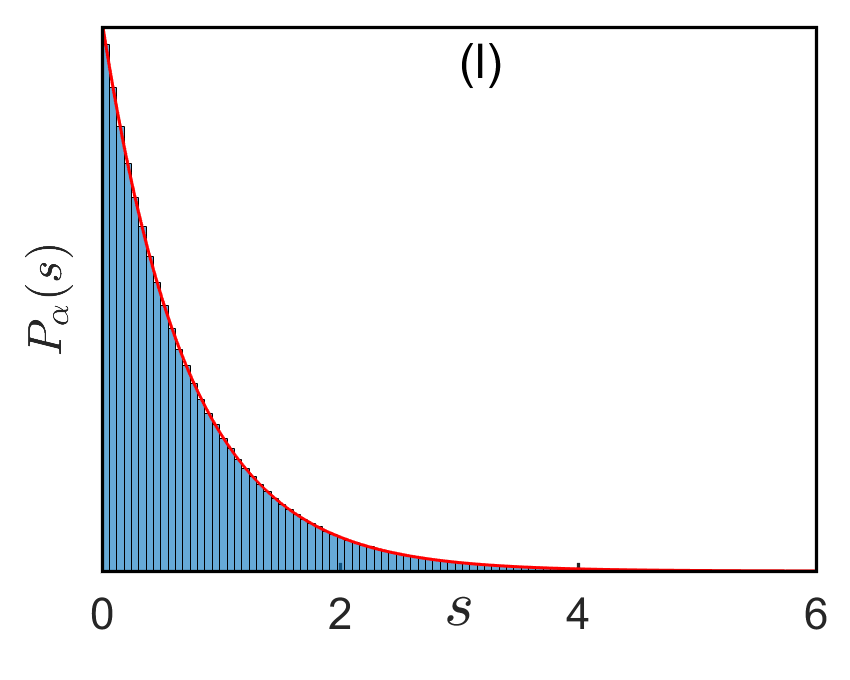,height=4cm,width=5cm,angle=0}
\psfig{figure=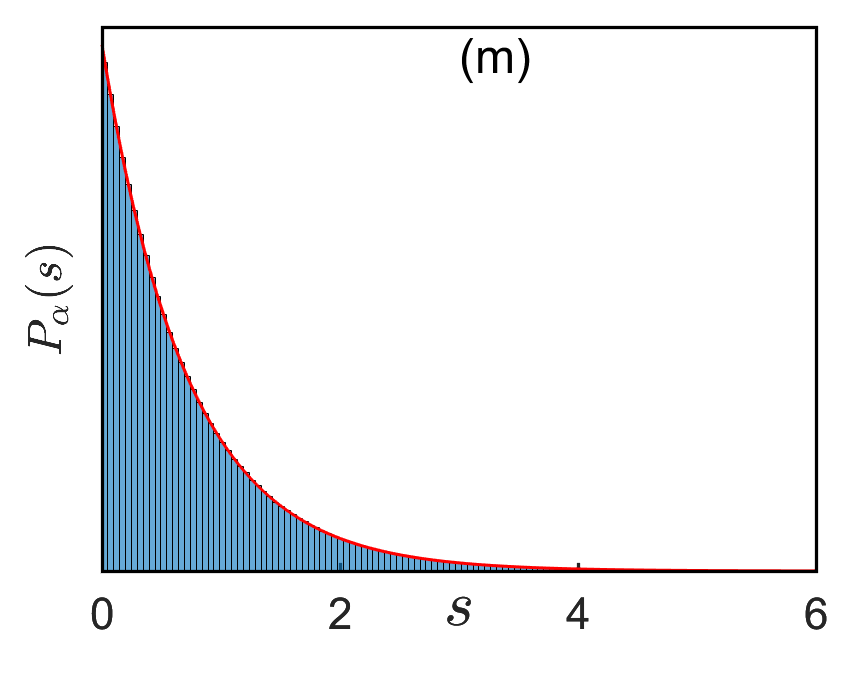,height=4cm,width=5cm,angle=0}}
\caption{The histogram of the normalized gaps  $\{\delta_{\rm norm}^\alpha(n)\ |\ 1\le n\le N=4000000\}$ of \eqref{fproblemhd}
with $d=2$, $L_1=1$ and $V(\bx)\equiv 0$ for different $0<\alpha\le 2$ and $0<L_2\le 1$: (a) $\alpha=2.0$ and $L_2=1$,
(b) $\alpha=2.0$ and $L_2=2/3$, (c) $\alpha=2.0$ and $L_2=\frac{\sqrt[3]{2}}{2}$;
(d) $\alpha=1.5$ and $L_2=1$,
(e) $\alpha=1.5$ and $L_2=2/3$, (f) $\alpha=1.5$ and $L_2=\frac{\sqrt[3]{2}}{2}$;
(g) $\alpha=1.0$ and $L_2=1$,
(h) $\alpha=1.0$ and $L_2=2/3$, (i) $\alpha=1.0$ and $L_2=\frac{\sqrt[3]{2}}{2}$;
(j) $\alpha=0.5$ and $L_2=1$,
(l) $\alpha=0.5$ and $L_2=2/3$, (m) $\alpha=0.5$ and $L_2=\frac{\sqrt[3]{2}}{2}$.  Solid lines are fitting curves for the gaps distribution statistics $P_\alpha(s)$.}
\label{fig:gaps2da}
\end{figure}

\subsection{Numerical results in two dimensions (2D) without potential}

We take $d=2$, $L_1=1$  and $V(\bx)\equiv 0$ in \eqref{fproblemhd}.
In this case, noting \eqref{eigen2D} and \eqref{eigenv2D} with $d=2$,
instead of using the JSM in 2D to compute eigenvalues and their corresponding eigenfunctions of \eqref{fproblemhd}, a simple and more efficient and accurate way is
to first use the JSM  in 1D to compute the eigenvalues and their corresponding eigenfunctions of
\eqref{fproblem} with $\Omega=(-1,1)$ and $V(x)\equiv 0$, and then
to get the eigenvalues and their corresponding eigenfunctions of \eqref{fproblemhd} with $d=2$ and $V(\bx)\equiv 0$
via \eqref{eigen2D} and \eqref{eigenv2D} with $d=2$.

In our computations, we first use the JSM  in 1D with $M=8192$ to compute numerically the eigenvalues of \eqref{fproblem} with $\Omega=(-1,1)$ and $V(x)\equiv 0$. Then we use the first $N=4096$ computed eigenvalues to
get the eigenvalues of \eqref{fproblemhd} with $d=2$ and $V(\bx)\equiv 0$
via \eqref{eigen2D} with $d=2$ and then rank (or order)
the total $4096\times 4096$ eigenvalues of \eqref{fproblemhd} as \eqref{eigenorder}. Finally, we take (up to) the first $N=4000000$ eigenvalues to compute the gaps and their distribution statistics.

%   We remark here that, since we are mainly interested in
%gaps and their distribution statistics, we choose $L_2$ an irrational number
%such that all eigenvalues of \eqref{fproblemhd} are distinct.

Figure \ref{fig:eig1dnphd} displays eigenvalues (in increasing order)
of \eqref{fproblemhd} for different $L_2$ and $\alpha$, which suggests
that $\lambda_n^\alpha\sim n^{\alpha/2}$ when $n\gg1$ for $0<\alpha\le 2$.
Then we fit numerically $\lambda_n^\alpha$ when $n\gg1$ by $C_2^\alpha n^{\alpha/2}$. Figure \ref{fig:eig2dfit} displays the fitting results
of $C_2^\alpha$ with respect to the area $S=4L_2$ of $\Omega$ and $\alpha$,
which suggests that
\be\label{c2alpw}
C_2^\alpha=\frac{4}{2+\alpha} \left(\frac{4\pi}{S}\right)^{\alpha/2} ,
\qquad 0<\alpha\le 2, \qquad S=4L_2>0.
\ee
These results immediately suggest that
\be\label{eigasp2d}
\lambda_n^\alpha = \frac{4}{2+\alpha}\left(\frac{4\pi}{S}\right)^{\alpha/2} \; n^{\alpha/2}+o(n^{\alpha/2}), \qquad n\gg1.
\ee
Specifically, when $\alpha=2$, our numerical results suggest that
\be\label{eigasp2d2}
\lambda_n^{\alpha=2} = \frac{4\pi}{S}\left[n+C_1n^{1/2}+O(1)\right], \qquad n\gg1,
\ee
where $C_1\approx 0.5943$ from our numerical results. In fact,
\eqref{eigasp2d2} can be regarded as an improved Weyl law when $\alpha=2$ \cite{Weyl}, and \eqref{eigasp2d} can be regarded as an extension of the Weyl law for $\alpha=2$ \cite{Weyl} to $0<\alpha\le2$, and we call \eqref{eigasp2d} as the generalized Weyl law on
the asymptotics of the eigenvalues of
the D-FSO in 2D.

In fact, combining \eqref{eigasp2d} and \eqref{avegap}, we
can obtain the asymptotic of the average gaps of the D-FSO
in \eqref{fproblemhd} as
\bea\label{ave2dasyw}
\delta_{\rm ave}^\alpha(N)&=&\frac{\lambda_{N+1}^\alpha-\lambda_1^\alpha}{N}\nonumber\\
&=&\frac{1}{N}\left[\frac{4}{2+\alpha}\left(\frac{4\pi}{S}\right)^{\alpha/2} \; (N+1)^{\alpha/2}+o((N+1)^{\alpha/2})-\lambda_1^\alpha\right]\nonumber\\
&=&\frac{4}{2+\alpha}\left(\frac{4\pi}{S}\right)^{\alpha/2} N^{(\alpha-2)/2}+o(N^{(\alpha-2)/2})\nonumber\\
&=&O(N^{(\alpha-2)/2}), \qquad N\gg1,
\eea
which immediately implies that, when $\alpha=2$, $\delta_{\rm ave}^\alpha(N)
\sim 1$ (i.e. almost a constant) when $N\gg1$, and respectively, when $0<\alpha<2$,
$\delta_{\rm ave}^\alpha(N) \sim N^{(\alpha-2)/2}$ (decrease with respect to $N$) when $N\gg1$.

In addition, Figure \ref{fig:gapshd} plots different eigenvalue gaps of \eqref{fproblemhd}
with $d=2$, $L_1=1$, $V(\bx)\equiv 0$, $L_2=\sqrt[3]{2}/2$ and different
$\alpha$. Figure \ref{fig:gaps2da} displays the histogram of the normalized  gaps $\{\delta_{\rm norm}^\alpha(n)\ | \ 1\le n\le N=4000000\}$ for different
$\alpha$ and  $L_2$. Figure \ref{fig:gaps2Ra} plots $1-R^\alpha(N)$ vs $N$ ($N\gg1$) for different $\alpha$ and $L_2$.

\begin{figure}[t!]
\centerline{
\psfig{figure=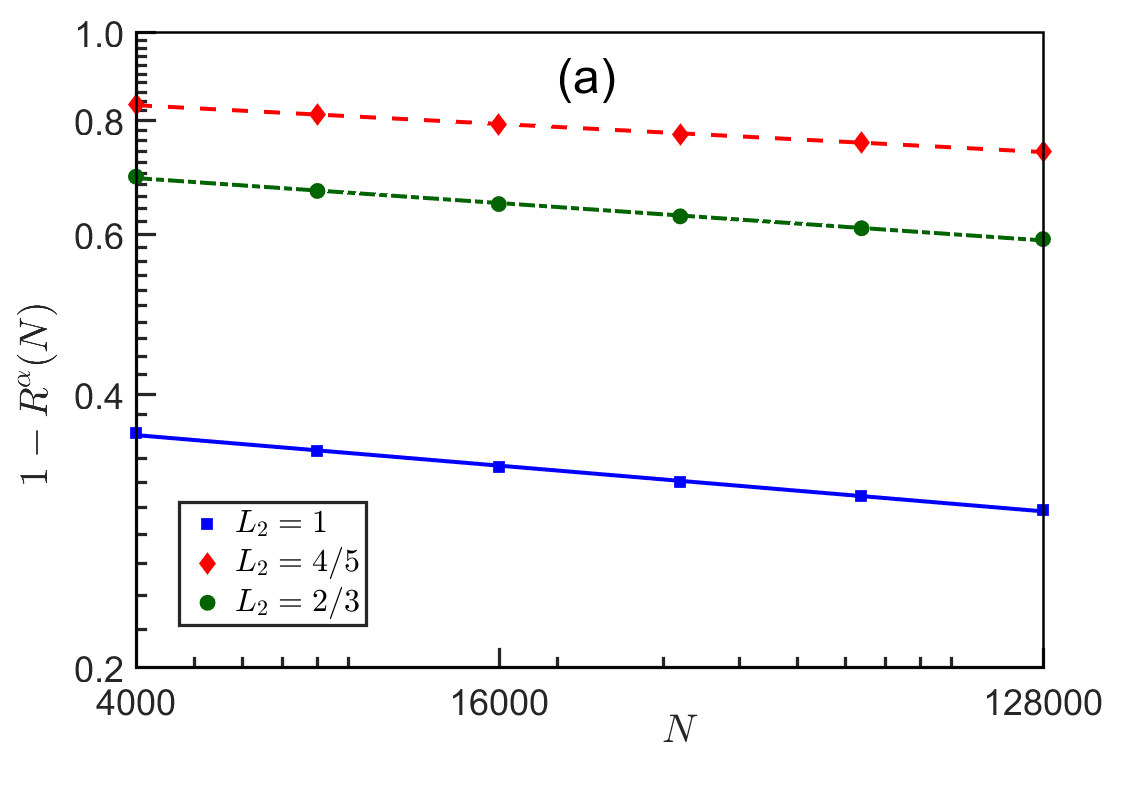,height=4cm,width=5cm,angle=0}
\psfig{figure=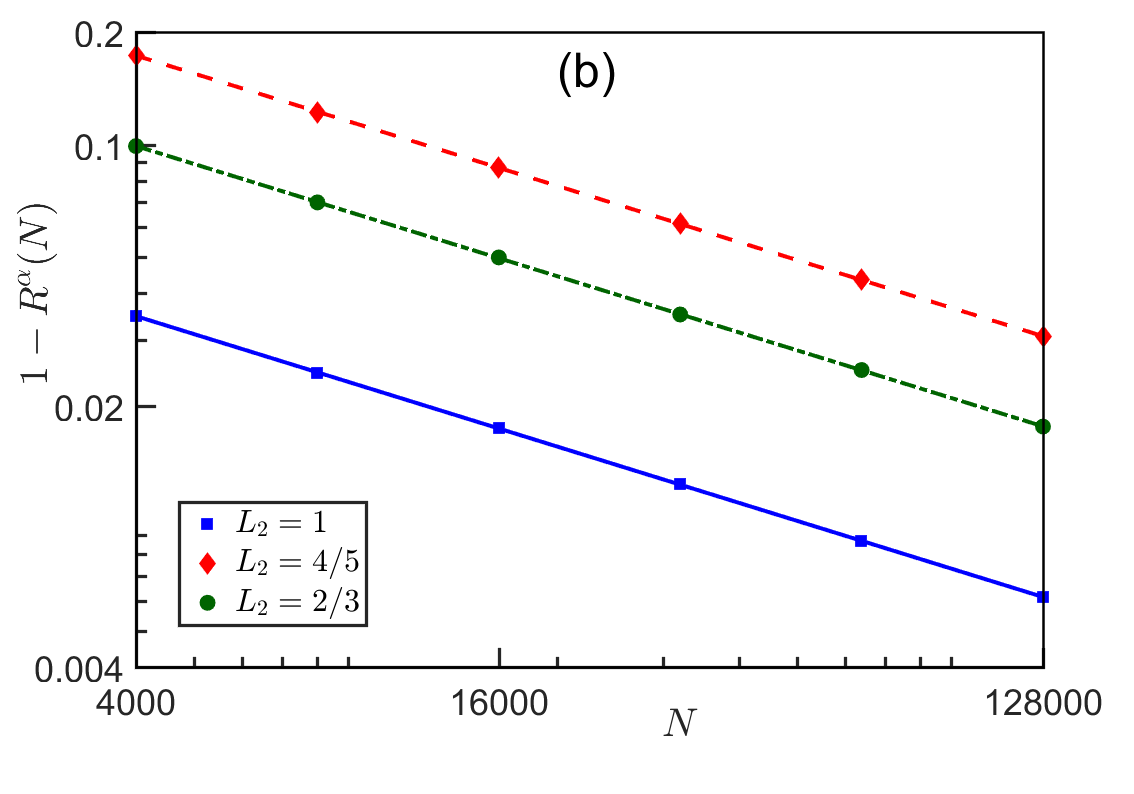,height=4cm,width=5cm,angle=0}
\psfig{figure=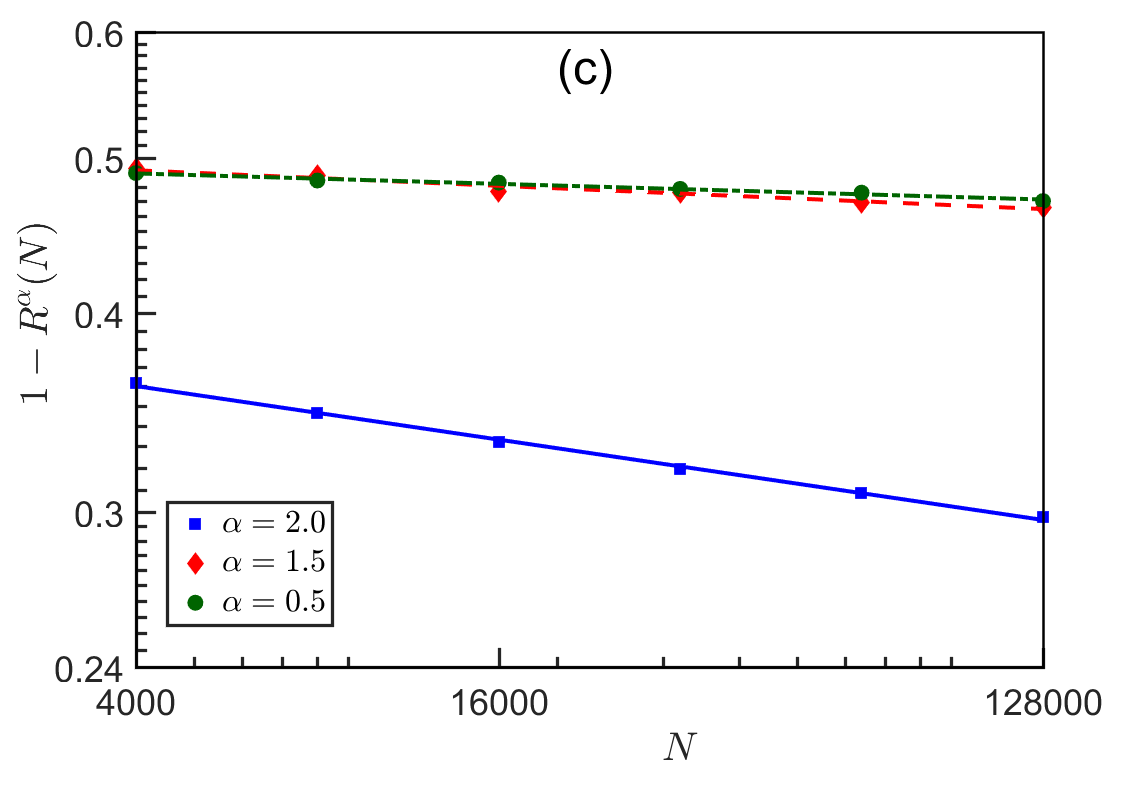,height=4cm,width=5cm,angle=0}
}
\caption{Plots of $1-R^\alpha(N)$ vs $N$ ($N\gg1$)
for different $\alpha$ and $L_2$: (a) $\alpha=2$ for different $L_2\in {\mathbb Q}$;
(b) $\alpha=1$ for different $L_2\in {\mathbb Q}$; and (c) $L_2=1$ for different $0<\alpha\le 2$.
 }
\label{fig:gaps2Ra}
\end{figure}

From Figs. \ref{fig:gapshd}-\ref{fig:gaps2Ra}, we can draw the following conclusions:

(i) The minimum gaps $\delta_{\rm min}(N)\to 0$ when $N\to+\infty$
(cf. Fig. \ref{fig:gapshd}b); and the average gaps $\delta_{\rm ave}(N)\sim 1$ when $N\gg1$
for $\alpha=2$, and respectively, $\delta_{\rm ave}(N)\sim N^{(\alpha-2)/2} $ when $N\gg1$ for $0<\alpha<2$ (cf. Fig. \ref{fig:gapshd}c),
which confirm the asymptotic results in \eqref{ave2dasyw}.

(ii) When $L_2=1$ and $0<\alpha\le 2$ or $\alpha=2$ and $L_2\in {\mathbb Q}$ or $\alpha=1$ and $L_2\in {\mathbb Q}$, the gaps distribution statistics $P_\alpha(s)=\delta(s)$
(cf. Fig. \ref{fig:gaps2da}a,b,d,g,h,j and Fig. \ref{fig:gaps2Ra}).
In these cases, $R^\alpha(N)\to 1$ when $N\to\infty$ (cf. Fig. \ref{fig:gaps2Ra}) and our numerical results suggest the following asymptotics: $R^\alpha(N)=1-N^{-\tau_2(L_2)}$ when $\alpha=2$ for different $L_2\in {\mathbb Q}$ (cf. Fig. \ref{fig:gaps2Ra}a); $R^\alpha(N)=1-N^{-1/2}$ when $\alpha=1$ for different $L_2\in {\mathbb Q}$ (cf. Fig. \ref{fig:gaps2Ra}b); and $R^\alpha(N)=1-N^{-\tau_3(\alpha)}$ when $L_2=1$ for different $0<\alpha\le2$ (cf. Fig. \ref{fig:gaps2Ra}c). In addition, Figure \ref{fig:fitls} plots
$\tau_2(L_2)$ and $\tau_3(\alpha)$ based on our numerical results.

\begin{figure}[t!]
\centerline{
\psfig{figure=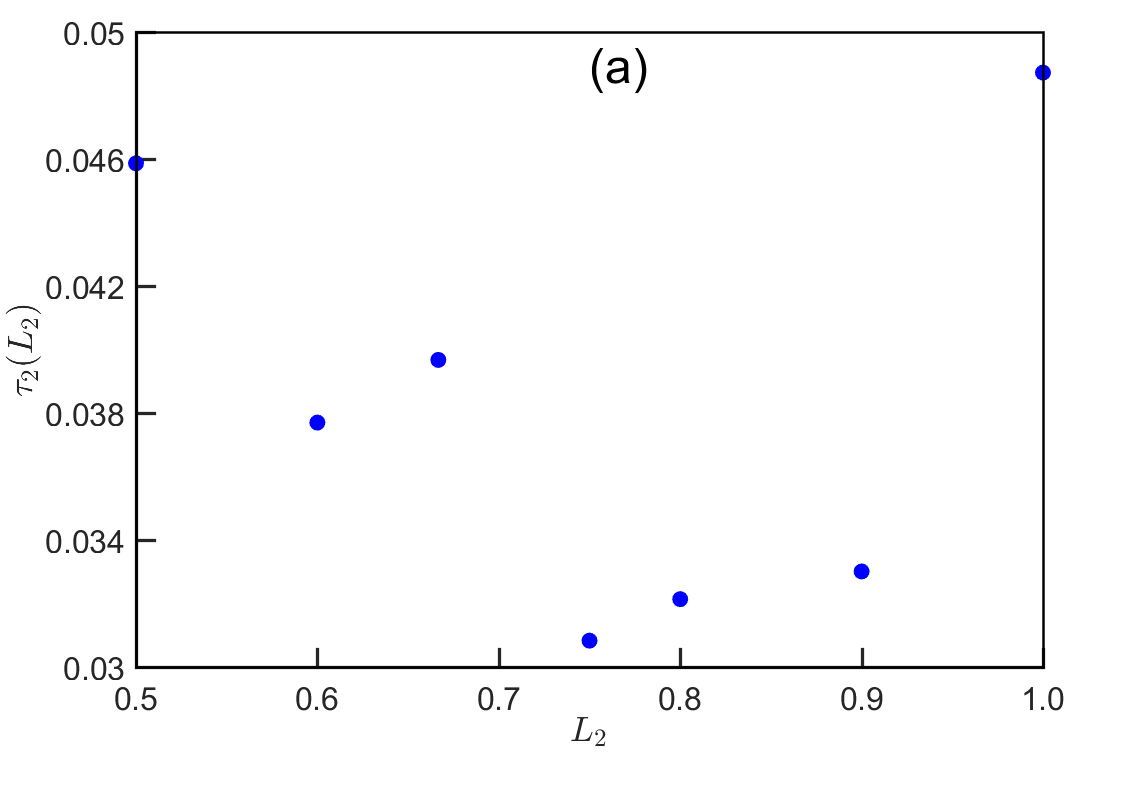,height=4cm,width=7cm,angle=0}
\psfig{figure=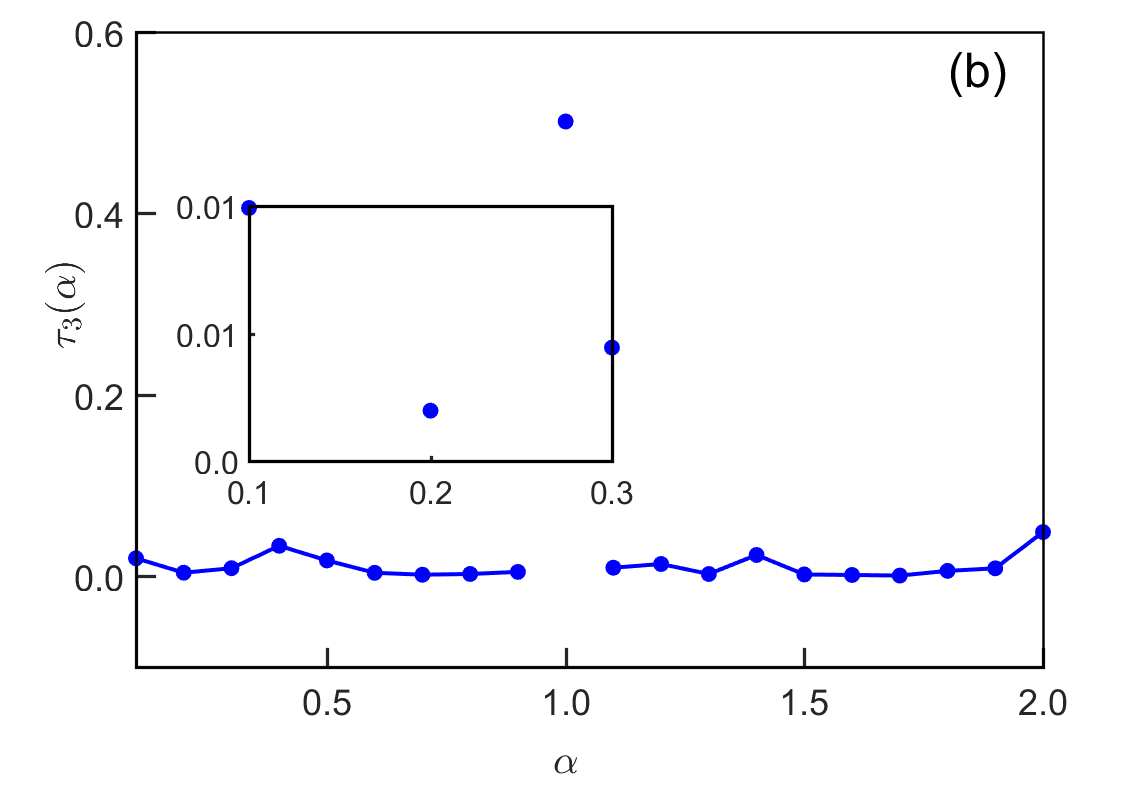,height=4cm,width=7cm,angle=0}
}
\caption{Fitting results of $\tau_2(L_2)$ for different
$L_2\in {\mathbb Q}$ (left) and
$\tau_3(\alpha)$ for different $\alpha$ (right). }
\label{fig:fitls}
\end{figure}

(iii) When $L_2\notin {\mathbb Q}$ and $0<\alpha<1$ or $1<\alpha\le 2$,  $P_\alpha(s)$ can be well approximated
by a Poisson distribution (cf. Fig. \ref{fig:gaps2da}c,e,f,l,m), i.e.
\be
P_\alpha(s)=\tau(\alpha) e^{-\tau(\alpha)\,s}, \qquad
s\ge0.
\ee
In addition, Figure \ref{fittau} plots $\tau(\alpha)$, which suggests that
\be
\tau(\alpha)\approx
\left\{\ba{ll}
1,  &1< \alpha \le 2,\\
1.057 \alpha^{-0.385},  &0< \alpha < 1.\\
\ea\right.
\ee

\begin{figure}[t!]
\centerline{
\psfig{figure=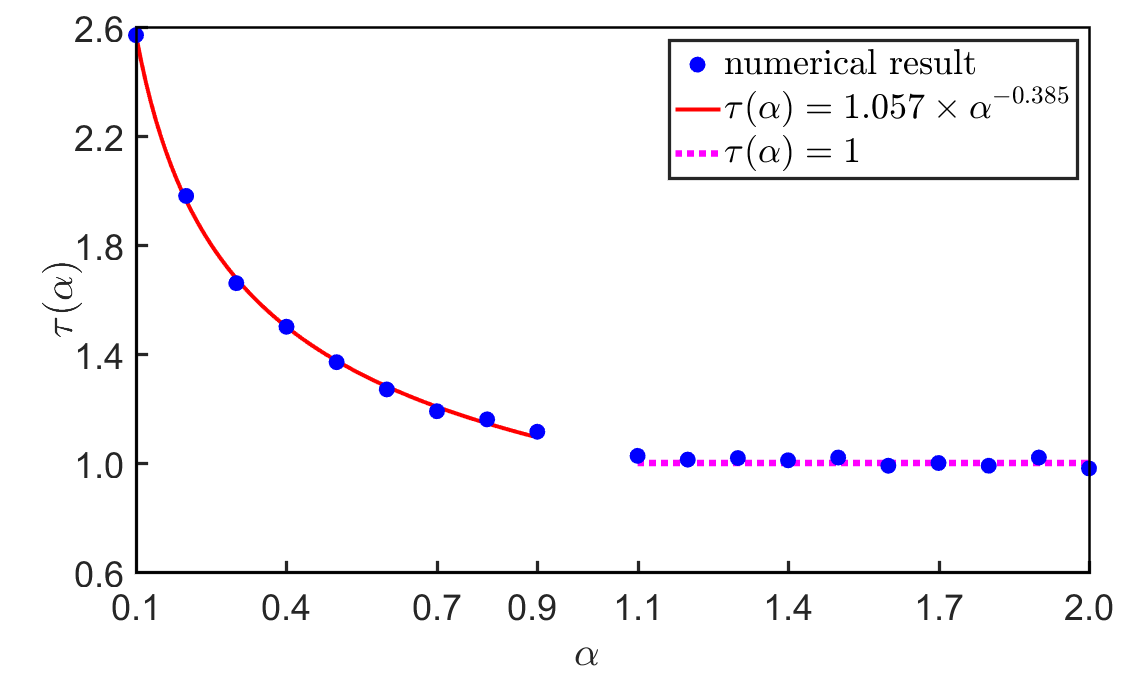,height=6cm,width=10cm,angle=0}
}
\caption{Fitting results of $\tau(\alpha)$ for different $\alpha$. }
\label{fittau}
\end{figure}

(iv) When $\alpha=1$ and $L_2\notin {\mathbb Q}$, $P_\alpha(s)$ can be well approximated by a bimodal distribution \cite{SW08} (cf. Fig. \ref{fig:gaps2da}i).

(v) The classification of the gaps distribution statistics $P_\alpha(s)$ for different $0<\alpha\le 2$ and $L_1>0$ and $L_2>0$
is summarized in Table \ref{eigsd}.

\begin{table}[t!]
\centering
\begin{tabular}{ |c |c|c|c|} \hline
& $L_2/L_1=1$          &  $1\ne L_2/L_1 \in {\mathbb Q}$          &  $1\ne L_2/L_1\notin {\mathbb Q}  $      \\ \hline
$\alpha=2$    &$\delta(s)$    &  $\delta(s)$ &  Poisson   \\
$1<\alpha<2$    &$\delta(s)$    &  Poisson &  Poisson  \\
$\alpha=1$    &$\delta(s)$    &  $\delta(s)$ &  Bimodal distribution \\
$0<\alpha<1$    &$\delta(s)$    & Poisson&  Poisson \\
\hline
\end{tabular}
\caption{Summary of the eigenvalue gap distribution statistics of \eqref{fproblemhd} with $d=2$ and $V(\bx)\equiv 0$
 for different $0<\alpha\le 2$ and $L_1>0$ and $L_2>0$.}
\label{eigsd}
\end{table}

\subsection{Numerical results in 2D with potential}

Here we use the JSM in 2D to compute numerically the eigenvalues
and their corresponding eigenfunctions of \eqref{fproblemhd} with $d=2$ and
a non-zero potential $V(x,y)$. In our computations, we choose the total DOF $M=144 \times 144$, i.e.
with DOFs $M_1=144$ and $M_2=144$ in
$x_1$ and $x_2$ directions, respectively. With the $M$ eigenvalues computed,
we only use $M/4$ (or even less) numerical eigenvalues to compute gaps and their distribution statistics. We take $L_1=1$ and $V(x,y) = \frac{x^2+y^2}{2}$ in \eqref{fproblemhd}.

%Case II. ~$V(x,y)=4 (x^2+y^2)$; % \frac{(x-0.1)^2}{2}$;

%Case III. $V(x,y)=4(x^2+y^2)+\sin(\frac{\pi}{2} x)+\sin(\frac{\pi}{2} y)$;   %x^2+\sin(x)$;

%Case IV. $V(x,y)=50(x^2+y^2)+\sin(2\pi x)+\sin(2\pi y)$.

Figure \ref{fig2d:eigdiffv} plots different eigenvalue gaps of \eqref{fproblemhd} with  $L_2=\sqrt[3]{2}/2$ for different
$\alpha$, and Figure \ref{fig2d:eigdiffvdv}  displays the histogram of the normalized  gaps $\{\delta_{\rm norm}^\alpha(n)\ | \ 1\le n\le N=4096\}$ for different $\alpha$ and  $L_2$.

We also carry out numerical simulations on eigenvalues and their different gaps as well as their distribution statistics of \eqref{fproblemhd} in 2D with different other potentials. Our numerical results suggest that the asymptotic
behavior of the eigenvalue $\lambda_n^\alpha$ in \eqref{eigasp2d} and
\eqref{eigasp2d2} are still valid when \eqref{fproblemhd} is with a potential
$V(\bx)\in C(\Omega)$. In addition,  similar to the 1D case, the gaps and their distribution statistics of \eqref{fproblemhd} with a potential are quite similar to those without potential, which are reported in Figs. \ref{fig:gapshd}\&\ref{fig:gaps2da}. Those numerical results are omitted here for brevity.

\begin{figure}[h!]
\centerline{
\psfig{figure=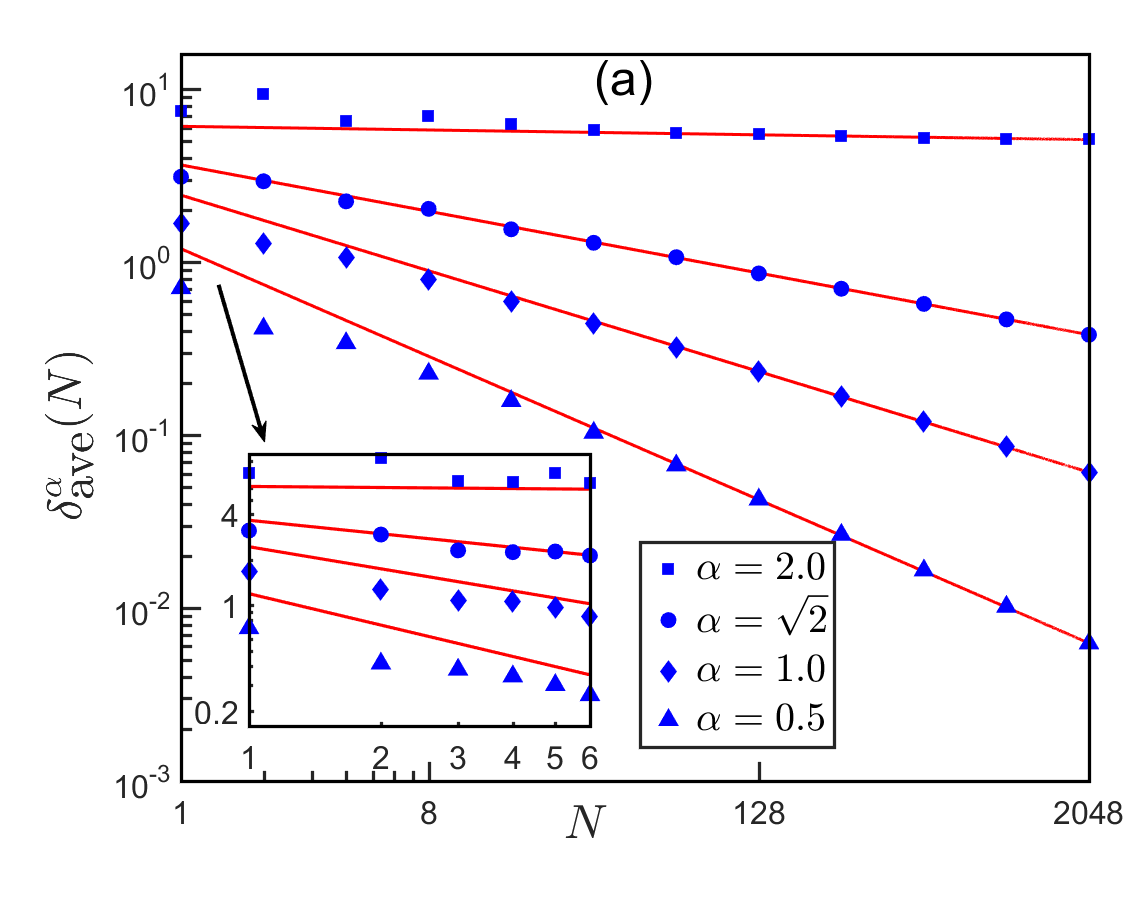,height=5cm,width=7cm,angle=0}
\psfig{figure=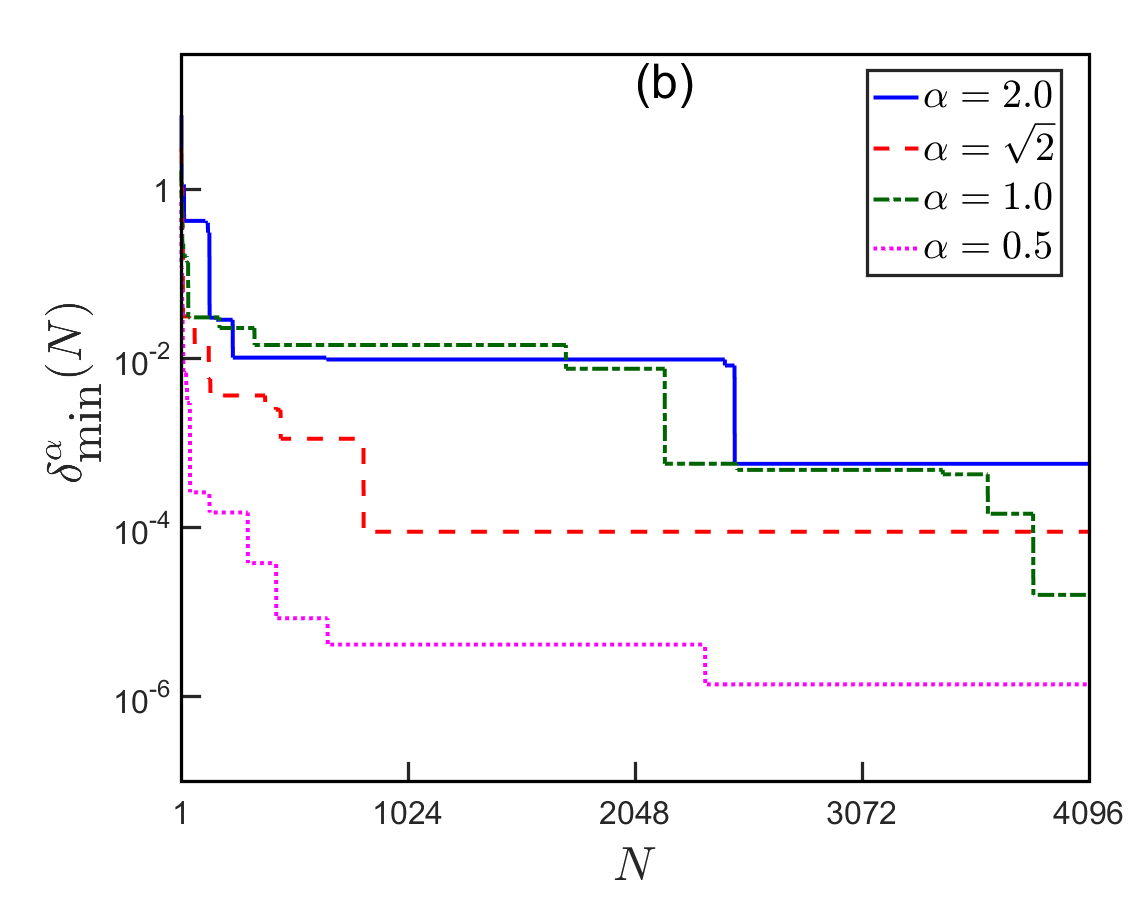,height=5cm,width=7cm,angle=0}}
\caption{Different gaps of \eqref{fproblemhd}
with $d=2$, $L_1=1$, $L_2=\sqrt[3]{2}/2$ and $V(x,y)=\frac{x^2+y^2}{2}$:
(a) the average gaps $\delta_{\textrm{ave}}^\alpha(N)$, and
 (b) the minimum gaps  $\delta_{\textrm{min}}^\alpha(N)$ (symbols denote numerical results and solids lines are from fitting formula when $N\gg1$).}
\label{fig2d:eigdiffv}
\end{figure}

\begin{figure}[h!]
\centerline{
\psfig{figure=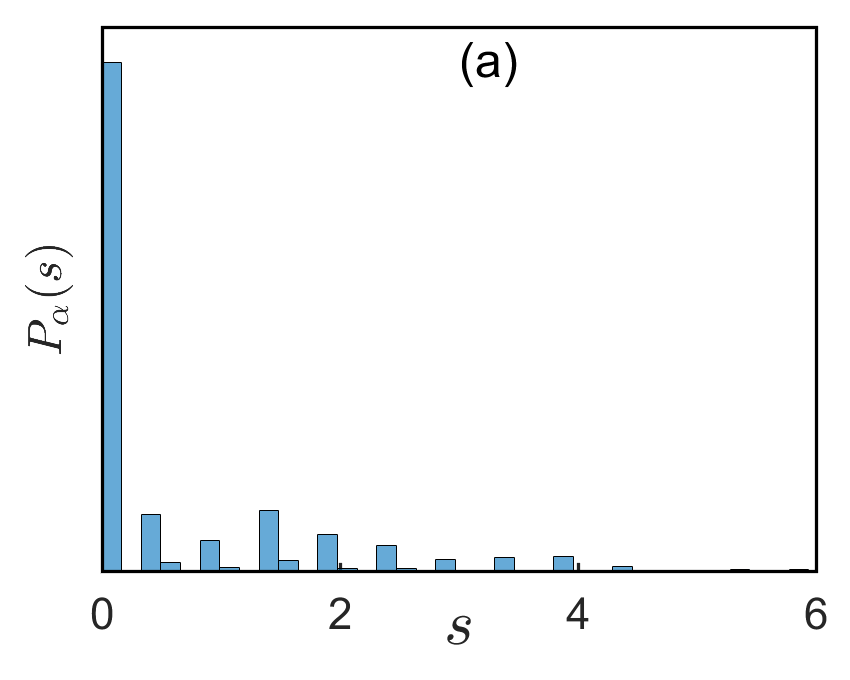,height=5cm,width=7cm,angle=0}
\psfig{figure=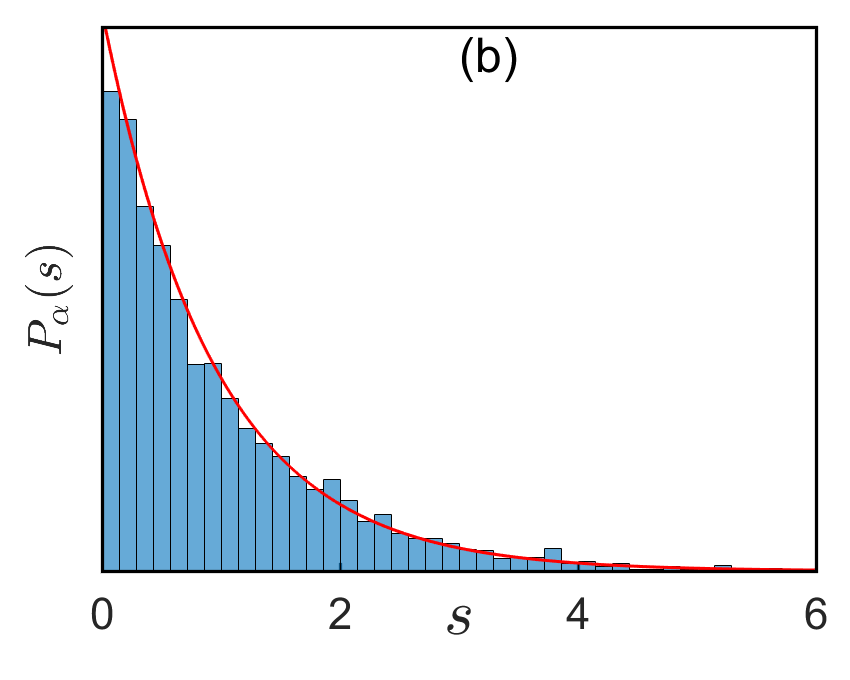,height=5cm,width=7cm,angle=0}}
\caption{The histogram of the normalized gaps  $\{\delta_{\rm norm}^\alpha(n)\ |\ 1\le n\le N=4096\}$ of \eqref{fproblemhd}
with $d=2$ and $V(x,y)=\frac{x^2+y^2}{2}$:
(a) $\alpha=2$ and $L_2=1$; and (b) $\alpha=\sqrt{2}$ and $L_2=\sqrt[3]{2}/2$ (the solid line is a fitting curve by the Poisson distribution).}
\label{fig2d:eigdiffvdv}
\end{figure}

\bigskip

Finally, based on our extensive numerical results and observations, we speculate the following observation (or conjecture) for the D-FSO in
\eqref{fproblemhd} without/with potential:

\smallskip

{\bf Conjecture III} (Gaps and their distribution statistics of D-FSO in
\eqref{fproblemhd} with $d=2$) Assume $0<\alpha\le 2$ and $V(x)\in C(\bar\Omega)$ in \eqref{fproblemhd}, then we have the following asymptotics of its eigenvalues:
\be\label{eigasp2dd}
\lambda_n^\alpha = \frac{4}{2+\alpha}\left(\frac{4\pi}{S}\right)^{\alpha/2} \; n^{\alpha/2}+o(n^{\alpha/2}), \qquad n\gg1,
\ee
where $S$ is the area of $\Omega$. In addition, we have the following asymptotics of different gaps:
\be
\begin{split}
&\delta_{\rm min}^\alpha(N) \to 0, \qquad N\to+\infty,\\
&\delta_{\rm ave}^\alpha(N)
=\frac{4}{2+\alpha}\left(\frac{4\pi}{S}\right)^{\alpha/2} N^{(\alpha-2)/2}+o(N^{(\alpha-2)/2}), \qquad N\gg1.
\end{split}
\ee
In addition, the gap distribution statistics summarized in Tab. \ref{eigsd}
is also valid for \eqref{fproblemhd} in 2D with the potential $V$.

\section{Conclusion}\label{sec:conclusion}
We proposed a Jacobi-Galerkin spectral method  for accurately computing
a large amount of eigenvalues  of the fractional Schr\"{o}dinger operator (FSO). A very important advantage of the proposed numerical method is that,
under a fixed number of degree of freedoms $M$, the Jacobi spectral method
can calculate accurately a large number of eigenvalues with the number
proportional to $M$. Based on the eigenvalues obtained numerically by the proposed method, we obtained several important and interesting results for
the eigenvalues and their different gaps of FSO in 1D and directional FSO in 2D. Based on the gaps, the distribution statistics of the normalized gaps
were obtained numerically for the FSO.

\begin{center}
Acknowledgment
\end{center}

The first author thanks very stimulating discussion with Professor
Zeev Rudnick. This work was partially support by the Ministry of Education of Singapore grant R-146-000-290-114 (W. Bao), the National Natural Science Foundation of China Grants  11671166 (L.Z. Chen), U1930402 (L.Z. Chen \& Y. Ma) and 11771254 (X.Y. Jiang). Part of the work was done when the authors
were visiting the Institute for Mathematical Sciences at the National University of Singapore in 2019 (W. Bao and Y. Ma).

%%%%%%%%%%%%%%%%%%%%%%%%%%%%%%%%%%%

\end{document}